\documentclass[11pt]{amsart}

\usepackage{kotex}
\usepackage{amsmath}
\usepackage{mathrsfs}
\usepackage{amsfonts}
\usepackage{amscd}
\usepackage{slashed}

\usepackage{amssymb}
\usepackage{graphicx}
\usepackage{caption}
\usepackage{subcaption}
\usepackage{dsfont}
\usepackage{color}
\usepackage{hyperref}
\usepackage{bbm}
\usepackage{a4wide}
\usepackage{sseq}
\usepackage{tikz-cd}
\usepackage[abs]{overpic}

\newtheorem{theorem}{Theorem}[section]
\newtheorem{lemma}[theorem]{Lemma}
\newtheorem{proposition}[theorem]{Proposition}
\newtheorem{corollary}[theorem]{Corollary}

\theoremstyle{definition}
\newtheorem*{definition*}{Definition}
\newtheorem{definition}[theorem]{Definition}

\theoremstyle{remark}
\newtheorem*{remark*}{Remark}
\newtheorem{remark}[theorem]{Remark}

\theoremstyle{question}
\newtheorem*{question*}{Question}

\newtheorem{maintheorem}{Setup}

\numberwithin{equation}{section}

\newcommand{\myproof}[2]{Proof of {#1} {#2}}

\begin{document}

\title[Quantization in polarization via transverse PBW theorem]{Quantization in polarization via transverse Poincar\'e--Birkhoff--Witt theorem}

\author[Wang]{Dan Wang}
\address{Department of Mathematics, The Chinese University of Hong
	Kong, Shatin, Hong Kong}
\email{dwang116@link.cuhk.edu.hk}

\author[Yau]{Yutung Yau}
\address{Yau Mathematical Sciences Center, Tsinghua University, Beijing, China}
\email{ytyau@tsinghua.edu.cn}

\thanks{}

\maketitle


\begin{abstract}
	We construct, for every prequantizable symplectic manifold $(M,\omega,L)$ equipped with a nonsingular polarization, a deformation quantization $(\mathcal{C}_M^\infty[[\hbar]],\star)$ together with a filtered sheaf of subalgebras $\mathcal{C}_{M,\hbar}^{(<\infty)}$ that admits an exact action on polarized sections of $L^{\otimes k}$ for each $k\in\mathbb Z^+$. This extends the Kostant--Souriau action of polarization-preserving functions.
	
	The main geometric ingredient is a transverse PBW theory associated with a Nirenberg integrable distribution $P$. For any locally free $P$-modules $E$ and $E'$, we construct a PBW isomorphism for the bundle of transverse differential operators from $E$ to $E'$. This extends the PBW theorem for Lie pairs of Laurent-Gengoux--Sti\'enon--Xu to the setting of $P$-modules.
	
	After specializing $\hbar= \tfrac{\sqrt{-1}}{k}$, we obtain a sheaf $\mathcal O_k^{(<\infty)}$ of level-$k$ quantizable functions and prove that its action identifies $\mathcal O_k^{(<\infty)}$ with the sheaf of transverse differential operators on $L^{\otimes k}$. In the K\"ahler case, this recovers the construction of Chan--Leung--Li.
	
	Finally, for standard symplectic tori with real polarizations, we compare
	our construction with Toeplitz-type operators and obtain their
	semiclassical asymptotic expansions.
\end{abstract}

\section{Introduction}
Geometric quantization selects quantum states through a polarization, but the Kostant--Souriau prescription preserves polarized sections only for polarization-preserving observables \cite{Kos1970,Woo1992}. We ask whether the geometry underlying these states can be used to construct, by a single procedure valid for all nonsingular polarizations, a deformation quantization together with an associative algebra of quantizable observables extending the polarization-preserving ones and acting exactly on polarized sections.\par
We give such a construction by combining Fedosov's method with a transverse Poincar\'e--Birkhoff--Witt isomorphism (Theorem \ref{Third main result}) developed in this paper. Inspired by the sheaf-theoretic construction of quantizable functions introduced by Chan--Leung--Li \cite{ChaLeuLi2023}, our approach works directly with the sheaf of polarized sections (cf. \cite{And1997, Kos1975}) and does not presuppose a
global operator-and-symbol calculus adapted to the polarization.\par
The basic setup considered in this paper is as follows:

\begin{maintheorem}
	\label{Setup: one}
	$(M, \omega)$ is a symplectic manifold, equipped with a prequantum line bundle $(L, \nabla^L)$ and a Nirenberg integrable complex Lagrangian subbundle $P \subset TM_\mathbb{C}$.
\end{maintheorem}
Henceforth, by a \emph{polarization}, we shall mean such a subbundle $P$ as in Setup \ref{Setup: one}. Let $\mathcal{C}_M^\infty$ be the sheaf of smooth complex valued functions on $M$, let $\mathcal{O}$ be the sheaf of smooth functions constant along $P$ and let $\mathcal{L}^{\otimes k}$ be the sheaf of $P$-polarized sections of $L^{\otimes k}$. For a smooth function $f$ on $M$, let $X_f$ be its Hamiltonian vector field. We say $f$ \emph{preserves} $P$ if $[X_f, \,\cdot\,]$ preserves smooth sections of $P$. Our first main result is:

\begin{theorem}[Theorems \ref{Theorem: star product}, \ref{Theorem: zeroth order quantizable functions} and \ref{Theorem: first order quantizable functions}]
	\label{First main result}
	Under Setup \ref{Setup: one}, there exist a deformation quantization $(\mathcal{C}_M^\infty[[\hbar]], \star)$, a sheaf of filtered subalgebras $\mathcal{C}_{M, \hbar}^{(<\infty)}$ with increasing filtration
	\begin{equation*}
		\mathcal{C}_{M,\hbar}^{(0)}
		\subset
		\mathcal{C}_{M,\hbar}^{(1)}
		\subset \cdots \subset
		\mathcal{C}_{M,\hbar}^{(r)}
		\subset \cdots \subset
		\mathcal{C}_{M,\hbar}^{(<\infty)},
	\end{equation*}
	and for each $k \in \mathbb{Z}^+$, an action $\rho_k$ of $\mathcal{C}_{M, \hbar}^{(<\infty)}$ on $\mathcal{L}^{\otimes k}$, with the following properties.
	\begin{itemize}
		\item The zeroth filtered piece $\mathcal{C}_{M,\hbar}^{(0)}$ is $\mathcal{O}$. On this piece, the star product $\star$ and the actions $\rho_k$ reduce to ordinary multiplication.
		\item The first filtered piece $\mathcal{C}_{M,\hbar}^{(1)}$ consists precisely of elements of the form $f = f_0 + \hbar f_1$, where $f_0 \in \mathcal{C}_M^\infty$ preserves $P$ and $f_1 \in \mathcal{O}$. The action $\rho_k$ of such an $f$ is given by the Kostant--Souriau operator associated with the function $f_{[k]} := f_0 + \tfrac{\sqrt{-1}}{k} f_1$:
		\[
		\rho_k(f)s
		= \frac{\sqrt{-1}}{k}
		\nabla^{L^{\otimes k}}_{X_{f_{[k]}}}s + f_{[k]} \cdot s.
		\]
	\end{itemize}
\end{theorem}

Since every $f\in\mathcal O$ preserves $P$, the inclusion $\mathcal C_{M,\hbar}^{(0)} \subset \mathcal C_{M,\hbar}^{(1)}$ is immediate from the above characterization. We call sections of $\mathcal{C}_{M, \hbar}^{(<\infty)}$ \emph{formal quantizable functions}; for $M = T^*B$ with the vertical polarization, these are precisely formal functions polynomial in $\hbar$ and fibre variables.\par

Our result gives an action of a distinguished sheaf of subalgebras, not of the full deformation algebra. This  action is \emph{exact} at each level $k$, as opposed to merely asymptotic, and one should not expect it to extend to the entire star-product algebra.\par

For each $k \in \mathbb{Z}^+$, the same construction also produces a sheaf $\mathcal{O}_k^{(<\infty)}$ of \emph{level-$k$ quantizable functions}, which is a non-formal counterpart of $\mathcal{C}_{M,\hbar}^{(<\infty)}$. Formal quantizable functions admit a natural evaluation map $\mathcal{C}_{M,\hbar}^{(<\infty)} \to \mathcal{O}_k^{(<\infty)}$ at $\hbar=\tfrac{\sqrt{-1}}{k}$, and the action $\rho_k$ is compatible with the resulting level-$k$ action on polarized sections. Our second main result gives:

\begin{theorem}[Theorem \ref{Theorem: isomorphism of transverse differential operators}]
	\label{Second main result}
	Under Setup \ref{Setup: one}, for each $k \in \mathbb{Z}^+$, the induced action of $\mathcal{O}_k^{(<\infty)}$ on $P$-polarized sections of $L^{\otimes k}$ defines an isomorphism of sheaves of filtered algebras
	\begin{equation*}
		\mathcal{O}_k^{(<\infty)} \to \widetilde{\mathcal{D}}_{L^{\otimes k}},
	\end{equation*}
	where $\widetilde{\mathcal D}_{L^{\otimes k}}$ is the sheaf of transverse differential operators on the $P$-module $L^{\otimes k}$.
\end{theorem}

Thus the level-$k$ algebra $\mathcal{O}_k^{(<\infty)}$ is not merely an auxiliary finite-level object: it is precisely the full sheaf of transverse differential operators, with the filtration corresponding to differential-operator order. For a K\"ahler polarization, this recovers the quantizable-function and differential-operator picture of Chan--Leung--Li \cite{ChaLeuLi2023}, with the corresponding ordering conventions. For the vertical polarization of a cotangent bundle, it recovers the correspondence between polynomial symbols and differential operators. Our construction also extends to $L^{\otimes k} \otimes \mathbf{L}$, where $\mathbf{L}$ is an auxiliary complex line bundle with a connection flat along $P$, including half-form bundles when available; see Section \ref{Section: construction of star products via action on polarized sections}.\par
A key feature of our construction is that the same algebra--module formalism applies uniformly to real, K\"ahler, and mixed polarizations. It requires neither a global leaf-space description nor a Hilbert-space realization of the quantum states. Moreover, although a choice of a Lagrangian complement to $P$ enters our construction, this complement is not required to be involutive; in particular, it need not define a transverse polarization.\par

The key ingredient in the proofs of Theorems \ref{First main result} and \ref{Second main result} is our construction of transverse PBW isomorphisms, which constitutes the third main result of the paper and is developed in the next subsection.

\subsection{Our construction of transverse PBW isomorphisms}
\quad\par
\label{Subsection: general theory of transverse differential operators}
We develop a transverse differential calculus associated with a Nirenberg integrable distribution $P\subset TM_{\mathbb C}$, independently of any symplectic or prequantization data. These additional structures
enter only when we apply the theory to polarized deformation quantization.\par

A locally free $P$-module is a complex vector bundle $E$ with a flat $P$-connection $d^E_P$ whose $d^E_P$-closed-section sheaf is locally free over $\mathcal O$. Its transverse jet bundle $\widetilde J E\subset JE$ consists, at finite order, of jets represented by germs of $d^E_P$-closed sections. The Grothendieck connection preserves this subbundle.\par
For two such modules $E,E'$, let $\widetilde D(E,E')$ be the bundle of transverse differential operators, with its order filtration. Its degree-$r$ symbol bundle is $\operatorname{Sym}^r Q\otimes E^*\otimes E'$, where $Q:=TM_{\mathbb C}/P$. We construct a \emph{Poincar\'e--Birkhoff--Witt (PBW) map} $\operatorname{pbw}^{E, E'}$ after fixing the auxiliary data
\begin{equation}
	\tag{$\sharp$}\label{Equation: auxiliary data for pbw isomorphism}
	\begin{cases}
		\text{a splitting of the short exact sequence of vector bundles } P \to TM_\mathbb{C} \to Q;\\
		\text{a torsion-free connection on } Q;\\
		\text{a connection on } E \text{ extending the flat } P \text{-connection } d_P^E.
	\end{cases}
\end{equation}
Here torsion-freeness is understood in the quotient-bundle sense of Remark \ref{Remark: torsion free connection on quotient bundles}. Such data always exist. This construction yields our third main result, a PBW-type theorem:

\begin{theorem}[Theorem \ref{Theorem: pbw theorem}]
	\label{Third main result}
	Let $M$ be a smooth manifold equipped with a Nirenberg integrable complex distribution $P$. Then for any locally free $P$-modules $E, E'$, the map
	\begin{equation*}
		\operatorname{pbw}^{E, E'}: \Gamma(-, \operatorname{Sym} Q \otimes E^* \otimes E') \to \Gamma(-, \widetilde{D}(E, E'))
	\end{equation*}
	is an isomorphism of filtered $\mathcal C_M^\infty$-modules whose associated graded map is the identity map on $\operatorname{Sym}Q\otimes E^*\otimes E'$.
\end{theorem}

When $E=E'=\underline{\mathbb C}$, this is the PBW map of
Laurent-Gengoux--Sti\'enon--Xu \cite{LauStiXu2021} for the Lie pair $(TM_{\mathbb C},P)$. Related PBW results for inclusions of Lie algebroids were established by Calaque \cite[Theorem 1.2]{Cal2014}. The present construction differs in that it provides a smooth filtered splitting for arbitrary locally free $P$-modules, which is precisely the form needed for polarization-adapted deformation quantization.\par
The significance of Theorem \ref{Third main result} for our construction lies in the flat connection $\nabla^{\operatorname{K},E}$ on $\widehat{\operatorname{Sym}}Q^* \otimes E$, obtained by transporting the Grothendieck connection on $\widetilde{J}E$ through the PBW-induced identification
\begin{equation*}
	\Gamma(M, \widehat{\operatorname{Sym}} Q^* \otimes E) \cong \Gamma(M, \widetilde{J}E).
\end{equation*}
We refer to $\nabla^{\operatorname{K}, E}$ as the \emph{Kapranov connection}, in recognition of Kapranov's original investigation of the $L_\infty$-structures defined by this connection in the K\"ahler setting \cite{Kap1999}. Its flat-section sheaf is naturally identified with the
$d^E_P$-closed sections of $E$ (Proposition \ref{Proposition: qausi-isomorphism for polarized sections}). In the next subsection, we explain how
$\nabla^{\operatorname{K},E}$ enters the construction of our deformation quantization.

\subsection{Geometric mechanism behind Theorems \ref{First main result} and \ref{Second main result}}
\quad\par
Returning to Setup \ref{Setup: one}, we adapt Fedosov's method to the polarization $P$. To this end, choose
\begin{equation}
	\tag{$\dagger$}\label{Equation: auxiliary data for deformation quantization}
	\begin{cases}
		\text{a Lagrangian complement } Q \subset TM_\mathbb{C} \text{ of } P;\\
		\text{a torsion-free connection on } TM_{\mathbb C}/P \cong Q.
	\end{cases}
\end{equation}
These choices determine a connection preserving $\omega$, $P$, and $Q$, which may have torsion. A torsion-free connection preserving both distributions would force $Q$ to be involutive; allowing torsion avoids imposing a second polarization transversal to $P$.\par
The splitting $TM_{\mathbb C}=P\oplus Q$ determines a
polarization-adapted fibrewise star product $\star^{\operatorname{F}}$ on the Weyl bundle $\mathcal{W} = \widehat{\operatorname{Sym}} T^*M_\mathbb{C}[[\hbar]]$. Its finite-polarized-weight part
\begin{equation*}
	\mathcal{W}^{<\infty} := \widehat{\operatorname{Sym}} Q^* \otimes \operatorname{Sym} P^*[\hbar]
\end{equation*}
admits a fibrewise action $\circledast_k^{\operatorname{F}}$ on $\widehat{\operatorname{Sym}}Q^*\otimes L^{\otimes k}$ after setting $\hbar=\tfrac{\sqrt{-1}}{k}$. 
Intuitively, under this action the $Q^*$-variables act by
multiplication, while the $P^*$-variables act, via the symplectic pairing between $P$ and $Q$, as derivatives with respect to the $Q^*$-variables. Here the polarized weight counts the $P^*$-degree and the $\hbar$-degree, but not the $Q^*$-degree.\par
Our key observation is that the Kapranov connections determine a single $\mathcal{W}^{< \infty}$-valued one-form $\gamma$ such that, for every $k$,
\begin{equation*}
	\nabla^{\operatorname{K}, L^{(k)}} = \nabla^{\widehat{\operatorname{Sym}} Q^* \otimes L^{(k)}} + \tfrac{k}{\sqrt{-1}} \gamma \circledast_k^{\operatorname{F}}.
\end{equation*}
Their flatness implies that $\gamma$ satisfies Fedosov's equation. Consequently, $D = \nabla + \tfrac{1}{\hbar} \left[ \gamma, \, \cdot\, \right]_{\star^{\operatorname{F}}}$ is a flat algebra connection compatible with the module connections. Rather than constructing the deformation quantization and subsequently seeking a representation, we obtain both from the transverse geometry of the polarization.\par
The $D$-flat lift $\mathbf{J}_f$ of a formal function $f$ defines its star product, and the condition $\mathbf{J}_f \in \mathcal{W}^{< \infty}$ defines formal quantizability. After specialization at $\hbar=\tfrac{\sqrt{-1}}{k}$, this condition ensures that the corresponding action has finite order in the transverse directions, giving a transverse differential operator; at first order, this
recovers the Kostant--Souriau prescription. The isomorphism in Theorem \ref{Second main result} is obtained from a filtered quasi-isomorphism between the finite-polarized-weight Weyl complex and the complex of transverse differential operators. These constructions constitute the geometric mechanism behind Theorems \ref{First main result} and \ref{Second main result}.

\subsection{Relation to previous work}
\quad\par
The extension of the Kostant--Souriau prescription to higher-order observables has an early predecessor in Hess \cite{Hes1981}. His construction is based on a pair of transverse involutive Lagrangian distributions. In general, however, the resulting higher-order operators need not preserve polarized sections. Nor does his construction produce
a filtered subalgebra of a deformation quantization admitting an exact action of the type considered here.\par
Chan--Leung--Li \cite{ChaLeuLi2022b,ChaLeuLi2023} provide the direct model for our approach. In the K\"ahler setting, they constructed sheaves of quantizable functions with non-formal specializations and identified the resulting level-$k$ algebras with sheaves of twisted differential operators. In \cite{Yau2025}, the relevant flat connection
on $\widehat{\operatorname{Sym}} T^{*(1, 0)}M \otimes L^{\otimes k}$ was identified precisely with the Kapranov connection. The present work develops a transverse PBW theory that extends this Kapranov-connection mechanism from the K\"ahler setting to arbitrary nonsingular polarizations.\par
Beyond this direct line of development, several complementary approaches relate deformation quantization to polarizations and to representations on quantum states. We briefly compare the ones most relevant to our construction.

\subsubsection{Polarized deformation quantization and modules}
\quad\par
Reshetikhin--Yakimov \cite{ResYak2000} constructed and classified star products adapted to Lagrangian fibrations, preserving the commutative algebra of functions pulled back from the base. Bressler--Donin \cite{BreDon2000} and Donin \cite{Don2003} developed a general theory of polarized deformation quantization, in which the polarization is encoded by a distinguished commutative subsheaf of the deformation quantization, with results on existence, characteristic classes, and classification.\par
Our emphasis is different: we require compatibility with the sheaf of polarized quantum states. The transverse PBW formalism incorporates the module structure on these states into the construction of the star product itself. Related viewpoints include Pol\'shin's formal modules on leaves of real polarizations \cite{Pol2008} and Tsygan's oscillatory modules \cite{Tsy2009}.

\subsubsection{Fedosov constructions and cotangent-bundle representations}
\quad\par
In the K\"ahler setting, Bordemann--Waldmann \cite{BorWal1997} constructed Wick-type star products by a modification of Fedosov's method. On cotangent bundles, Bordemann--Neumaier--Waldmann \cite{BorNeuWal1998} constructed homogeneous Fedosov star products together with differential-operator representations. Bordemann--Neumaier--Pflaum--Waldmann \cite{BorNeuPflWal2003} extended this framework to Hermitian line
bundles and magnetic terms, relating formal representations to global
pseudodifferential symbol calculus.\par
These constructions provide important models in which a deformation
quantization is accompanied by a compatible operator representation. Our
cotangent-bundle example recovers the polynomial-symbol correspondence,
while the transverse construction itself does not rely on a
cotangent-bundle structure.

\subsubsection{Global operator and symbol calculi}
\quad\par
Another important approach to deformation quantization starts from a global operator-and-symbol calculus. In the K\"ahler setting, Berezin's symbolic calculus and its geometric interpretation by Moreno--Ortega-Navarro \cite{MorOrt1983} and Cahen--Gutt--Rawnsley \cite{CahGutRaw1990, CahGutRaw1993} provide a route from operator composition to deformation quantization. The Berezin--Toeplitz construction gives another prominent example, where the asymptotic composition of Toeplitz operators determines a star product on compact quantizable K\"ahler manifolds \cite{BorMeiSch1994, Gui1995, Sch2000}.\par
These constructions rely on global analytic structures associated with
the K\"ahler polarization, whereas our starting point is the local
differential geometry of polarized sections. This distinction is
particularly relevant for mixed polarizations, for which a comparable global operator-and-symbol calculus is generally unavailable. Mixed polarizations and their associated quantum states have been studied from various geometric-quantization viewpoints; see \cite{BaiFerHilMouNun2025, BaiHilKayMouNun2025, LeuWan2023, LeuWan2024a, LeuWan2024b, MouNunPer2020, MouNunPerWan2024, Wan2024}.\par
When an independent operator calculus is available, comparing it with our
transverse construction becomes a separate analytic question. We
address this question in the next subsection for Toeplitz-type operators
in real polarization.

\subsection{Comparison with Toeplitz-type operators in real polarization}
\quad\par
For any $f \in \mathcal{C}^\infty(M)$, the flat-lift construction gives a series of differential operators $\sum_{r=0}^\infty T_{f, r, k}$ (Definition \ref{Definition: components of actions}). Individual terms of this series need not preserve polarized sections, and no convergence is asserted in general. If, furthermore, $f$ is formal quantizable, then the series terminates and its sum gives the exact action of Theorem \ref{First main result}. In the K\"ahler case, $\sum_{r=0}^\infty T_{f, r, k}$ was related to Toeplitz operators by Chan--Leung--Li and the second author \cite{ChaLeuLiYau2025}; see also Andersen \cite{And2024}.\par
Toeplitz-type operators $T_{f, k}$ acting on the space $\Gamma_P^{-\infty}(M, L^{\otimes k})$ of distributional quantum states in real polarization were constructed by Leung and the second author \cite{LeuYau2022,LeuYau2023}. 
We compare our operators with those of \cite{LeuYau2023} on the standard symplectic torus $M=\mathbb R^{2n}/\mathbb Z^{2n}$, equipped with its canonical prequantum line bundle $L$ and the standard pair $(P,Q)$ of transverse real polarizations. The comparison is expressed by a Taylor-type expansion with a uniform remainder estimate.

\begin{theorem}[Theorem \ref{Theorem: symplectic tori}]
	\label{Fourth main result}
	In this torus setting, let $f \in \mathcal{C}^\infty(M)$ and $N \in \mathbb{N}$. Then there exists a constant $C_{f, N} > 0$ such that for all $k \in \mathbb{Z}^+$, $s \in \Gamma_P^{-\infty}(M, L^{\otimes k})$ and $\tau \in \Gamma(M, (L^*)^{\otimes k})$,
	\begin{equation*}
		\left\lvert \left\langle T_{f, k} s - \sum_{r=0}^N T_{f, r, k} s, \tau \right\rangle \right\rvert \leq C_{f, N} \cdot \lVert s \rVert_{\ell^1} \cdot \lVert \tau \rVert_{N+1} \cdot \frac{1}{k^{N+1}},
	\end{equation*}
	where the norms $\lVert \, \cdot\, \rVert_{\ell^1}$ and $\lVert \, \cdot\, \rVert_{N+1}$ are defined in \eqref{Equation: l1 norm} and \eqref{Equation: seminorm}, respectively.
\end{theorem}

The estimate compares the operators through their distributional pairings. When $f$ is formal quantizable in this model, the finite sum agrees exactly with the Toeplitz-type operator. The torus example therefore provides a concrete bridge between the general transverse construction developed in this paper and an independently defined Toeplitz-type semiclassical calculus in real polarization.

\subsection{Structure of the paper}
\quad\par
This paper is organized into two parts.\par
The first part, consisting of Sections \ref{Section: transverse jet bundles and transverse differential operators}
and \ref{Section: transverse pbw theorem}, develops the theory of transverse jet bundles and transverse differential
operators. Related notions have appeared in \cite{CheXiaXu2022, Faz2025a, ZubShu2016}, but not in a form suited to the polarization-adapted deformation-quantization problem considered here. We fill this gap by establishing the transverse PBW isomorphisms of Theorem \ref{Third main result} and by carrying out a systematic study of the induced Kapranov connections. These constructions form the geometric foundation for the quantization theory developed in the second part and may also be of independent interest.

The second part develops our construction of polarization-adapted deformation quantization. Section \ref{Section: Fedosov quantization in the presence of a polarization} introduces the polarization-adapted fibrewise algebra and module structures. Section \ref{Section: construction of star products via action on polarized sections} constructs the compatible Fedosov connection, proves Theorems \ref{First main result} and \ref{Second main result}, and discusses several examples. Finally, Section \ref{Section: asymptotics of Toeplitz operators in real polarization} compares the resulting operators with Toeplitz-type operators in real
polarization and proves Theorem \ref{Fourth main result}.

\subsection{Notations}
\quad\par
\label{Subsection: notations}
Throughout this paper, $M$ is a smooth manifold equipped with a Nirenberg integrable complex distribution $P$ on it, i.e. a complex vector subbundle of $TM_\mathbb{C}$ such that $P \cap \overline{P}$ is of constant rank and both $P$ and $P + \overline{P}$ are involutive. We write $Q = TM_\mathbb{C}/  P$, identifying it with a complementary subbundle when a splitting has been chosen. The symbol $\underline{\mathbb{C}}$ denotes the trivial flat complex line bundle $(M \times \mathbb{C}, d)$ over $M$.\par
We write $\mathcal{C}^\infty(M)$ for the algebra of smooth $\mathbb{C}$-valued functions on $M$, $\mathcal{C}_M^\infty$ for the sheaf of smooth $\mathbb{C}$-valued functions on $M$, and
\begin{equation*}
	\mathcal{O}(U) := \{ f \in \mathcal{C}^\infty(U): \mathcal{L}_Y f = 0 \text{ for all } Y \in \Gamma(U, P) \}
\end{equation*}
for any open subset $U \subset M$. For a (pro-)vector bundle $V$ over $M$, set
\begin{equation*}
	C^{q, p}(M, V) := \Gamma(M, \textstyle \bigwedge^q Q^* \otimes \bigwedge^p P^* \otimes V).
\end{equation*}
When $V$ is a complex vector bundle of finite rank, the notations $\operatorname{Sym} V$ and $\widehat{\operatorname{Sym}} V$ denote the symmetric algebra bundle and its completion. Our convention for symmetric tensors is
\begin{equation}
	Z_1 \cdots Z_r = \frac{1}{r!} \sum_{\sigma \in \mathcal{S}_r} Z_{\sigma(1)} \otimes \cdots \otimes Z_{\sigma(r)},
\end{equation}
where $\mathcal{S}_r$ is the symmetric group on $r$ elements. An empty product is $1$.\par
For local frames $(v_i)$ of $Q$ and $(\check v_i)$ of $P$, we write $(v^i)$ and $(\check v^i)$ for their dual frames. When regarded as symmetric rather than exterior variables, $v^i$ and $\check v^i$ are denoted by $u^i$ and $\check u^i$.\par
For arbitrary complex vector bundles $E, E', E''$ over $M$, the bundles $J^rE$ and $JE := \varprojlim_r J^rE$ denote finite and infinite jets, with Grothendieck connection $\nabla^{G,E}$ on $JE$. We set
\[
D^r(E, E') := \operatorname{Hom}(J^rE, E'), \quad D(E, E') := \varinjlim_r D^r(E, E').
\]
A smooth section of $D^r(E,E')$ is identified with the corresponding differential operator of order $r$; composition is denoted by $\circ$. If $E$ has a connection $\nabla^E$, its curvature is written $R^E$.

\subsection*{Acknowledgement}
\quad\par
This work was supported by World Premier International Research Center Initiative (WPI), MEXT, Japan. The first author is grateful to Christian Blohmann, Eckhard Meinrenken, Ioan Mărcuț, Marius Crainic, Jean-Marie Lescure, and Qingyuan Jiang for helpful discussions, as well as to the Max Planck Institute for Mathematics (MPIM) in Bonn for their hospitality and financial support. The second author thanks Hsuan-Yi Liao, Mathieu Sti\'enon and Ping Xu for patiently answering his questions on the theory of Lie pairs during his visit to National Tsing Hua University, and gratefully acknowledges the support of the National Center for Theoretical Sciences, where part of this work was carried out.

\section{Transverse jet bundles and transverse differential operators}
\label{Section: transverse jet bundles and transverse differential operators}
A $P$-\emph{module} is a pair $(E, d_P^E)$ for which $E$ is a complex vector bundle over $M$ and
\begin{equation*}
	d_P^E: \Gamma(M, E) \to \Gamma(M, P^* \otimes E)
\end{equation*}
is a flat $P$-connection on $E$. Throughout this paper, we focus on $P$-modules $(E, d_P^E)$ that are \emph{locally free}, meaning that the sheaf of $d_P^E$-closed sections of $E$ forms a locally free $\mathcal{O}$-module.\par
In Subsection \ref{Subsection: transverse jet bundles} we introduce the \emph{transverse jet bundle} $\widetilde{J}E$ of $(E, d_P^E)$, a distinguished subbundle of the jet bundle $JE$ associated to the flat $P$-connection. Subsection \ref{Subsection: Grothendieck connection on transverse jet bundles} shows that the Grothendieck connection on $JE$ preserves $\widetilde{J}E$. Subsection \ref{Subsection: transverse differential operators} then defines the \emph{bundle of transverse differential operators} between two locally free $P$-modules. Finally, Subsection \ref{Subsection: Canonical module structures on transverse structures} equips these bundles with canonical $P$-module structures and introduces the notion of \emph{transverse differential operators}.\par
As a remark, when $P = T^{0, 1}M$ for a complex manifold $M$, every $P$-module is locally free. For a general Nirenberg integrable distribution $P$, one cannot generally expect this to hold. However, the following theorem --- essentially Theorem 3 from \cite{Raw1977} --- remains valid.

\begin{theorem}
	\label{Theorem: rank-1 line bundles are foliated}
	For Nirenberg integrable $P \subset TM_\mathbb{C}$, every rank-$1$ $P$-module is locally free.
\end{theorem}

\subsection{Transverse jet bundles}
\quad\par
\label{Subsection: transverse jet bundles}
Consider a locally free $P$-module $(E, d_P^E)$. There induces a $\mathcal{C}^\infty(M)$-linear map
\begin{equation*}
	\Gamma(M, P) \to \Gamma(M, D^1(E, E)), \quad Y \mapsto \iota_Y \circ d_P^E,
\end{equation*}
where $D^1(E, E)$ is the bundle of first order differential operators from $E$ to itself (see Subsection \ref{Subsection: notations} for relevant notations). For $r \in \mathbb{N}$, denote by $\langle \, \cdot\,, \, \cdot\, \rangle$ the natural pairing
\begin{equation*}
	\Gamma(M, D^r(E, \underline{\mathbb{C}})) \times \Gamma(M, J^rE) \to \mathcal{C}^\infty(M).
\end{equation*}
\begin{definition}
	For each $r \in \mathbb{N}$, we define a $\mathcal{C}_M^\infty$-submodule $\Gamma(-, \widetilde{J}^rE)$ of $\Gamma(-, J^rE)$ as follows: $\Gamma(-, \widetilde{J}^0E)$ is the sheaf of smooth sections of $E$; for $r > 0$ and an open subset $U$ of $M$, $\Gamma(U, \widetilde{J}^rE)$ is the module of sections $\sigma \in \Gamma(U, J^rE)$ satisfying the condition that
	\begin{equation*}
		\left\langle \Phi \circ \iota_Y \circ d_P^E, \sigma \right\rangle = 0 \quad \text{for all } Y \in \Gamma(U, P) \text{ and } \Phi \in \Gamma(U, D^{r-1}(E, \underline{\mathbb{C}})).
	\end{equation*}
\end{definition}

\begin{remark}
	It immediately follows from the above definition that, for every $d_P^E$-closed section $s \in \Gamma(U, E)$, its $r$-jet $\mathfrak{j}^r s$ lies in $\Gamma(U, \widetilde{J}^rE)$.
\end{remark}

The following proposition states that $\Gamma(-, \widetilde{J}^rE)$ defines a smooth subbundle $\widetilde{J}^rE$ of $J^rE$, which we call the \emph{transverse} $r$\emph{-jet bundle} of $E$, verifying that $\Gamma(-, \widetilde{J}^rE)$ is a suitable notation (for the ease of notations, we let $J^{-1}E = \widetilde{J}^{-1}E$ be the zero vector bundle over $M$).

\begin{proposition}
	\label{Proposition: transverse jet bundles are smooth subbundles}
	For all $r \in \mathbb{N}$, the sheaf $\Gamma(-, \widetilde{J}^rE)$ is a locally free $\mathcal{C}_M^\infty$-module of finite rank, and the canonical short exact sequence
	\begin{center}
		\begin{tikzcd}
			0 \ar[r] & \operatorname{Sym}^r T^*M_\mathbb{C} \otimes E \ar[r] & J^rE \ar[r] & J^{r-1}E \ar[r] & 0
		\end{tikzcd}
	\end{center}
	restricts to a short exact sequence of complex vector bundles:
	\begin{equation}
		\label{Equation: short exact sequence of transverse jet bundles}
		\begin{tikzcd}
			0 \ar[r] & \operatorname{Sym}^r Q^* \otimes E \ar[r] & \widetilde{J}^rE \ar[r] & \widetilde{J}^{r-1}E \ar[r] & 0
		\end{tikzcd}
	\end{equation}
	Here, $\widetilde{J}^rE \subset J^rE$ denotes the subbundle whose smooth sections form the sheaf $\Gamma(-, \widetilde{J}^rE)$.
\end{proposition}

This implies that we have a sequence of surjective complex vector bundle morphisms
\begin{center}
	\begin{tikzcd}
		\cdots \ar[r] & \widetilde{J}^{r+1}E \ar[r] & \widetilde{J}^rE \ar[r] & \cdots \ar[r] & \widetilde{J}^1E \ar[r] & \widetilde{J}^0E = E
	\end{tikzcd}
\end{center}

\begin{definition}
	The \emph{transverse jet bundle} of $(M, P)$, denoted by $\widetilde{J}E$, is defined as the projective limit of the above sequence.
\end{definition}

The transverse jet bundle $\widetilde{J}E$ is then a pro-vector bundle and a subbundle of $JE$. From Proposition \ref{Proposition: transverse jet bundles are smooth subbundles} we see that $\widetilde{J}E$ has a decreasing filtration whose $r$th filtered piece is given by $\ker \left( \widetilde{J}E \to \widetilde{J}^rE \right)$, and the associated graded pro-vector bundle of $\widetilde{J}E$ is
\begin{equation*}
	\widehat{\operatorname{Sym}} Q^* \otimes E = \varprojlim_r \operatorname{Sym}^{\leq r} Q^* \otimes E.
\end{equation*}

As the proof of Proposition \ref{Proposition: transverse jet bundles are smooth subbundles} involves tedious technical details which are unimportant for our main theme of this section, we will postpone it to Appendix \ref{Section: Proof of locally freeness of transverse jet bundles}.

\subsection{The Grothendieck connection}
\quad\par
\label{Subsection: Grothendieck connection on transverse jet bundles}
Let $(E, d_P^E)$ be a locally free $P$-module. The jet bundle $JE$ has a canonical flat connection $\nabla^{\operatorname{G}, E}$ on $JE$, known as the \emph{Grothendieck connection}, satisfying the condition that a smooth section $s \in \Gamma(M, JE)$ is $\nabla^{\operatorname{G}, E}$-flat if and only if it is the infinite jet of a smooth section of $E$.

\begin{proposition}
	\label{Proposition: Grothendieck connection preserves transverse jet bundles}
	The Grothendieck connection $\nabla^{\operatorname{G}, E}$ preserves the subbundle $\widetilde{J}E$.
\end{proposition}

We recall the definition of $\nabla^{\operatorname{G}, E}$, adopting a formulation that is suited to the proof of the preceding proposition. Observe that any local smooth section $\sigma \in \Gamma(U, JE)$ naturally induces a sequence of $\mathcal{C}^\infty(U)$-linear maps
\begin{equation*}
	\Gamma(U, D^r(E, \underline{\mathbb{C}})) \to \mathcal{C}^\infty(U) \quad (\text{for } r \in \mathbb{N}),
\end{equation*}
all denoted by $\langle \, \cdot\,, \sigma \rangle$ by abuse of notation, such that the diagram below is commutative:
\begin{center}
	\begin{tikzcd}
		\Gamma(U, D^r(E, \underline{\mathbb{C}})) \ar[rr] \ar[dr, "\langle \, \cdot\,{,} \sigma \rangle"'] && \Gamma(U, D^{r+1}(E, \underline{\mathbb{C}})) \ar[dl, "\langle \, \cdot\,{,} \sigma \rangle"]\\
		& \mathcal{C}^\infty(U)
	\end{tikzcd}
\end{center}
Here, the horizontal arrow is the inclusion map. Indeed, this gives a one-to-one correspondence between sections in $\Gamma(U, JE)$ and such sequences of $\mathcal{C}^\infty(U)$-linear maps. Then $\nabla^{\operatorname{G}, E}$ is defined as follows: for all $X \in \Gamma(U, TM_\mathbb{C})$, $\sigma \in \Gamma(U, JE)$, $r \in \mathbb{N}$ and $\Phi \in \Gamma(U, D^r(E, \underline{\mathbb{C}}))$,
\begin{equation}
	\label{Equation: formula of Grothendieck connection}
	\langle \Phi, \nabla_X^{\operatorname{G}, E} \sigma \rangle = \mathcal{L}_X \langle \Phi, \sigma \rangle - \langle \mathcal{L}_X \circ \Phi, \sigma \rangle.
\end{equation}

\begin{proof}[\myproof{Proposition}{\ref{Proposition: Grothendieck connection preserves transverse jet bundles}}]
	Consider any open $U \subset M$, $X \in \Gamma(U, TM_\mathbb{C})$ and $\sigma \in \Gamma(U, \widetilde{J}E)$. Fix $Y \in \Gamma(U, P)$, $r \in \mathbb{N}$ and $\Phi \in \Gamma(U, D^r(E, \underline{\mathbb{C}}))$. Then
	\begin{align*}
		\langle \Phi \circ \iota_Y \circ d_P^E, \nabla_X^{\operatorname{G}, E} \sigma \rangle = \mathcal{L}_X \langle \Phi \circ \iota_Y \circ d_P^E, \sigma \rangle - \langle \mathcal{L}_X \circ \Phi \circ \iota_Y \circ d_P^E, \sigma \rangle = 0.
	\end{align*}
	This concludes the proof.
\end{proof}

\subsection{Bundles of transverse differential operators}
\quad\par
\label{Subsection: transverse differential operators}
Suppose that $(E, d_P^E), (E', d_P^{E'})$ are two locally free $P$-modules. For each $r \in \mathbb{N}$, define the \emph{bundle of} $r$\emph{th order transverse differential operators from} $E$ \emph{to} $E'$ as
\begin{equation*}
	\widetilde{D}^r(E, E') := \operatorname{Hom}(\widetilde{J}^rE, E').
\end{equation*}
Every section $\Phi \in \Gamma(M, \widetilde{D}^r(E, E'))$ determines a morphism of sheaves
\begin{equation}
	\label{Equation: associated transverse differential operators}
	\Gamma_P(-, E) \to \Gamma(-, E'), \quad s \mapsto \langle \Phi, \mathfrak{j}^r s \rangle.
\end{equation}
Here, $\Gamma_P(-, E)$ is the space of $d_P^E$-closed sections of $E$ and $\langle \, \cdot\,, \, \cdot\, \rangle$ denotes the natural pairing
\begin{equation*}
	\Gamma(M, \widetilde{D}^r(E, E')) \times \Gamma(M, \widetilde{J}^rE) \to \Gamma(M, E').
\end{equation*}
There is a natural surjective map $D^r(E, E') \to \widetilde{D}^r(E, E')$. Applying $\operatorname{Hom}(-, E')$ to \eqref{Equation: short exact sequence of transverse jet bundles}, we obtain a short exact sequence of complex vector bundles over $M$:
\begin{center}
	\begin{tikzcd}
		0 \ar[r] & \widetilde{D}^{r-1}(E, E') \ar[r] & \widetilde{D}^r(E, E') \ar[r, "\widetilde{\operatorname{ps}}_r^{E, E'}"] & \operatorname{Sym}^r Q \otimes \operatorname{Hom}(E, E') \ar[r] & 0
	\end{tikzcd}
\end{center}
For any $\Phi \in \Gamma(M, \widetilde{D}^r(E, E'))$, we call $\widetilde{\operatorname{ps}}_r^{E, E'}(\Phi)$ the $r$\emph{th transverse principal symbol} of $\Phi$.

\begin{proposition}
	\label{Proposition 2.6}
	For each $r \in \mathbb{N}$, the following diagram is commutative:
	\begin{center}
		\begin{tikzcd}
			D^{r-1}(E, E') \ar[r] \ar[d] &  D^r(E, E') \ar[rr, "\operatorname{ps}_r^{E, E'}"] \ar[d] && \operatorname{Sym}^r TM_\mathbb{C} \otimes \operatorname{Hom}(E, E') \ar[d]\\
			\widetilde{D}^{r-1}(E, E') \ar[r] &  \widetilde{D}^r(E, E') \ar[rr, "\widetilde{\operatorname{ps}}_r^{E, E'}"] && \operatorname{Sym}^rQ \otimes \operatorname{Hom}(E, E')
		\end{tikzcd}
	\end{center}
	where $\operatorname{ps}_r^{E, E'}$ denotes the $r$th principal symbol map.
\end{proposition}
\begin{proof}
	By Proposition \ref{Proposition: transverse jet bundles are smooth subbundles}, the following diagram is commutative:
	\begin{center}
		\begin{tikzcd}
			\operatorname{Sym}^rQ^* \otimes E \ar[r] \ar[d] & \widetilde{J}^rE \ar[rr] \ar[d] && \widetilde{J}^{r-1}E \ar[d]\\
			\operatorname{Sym}^r T^*M_\mathbb{C} \otimes E \ar[r] &  J^rE \ar[rr] && J^{r-1}E
		\end{tikzcd}
	\end{center}
	Applying $\operatorname{Hom}(-, E')$ to the above commutative diagram, we are done.
\end{proof}

\begin{definition}
	The \emph{bundle} $\widetilde{D}(E, E')$ \emph{of transverse differential operators from} $E$ \emph{to} $E'$ is defined as the inductive limit of the following sequence of vector bundle morphisms:
	\begin{center}
		\begin{tikzcd}
			\operatorname{Hom}(E, E') = \widetilde{D}^0(E, E') \ar[r] & \widetilde{D}^1(E, E') \ar[r] & \cdots \ar[r] & \widetilde{D}^r(E, E') \ar[r] & \cdots
		\end{tikzcd}
	\end{center}
\end{definition}

Hence, $\widetilde{D}(E, E')$ is an ind-vector bundle and its associated graded ind-vector bundle is
\begin{equation*}
	\operatorname{Sym} Q \otimes \operatorname{Hom}(E, E') = \varinjlim_r \left( \operatorname{Sym}^{\leq r} Q \otimes \operatorname{Hom}(E, E') \right).
\end{equation*}

Let $(E'', d_P^{E''})$ be another locally free $P$-module and $r, r' \in \mathbb{N}$. There is a binary operation
\begin{equation*}
	\Gamma(M, D^{r'}(E', E'')) \times \Gamma(M, D^r(E, E')) \to \Gamma(M, D^{r+r'}(E, E'')), \quad (\Phi', \Phi) \mapsto \Phi' \circ \Phi,
\end{equation*}
which corresponds to the usual composition of  the associated differential operators (c.f. Subsection \ref{Subsection: notations}). Explicitly, $\Phi' \circ \Phi$ is defined by the formula
\begin{equation}
	\label{Equation: composition formula}
	\langle \Phi' \circ \Phi, \sigma \rangle = \langle \Phi', \mathfrak{j}_x^{r'} \langle \Phi, \mathfrak{j}^r s\rangle \rangle,
\end{equation}
for all $x \in M$ and $\sigma \in J_x^{r+r'}E$, where $s$ is any local smooth section of $E$ such that $\sigma = \mathfrak{j}_x^{r+r'} s$. Indeed, this binary operation descends to another binary operation
\begin{equation}
	\label{Equation: composition of differential with transverse differential operators}
	\Gamma(M, D^{r'}(E', E'')) \times \Gamma(M, \widetilde{D}^r(E, E')) \to \Gamma(M, \widetilde{D}^{r+r'}(E, E'')), \quad (\Phi', \Phi) \mapsto \Phi' \circ \Phi.
\end{equation}
The morphism \eqref{Equation: associated transverse differential operators} associated with $\Phi' \circ \Phi$ agrees with the usual composition of the differential operator corresponding to $\Phi'$ with the morphism \eqref{Equation: associated transverse differential operators} associated with $\Phi$. The well definedness of \eqref{Equation: composition of differential with transverse differential operators} can be verified by working in local trivializations of $J^rE, J^{r'}E'$ of the form \eqref{Equation: local trivialization of jet bundle}.

\subsection{Canonical $P$-module structures on transverse structures}
\quad\par
\label{Subsection: Canonical module structures on transverse structures}
In this subsection we describe the canonical $P$-module structures carried by the bundles of transverse differential operators and by the transverse jet bundles. This will allow us to introduce the notion of \emph{transverse differential operators}. These constructions will only be used later, in Subsection \ref{Subsection: non formal quantizable functions}, where they play a key role in proving our second main result. Readers already familiar with the corresponding material may safely skip this subsection.\par
For each $r \in \mathbb{N}$, the bundle $\widetilde{D}^r(E, E')$ carries a canonical flat $P$-connection:
\begin{equation}
	d_P^{\widetilde{D}^r(E, E')}(\Phi) = d_P^{E'} \circ \Phi.
\end{equation}
This is well defined because for every $Y \in \Gamma(M, P)$, $(\iota_Y \circ d_P^{E'}) \circ \Phi$ is again a section of $\widetilde{D}^r(E, E')$; this will follow from Proposition \ref{Lemma 3.4} in the next section.\par
Similarly, the transverse $r$-jet bundle $\widetilde{J}^rE$ carries a canonical flat $P$-connection $d_P^{\widetilde{J}^rE}$ defined by the analogue of \eqref{Equation: formula of Grothendieck connection}: for all $Y \in \Gamma(M, P)$, $\sigma \in \Gamma(M, \widetilde{J}^rE)$ and $\Phi \in \Gamma(M, \widetilde{D}^r(E, \underline{\mathbb{C}}))$,
\begin{equation}
	\left\langle \Phi, \iota_Y d_P^{\widetilde{J}^rE} \sigma \right\rangle = \mathcal{L}_Y \langle \Phi, \sigma \rangle - \langle \mathcal{L}_Y \circ \Phi, \sigma \rangle.
\end{equation}
From this one obtains, for all $Y \in \Gamma(M, P)$, $\Phi \in \Gamma(M, \widetilde{D}^r(E, E'))$ and $\sigma \in \Gamma(M, \widetilde{J}^rE)$,
\begin{equation}
	\label{Equation: two descriptions of flat P-connections}
	\left\langle \iota_Y d_P^{\widetilde{D}^r(E, E')}(\Phi), \sigma \right\rangle = \iota_Y d_P^{E'} \langle \Phi, \sigma \rangle - \left\langle \Phi, \iota_Y d_P^{\widetilde{J}^rE} \sigma \right\rangle.
\end{equation}

\begin{proposition}
	Let $r \in \mathbb{N}$ and suppose that $\Phi \in \Gamma(M, \widetilde{D}^r(E, E'))$ is $d_P^{\widetilde{D}^r(E, E')}$-closed. Then \eqref{Equation: associated transverse differential operators} defines a morphism of sheaves $\Gamma_P(-, E) \to \Gamma_P(-, E')$.
\end{proposition}
\begin{proof}
	Let $s$ be a local $d_P^E$-closed smooth section of $E$. Since the infinite jet of $s$ is $\nabla^{\operatorname{G}, E}$-closed, we obtain $d_P^{\widetilde{J}^rE} (\mathfrak{j}^r s) = 0$. It implies that $\iota_Y d_P^{E'} \langle \Phi, \mathfrak{j}^r s \rangle = \left\langle \iota_Y d_P^{\widetilde{D}^r(E, E')}(\Phi), \mathfrak{j}^r s \right\rangle = 0$ for all $Y \in \Gamma(M, P)$. Therefore, $\langle \Phi, \mathfrak{j}^r s \rangle$ is $d_P^{E'}$-closed.
\end{proof}

\begin{definition}
	An $r$\emph{th order transverse differential operator from} $E$ \emph{to} $E'$ (with respect to $P$) is a $d_P^{\widetilde{D}^r(E, E')}$-closed section of $\widetilde{D}^r(E, E')$.
\end{definition}

Finally, let $(E'', d_P^{E''})$ be another locally free $P$-module , let $r, r' \in \mathbb{N}$, let $\Phi \in \Gamma(M, \widetilde{D}^r(E, E'))$ be $d_P^{\widetilde{D}^r(E, E')}$-closed and $\Phi' \in \Gamma(M, \widetilde{D}^{r'}(E', E''))$ be $d_P^{\widetilde{D}^{r'}(E', E'')}$-closed. Then one may define a $d_P^{\widetilde{D}^{r+r'}(E, E'')}$-closed section
\begin{equation}
	\label{Equation: composition of transverse differential operators}
	\Phi' \circ \Phi \in \Gamma(M, \widetilde{D}^{r+r'}(E', E''))
\end{equation}
by the formula \eqref{Equation: composition formula}: for each $x \in M$ and $\sigma \in \widetilde{J}_x^{r+r'}E$, choose any local $d_P^E$-closed section $s$ of $E$ with $\sigma = \mathfrak{j}_x^{r+r'} s$ (whose existence is ensured by Lemma \ref{Lemma: coordinate description of transverser jet bundle}) and set $\langle \Phi' \circ \Phi, \sigma \rangle$ accordingly. We can argue, by working in local trivializations of the form \eqref{Equation: local trivialization of jet bundle}, that this is well defined. Moreover, the morphism of sheaves \eqref{Equation: associated transverse differential operators} associated with $\Phi' \circ \Phi$ coincides with the composition of the morphisms associated with $\Phi'$ and with $\Phi$.

\section{Transverse Poincar\'e--Birkhoff--Witt theorem}
\label{Section: transverse pbw theorem}
Continuing from the previous section, we now establish Theorem \ref{Third main result} in Subsection \ref{Subsection: Transverse pbw maps}. The theorem asserts the existence of a distinguished splitting
\begin{equation}
	\label{Equation: splitting of transverse differential operators}
	\widetilde{D}(E, E') \cong \operatorname{Sym} Q \otimes \operatorname{Hom}(E, E'),
\end{equation}
obtained upon choosing the auxiliary data \eqref{Equation: auxiliary data for pbw isomorphism}, namely a splitting $TM_\mathbb{C} = P \oplus Q$, a torsion-free connection $\nabla^Q$ on $Q$, and a connection $\nabla^E$ on $E$ extending $d_P^E$.

\begin{remark}
	\label{Remark: torsion free connection on quotient bundles}
	A connection $\nabla^Q$ on $Q$ is called \emph{torsion-free} (see Subsection 5.1 of \cite{LauStiXu2021}) if
	\begin{equation*}
		\nabla_X^Q \pi(X') - \nabla_{X'}^Q \pi(X) = \pi([X, X'])
	\end{equation*}
	for all $X, X' \in \Gamma(M, TM_\mathbb{C})$. This condition ensures that $\nabla^Q$ extends the Bott $P$-connection $d_P^Q$ on $Q$, i.e. for all $Y \in \Gamma(M, P)$ and $X \in \Gamma(M, TM_\mathbb{C})$ with $Z = \pi(X)$,
	\begin{equation*}
		\nabla_Y^Q Z = \pi([Y, X]) =: \iota_Y (d_P^Q Z).
	\end{equation*}
	Since $P$ is involutive, $(TM_\mathbb{C}, P)$ forms a Lie pair in the sense of \cite{LauStiXu2012}. By Proposition 5.3 of \cite{LauStiXu2021}, such a torsion-free connection on $Q$ always exists.
\end{remark}

In the special case $E = E' = \underline{\mathbb{C}}$, such a splitting was previously established by Laurent-Gengoux--Sti\'enon--Xu \cite{LauStiXu2021} using the general framework of Lie pairs. Since our applications to deformation quantization require only Lie pairs of the form $(TM_\mathbb{C}, P)$, we restrict to this setting when constructing the splitting \eqref{Equation: splitting of transverse differential operators} for general $E, E'$.\par
In Subsection \ref{Subsection: splitting of transverse jet bundles}, we identify the transverse jet bundle $\widetilde{J}E$ with $\widehat{\operatorname{Sym}} Q^* \otimes E$ via \eqref{Equation: splitting of transverse differential operators}. In Subsection \ref{Subsection: a decomposition of the Grothendieck connection}, we study the induced decomposition of the Grothendieck connection $\nabla^{\operatorname{G}, E}$ on $\widetilde{J}E$, which plays a crucial role in our construction of a special class of deformation quantizations. Consequently, the Chevalley-Eilenberg cochain complex $(C^{0, *}(M, E), d_P^E)$ is quasi-isomorphic to the complex $( \Omega^*(M, \widetilde{J}E), \nabla^{\operatorname{G}, E} )$, which we establish in Subsection \ref{Subsection: Chevalley--Eilenberg cohomology}.\par
In this section, we identify $Q$ with its image under the chosen splitting $Q \to TM_\mathbb{C}$. Then
\begin{equation*}
	\Omega^*(M, \widehat{\operatorname{Sym}} Q^* \otimes E) = C^{*, *}(M, \widehat{\operatorname{Sym}} Q^* \otimes E),
\end{equation*}
where $C^{q, p}(M, \widehat{\operatorname{Sym}} Q^* \otimes E) := \Gamma(M, \textstyle \bigwedge^q Q^* \otimes \bigwedge^p P^* \otimes \widehat{\operatorname{Sym}} Q^* \otimes E)$ for each $p, q \in \mathbb{N}$. We also define $\mathcal{C}^\infty(M)$-linear operators $\delta_Q, \delta_Q^{-1}, \pi_Q$ on $\Omega^*(M, \widehat{\operatorname{Sym}} Q^* \otimes E)$ as follows: for any $\alpha \in C^{q, p}(M, \operatorname{Sym}^l Q^* \otimes E)$, set
\begin{align*}
	\delta_Q \alpha := v^i \wedge \frac{\partial \alpha}{\partial u^i}, \quad \delta_Q^{-1} \alpha := 
	\begin{cases} 
		\frac{1}{q+l} u^i \cdot \iota_{v_i} \alpha, & q+l>0,\\
		0, & q+l=0,
	\end{cases}
	\quad 
	\pi_Q(\alpha) := 
	\begin{cases}
		0, & q+l>0,\\
		\alpha, & q+l=0.
	\end{cases}
\end{align*}
Here, $(v_1, ..., v_q)$ is any local frame of $Q$, and the notation $v^i, u^i$ follows the conventions specified in Subsection \ref{Subsection: notations}. These operators satisfy the formal Hodge decomposition:
\begin{equation}
	\label{Equation: Hodge decomposition for Dolbeault differential}
	\operatorname{Id} - \pi_Q = \delta_Q \circ \delta_Q^{-1} + \delta_Q^{-1} \circ \delta_Q.
\end{equation}

\subsection{Poincar\'e--Birkhoff--Witt maps for transverse differential operators}
\quad\par
\label{Subsection: Transverse pbw maps}
The connection $\nabla^E: \Gamma(M, TM_\mathbb{C}) \to \Gamma(M, D^1(E, E))$ restricts to a $\mathcal{C}^\infty(M)$-linear map
\begin{equation}
	\label{Equation: restriction of connection}
	\Gamma(M, Q) \to \Gamma(M, \widetilde{D}^1(E, E)).
\end{equation}
We extend \eqref{Equation: restriction of connection} to a morphism of sheaves of $\mathcal{C}_M^\infty$-modules
\begin{equation*}
	\Gamma(-, \operatorname{Sym} Q) \to \Gamma(-, \widetilde{D}(E, E)),
\end{equation*}
still denoted by the symbol $\nabla^E$, by a recursive formula: for $f \in \mathcal{C}^\infty(M)$, $\nabla_f^E = f\operatorname{Id}_E$; for all $r \in \mathbb{Z}^+$ and $Z_1, ..., Z_r \in \Gamma(M, Q)$,
\begin{equation*}
	\nabla_{Z_1 \cdots Z_r}^E = \frac{1}{r} \sum_{i=1}^r \left( \nabla_{Z_i}^E \circ \nabla_{Z_1 \cdots \widehat{Z_i} \cdots Z_r}^E - \nabla_{\nabla_{Z_i}^Q \left( Z_1 \cdots \widehat{Z_i} \cdots Z_r \right)}^E \right).
\end{equation*}
We now introduce the following construction, using the chosen data \eqref{Equation: auxiliary data for pbw isomorphism}.

\begin{definition}
	The \emph{Poincar\'e--Birkhoff--Witt map} is the morphism of $\mathcal{C}_M^\infty$-modules
	\begin{equation*}
		\operatorname{pbw}^{E, E'}: \Gamma(-, \operatorname{Sym} Q \otimes E^* \otimes E') \to \Gamma(-, \widetilde{D}(E, E'))
	\end{equation*}
	defined as follows. For all $r \in \mathbb{N}$, $G \in \Gamma(M, \operatorname{Sym}^r Q)$ and $A \in \Gamma(M, E^* \otimes E')$,
	\begin{equation*}
		\operatorname{pbw}^{E, E'}(G \otimes A) = A \circ \nabla_G^E.
	\end{equation*}
\end{definition}

\begin{theorem}[Theorem \ref{Third main result}]
	\label{Theorem: pbw theorem}
	Let $M$ be a smooth manifold equipped with a Nirenberg integrable complex distribution $P$. Then for any locally free $P$-modules $E, E'$, the map
	\begin{equation*}
		\operatorname{pbw}^{E, E'}: \Gamma(-, \operatorname{Sym} Q \otimes E^* \otimes E') \to \Gamma(-, \widetilde{D}(E, E'))
	\end{equation*}
	is an isomorphism of filtered $\mathcal C_M^\infty$-modules whose associated graded map is the identity map on $\operatorname{Sym}Q\otimes E^*\otimes E'$.
\end{theorem}
\begin{proof}
	It is clear that $\operatorname{pbw}^{E, E'}$ sends $\Gamma(U, \operatorname{Sym}^{\leq r} Q \otimes E^* \otimes E')$ to $\Gamma(U, \widetilde{D}^r(E, E'))$ for each open subset $U$ of $M$ and each $r \in \mathbb{N}$. Thus, it suffices to prove that the composition
	\begin{center}
		\begin{tikzcd}
			\operatorname{Sym}^r Q \otimes E^* \otimes E' \ar[rr, "\operatorname{pbw}^{E, E'}"] && \widetilde{D}^r(E, E') \ar[rr, "\widetilde{\operatorname{ps}}_r^{E, E'}"] && \operatorname{Sym}^r Q \otimes E^* \otimes E'
		\end{tikzcd}
	\end{center}
	is the identity map. The statement is trivially true when $r = 0$. Let $r \geq 1$. Fix $Z \in \Gamma(M, Q)$ and $A \in \Gamma(M, E^* \otimes E')$. We can show by induction that
	\begin{align*}
		\operatorname{pbw}^{E, E'}(Z^r \otimes A) - A \circ \nabla_Z^E \circ \cdots \circ \nabla_Z^E \in \Gamma(M, \widetilde{D}^{r-1}(E, E')).
	\end{align*}
	On the other hand, we know that $\operatorname{ps}_r^{E, E'} \left( A \circ \nabla_Z^E \circ \cdots \circ \nabla_Z^E \right) = Z^r \otimes A$. Then
	\begin{equation*}
		\widetilde{\operatorname{ps}}_r^{E, E'} \left( \operatorname{pbw}^{E, E'}(Z^r \otimes A) \right) = \operatorname{ps}_r^{E, E'} \left( A \circ \nabla_Z^E \circ \cdots \circ \nabla_Z^E \right) = Z^r \otimes A
	\end{equation*}
	by the commutative diagram in Proposition \ref{Proposition 2.6}.
\end{proof}

We now establish a proposition that will be a key tool in Proposition \ref{Lemma 3.7} for analysing the decomposition of the Grothendieck connection on $\widetilde{J}E$.

\begin{proposition}
	\label{Lemma 3.4}
	Let $\nabla^{E'}$ be a connection on $E'$ extending $d_P^{E'}$. For $X \in \Gamma(M, TM_\mathbb{C})$, $r \in \mathbb{N}$ and $\phi \in \Gamma(M, \operatorname{Sym}^r Q \otimes E^* \otimes E')$, one has $\nabla_X^{E'} \circ \operatorname{pbw}^{E, E'}(\phi) \in \Gamma(M, \widetilde{D}^{r+1}(E, E'))$ and
	\begin{align*}
		\widetilde{\operatorname{ps}}_{r+1}^{E, E'} \left( \nabla_X^{E'} \circ \operatorname{pbw}^{E, E'}(\phi) \right) = & \pi(X) \cdot \phi,\\
		\widetilde{\operatorname{ps}}_r^{E, E'} \left( \nabla_X^{E'} \circ \operatorname{pbw}^{E, E'}(\phi) - \operatorname{pbw}^{E, E'}( \pi(X) \cdot \phi ) \right) = & \nabla_X^{\operatorname{Sym} Q \otimes E^* \otimes E'} \phi.
	\end{align*}
	In particular, if $Y \in \Gamma(M, P)$, then $\nabla_Y^{E'} \circ \operatorname{pbw}^{E, E'}(\phi) \in \Gamma(M, \widetilde{D}^r(E, E'))$.
\end{proposition}

The proof of Proposition \ref{Lemma 3.4} proceeds by induction on $r$. We being with the base case:

\begin{lemma}
	Let $\nabla^{E'}$ be a connection on $E'$ extending $d_P^{E'}$. Suppose that $X \in \Gamma(M, TM_\mathbb{C})$ and $A \in \Gamma(M, E^* \otimes E')$. Then for any $\sigma \in \Gamma(M, \widetilde{J}E)$,
	\begin{equation*}
		\langle \nabla_X^{E'} \circ A, \sigma \rangle = \langle \operatorname{pbw}^{E, E'}(\pi(X) \otimes A + \nabla_X^{E^* \otimes E'} A), \sigma \rangle.
	\end{equation*}
\end{lemma}
\begin{proof}
	Let $Y = X - \pi(X) \in \Gamma(M, P)$. We have $\nabla_Y^E \sigma = 0$. Then
	\begin{equation*}
		\langle \nabla_X^{E'} \circ A, \sigma \rangle = \langle A \circ \nabla_X^E, \sigma \rangle + \langle \nabla_X^{E^* \otimes E'} A, \sigma \rangle = \langle A \circ \nabla_{\pi(X)}^E, \sigma \rangle + \langle \nabla_X^{E^* \otimes E'} A, \sigma \rangle,
	\end{equation*}
	which concludes the proof.
\end{proof}

To carry out the induction step, we need the following formula.

\begin{lemma}
	\label{Lemma 3.5}
	Let $\nabla^{E'}$ be a connection on $E'$ extending $d_P^{E'}$. Suppose that $X \in \Gamma(M, TM_\mathbb{C})$, $A \in \Gamma(M, E^* \otimes E')$, $Z \in \Gamma(M, Q)$ and $r \geq 1$. Then
	\begin{equation}
		\nabla_X^{E'} \circ A \circ \nabla_{Z^r}^E = A \circ \nabla_Z^E \circ \nabla_X^E \circ \nabla_{Z^{r-1}}^E + \Phi_{X, 1} + \Phi_{X, 2},
	\end{equation}
	where
	\begin{align*}
		\Phi_{X, 1} = & (\nabla_X^{E^* \otimes E'} A) \circ \nabla_{Z^r}^E + A \circ \nabla_{[X, Z]}^E \circ \nabla_{Z^{r-1}}^E - A \circ \nabla_X^E \circ \nabla_{\nabla_Z^Q Z^{r-1}}^E,\\
		\Phi_{X, 2} = & A \circ R^E( X, Z ) \circ \nabla_{Z^{r-1}}^E.
	\end{align*}
\end{lemma}
\begin{proof}
	The following direct computations conclude the proof:
	\begin{align*}
		\nabla_X^{E'} \circ A \circ \nabla_{Z^r}^E = & A \circ \nabla_X^E \circ \nabla_{Z^r}^E + (\nabla_X^{E^* \otimes E'} A) \circ \nabla_{Z^r}^E,\\
		\nabla_X^E \circ \nabla_{Z^r}^E = & \nabla_X^E \circ \nabla_Z^E \circ \nabla_{Z^{r-1}}^{E, (r-1)} - \nabla_X^E \circ \nabla_{\nabla_Z^Q Z^{r-1}}^{E, (r-1)}\\
		= & \left( \nabla_Z^E \circ \nabla_X^E + \nabla_{[X, Z]}^E + R^E( X, Z ) \right) \circ \nabla_{Z^{r-1}}^{E, (r-1)} - \nabla_X^E \circ \nabla_{\nabla_Z^Q Z^{r-1}}^{E, (r-1)}.
	\end{align*}
\end{proof}

\begin{proof}[\myproof{Proposition}{\ref{Lemma 3.4}}]
	Consider the case for $r \geq 1$, assuming that the statement is true for $0, ..., r - 1$. Without loss of generality, we assume that $\phi = Z^r \otimes A$ for some $Z \in \Gamma(M, Q)$ and $A \in \Gamma(M, E^* \otimes E')$. Write
	\begin{equation*}
		\nabla_X^{E'} \circ \operatorname{pbw}^{E, E}(\phi) = A \circ \nabla_Z^E \circ \nabla_X^E \circ \nabla_{Z^{r-1}}^{E, (r-1)} + \Phi_{X, 1} + \Phi_{X, 2}
	\end{equation*}
	as in Lemma \ref{Lemma 3.5}. Note that $\Phi_{X, 2} \in \Gamma(M, \widetilde{D}^{r-1}(E, E'))$ and $\Phi_{X, 1} \in \Gamma(M, \widetilde{D}^r(E, E'))$. By the induction hypothesis,
	\begin{equation*}
		\widetilde{\operatorname{ps}}_r^{E, E'} \left( \Phi_{X, 1} \right) = Z^r \otimes (\nabla_X^{E^* \otimes E'} A) + (\nabla_X^Q Z) \cdot Z^{r-1} \otimes A - \nabla_Z^Q \left( Z' \cdot Z^{r-1} \right) \otimes A,
	\end{equation*}
	where $Z' = \pi(X)$. Here, we have used $\pi([X, Z]) = \nabla_X^Q Z - \nabla_Z^Q Z'$ since $\nabla^Q$ is torsion-free.\par
	For any $\Phi, \Phi' \in \Gamma(M, \widetilde{D}(E, E'))$, we write $\Phi \equiv \Phi'$ if and only if $\Phi - \Phi' \in \Gamma(M, \widetilde{D}^{r-1}(E, E'))$. By the induction hypothesis again, there exists $\Phi_{X, 4} \in \Gamma(M, \widetilde{D}^{r-2}(E, E))$ such that
	\begin{align*}
		A \circ \nabla_Z^E \circ \nabla_X^E \circ \nabla_{Z^{r-1}}^{E, (r-1)} \equiv & A \circ \nabla_Z^E \circ \left( \nabla_{Z' \cdot Z^{r-1}}^E + \nabla_{\nabla_X^Q Z^{r-1}}^{E, (r-1)} + \Phi_{X, 4} \right)\\
		\equiv & A \circ \nabla_Z^E \circ \nabla_{Z' \cdot Z^{r-1}}^E + A \circ \nabla_{Z \cdot \nabla_X^Q Z^{r-1}}^E.
	\end{align*}
	Therefore,
	\begin{equation*}
		\nabla_X^{E'} \circ \operatorname{pbw}^{E, E'}(\phi) \equiv A \circ \left( \nabla_Z^E \circ \nabla_{Z' \cdot Z^{r-1}}^E - \nabla_{\nabla_Z^Q(Z' \cdot Z^{r-1})}^E \right) + \operatorname{pbw}^{E, E'} \left( \nabla_X^{\operatorname{Sym} Q \otimes E^* \otimes E'} \phi \right).
	\end{equation*}
	It implies that for any $Y \in \Gamma(M, P)$, $\nabla_Y^{E'} \circ \operatorname{pbw}^{E, E'}(\phi) \in \Gamma(M, \widetilde{D}^r(E, E'))$ and
	\begin{equation*}
		\widetilde{\operatorname{ps}}_r^{E, E'} \left( \nabla_Y^{E'} \circ \operatorname{pbw}^{E, E'}(\phi) \right) = \nabla_Y^{\operatorname{Sym} Q \otimes E^* \otimes E'} \phi.
	\end{equation*}
	Let $Y' = X - Z' \in \Gamma(M, P)$. We can furthermore see that
	\begin{align*}
		\nabla_X^{E'} \circ \operatorname{pbw}^{E, E'}(\phi) \equiv & A \circ \nabla_{Z'}^E \circ \nabla_{Z^r}^E + (\nabla_{Z'}^{E^* \otimes E'} A) \circ \nabla_{Z^r}^E + \operatorname{pbw}^{E, E'} \left( \nabla_{Y'}^{\operatorname{Sym} Q \otimes E^* \otimes E'} \phi \right)\\
		\equiv & A \circ \left( \nabla_{Z'}^E \circ \nabla_{Z^r}^E - \nabla_{\nabla_{Z'}^Q Z^r}^E \right) + \operatorname{pbw}^{E, E'} \left( \nabla_X^{\operatorname{Sym} Q \otimes E^* \otimes E'} \phi \right).
	\end{align*}
	Recall that $\nabla_{Z' \cdot Z^r}^{E, (r+1)}$ is equal to
	\begin{align*}
		\frac{1}{r+1} \left( \nabla_{Z'}^E \circ \nabla_{Z^r}^E + r \cdot \nabla_Z^E \circ \nabla_{Z' \cdot Z^{r-1}}^E - \nabla_{\nabla_{Z'}^Q (Z^r)}^E - r \cdot \nabla_{\nabla_Z^Q(Z' \cdot Z^{r-1})}^E \right),
	\end{align*}
	Eventually, we see that from the following computation that the statement is also true for $r$:
	\begin{align*}
		\nabla_X^{E'} \circ \operatorname{pbw}^{E, E'}(\phi) \equiv & \frac{1}{r+1} \left( \nabla_X^{E'} \circ \operatorname{pbw}^{E, E'}(\phi) + r \cdot \nabla_X^{E'} \circ \operatorname{pbw}^{E, E'}(\phi) \right)\\
		\equiv & \operatorname{pbw}^{E, E'} \left( Z' \cdot \phi + \nabla_X^{\operatorname{Sym} Q \otimes E^* \otimes E'} \phi \right).
	\end{align*}
\end{proof}

\subsection{Splitting of transverse jet bundles}
\quad\par
\label{Subsection: splitting of transverse jet bundles}
The Poincar\'e--Birkhoff--Witt map $\operatorname{pbw}^{E, \underline{\mathbb{C}}}$ induces an isomorphism of $\mathcal{C}_M^\infty$-modules
\begin{equation*}
	S^E: \Gamma(-, \widehat{\operatorname{Sym}} Q^* \otimes E) \to \Gamma(-, \widetilde{J}E)
\end{equation*}
defined as follows. For any open $U \subset M$, $s \in \Gamma(U, \widehat{\operatorname{Sym}} Q^* \otimes E)$, $r \in \mathbb{N}$ and $\Phi \in \Gamma(U, \widetilde{D}^r(E, \underline{\mathbb{C}}))$,
\begin{equation}
	\label{Equation: splitting of transverse jet bundle}
	\langle \Phi, S^E(s) \rangle = \langle (\operatorname{pbw}^{E, \underline{\mathbb{C}}})^{-1} (\Phi), s \rangle,
\end{equation}
where $\langle \, \cdot\,, \, \cdot\, \rangle$ on the right hand side is the natural pairing
\begin{equation*}
	\Gamma(U, \operatorname{Sym}^{\leq r} Q \otimes E^*) \times \Gamma(U, \widehat{\operatorname{Sym}} Q^* \otimes E) \to \mathcal{C}^\infty(U).
\end{equation*}
The notation $S^E$ emphasizes that this map provides a splitting of the projection from $\widetilde{J}E$ onto its associated graded pro-vector bundle $\widehat{\operatorname{Sym}} Q^* \otimes E$.\par
While the space $\Gamma(M, \widetilde{J}E)$ naturally carries a $\Gamma(M, \widetilde{J}\underline{\mathbb{C}})$-module structure, the other space $\Gamma(M, \widehat{\operatorname{Sym}} Q^* \otimes E)$ carries a natural $\Gamma(M, \widehat{\operatorname{Sym}} Q^*)$-module structure, which satisfies
\begin{equation}
	\label{Equation: module structure on symmetric algebra}
	\langle Z_1 \cdots Z_r \otimes A, a \cdot s \rangle = \sum_{i=0}^r \binom{r}{i} \sum_{\sigma \in \mathcal{S}_r} \frac{1}{r!} \langle Z_{\sigma(1)} \cdots Z_{\sigma(i)}, a \rangle \cdot \langle Z_{\sigma(i+1)} \cdots Z_{\sigma(r)} \otimes A, s \rangle,
\end{equation}
for all $a \in \Gamma(M, \widehat{\operatorname{Sym}} Q^*)$, $s \in \Gamma(M, \widehat{\operatorname{Sym}} Q^* \otimes E)$, $Z_1, ..., Z_r \in \Gamma(M, Q)$ and $A \in \Gamma(M, E^*)$, where $\mathcal{S}_r$ denotes the symmetric group on $r$ elements. The following proposition asserts that the splittings $S^{\underline{\mathbb{C}}}$ and $S^E$ are compatible with these module structures.

\begin{proposition}
	\label{Proposition: pbw maps preserve module structure}
	For any $a \in \Gamma(M, \widehat{\operatorname{Sym}} Q^*)$ and $s \in \Gamma(M, \widehat{\operatorname{Sym}} Q^* \otimes E)$,
	\begin{equation*}
		S^E(a \cdot s) = S^{\underline{\mathbb{C}}}(a) \cdot S^E(s).
	\end{equation*}
\end{proposition}
\begin{proof}
	It suffices to verify the identity after pairing with $\phi = Z^r \otimes A \in \Gamma(M,\operatorname{Sym}^r Q \otimes E^*)$. By definition, $\left\langle \operatorname{pbw}^{E,\underline{\mathbb{C}}}(\phi), S^E(s') \right\rangle = \langle \phi, s' \rangle$ for all $s' \in \Gamma(M, \widehat{\operatorname{Sym}} Q^* \otimes E)$. Then by \eqref{Equation: module structure on symmetric algebra},
	\begin{align*}
		\left\langle \operatorname{pbw}^{E, \underline{\mathbb{C}}} (\phi), S^E(a \cdot s) \right\rangle = & \sum_{i=0}^r \binom{r}{i} \langle \operatorname{pbw}^{\underline{\mathbb{C}}, \underline{\mathbb{C}}} (Z^i), S^{\underline{\mathbb{C}}} (a) \rangle \cdot \langle \operatorname{pbw}^{E, \underline{\mathbb{C}}} (Z^{r-i} \otimes A), S^E(s) \rangle.
	\end{align*}
	To complete the proof, we show that the right hand side equals $\left\langle \operatorname{pbw}^{E, \underline{\mathbb{C}}} (\phi), S^{\underline{\mathbb{C}}} (a) \cdot S^E(s) \right\rangle$. It suffices to prove that for all $Z_1, ..., Z_r \in \Gamma(M, Q)$, $f \in \mathcal{C}^\infty(M)$ and $s \in \Gamma(M, E)$,
	\begin{equation*}
		\nabla_{Z_1 \cdots Z_r}^E (f \cdot s) = \sum_{i=0}^r \binom{r}{i} \sum_{\sigma \in \mathcal{S}_r} \frac{1}{r!} \left( \mathcal{L}_{Z_{\sigma(1)} \cdots Z_{\sigma(i)}} f \right) \cdot \left( \nabla_{Z_{\sigma(i+1)} \cdots Z_{\sigma(r)}}^E s \right),
	\end{equation*}
	where $\mathcal{L}_{Z_{\sigma(1)} \cdots Z_{\sigma(i)}}$ is defined analogously to $\nabla_{Z_{\sigma(1)} \cdots Z_{\sigma(i)}}^E$ with respect to Lie derivatives of functions. The claim is immediate for $r = 0$.\par
	Now, assume that our claim is true for $r \in \mathbb{N}$. Then we obtain
	\begin{align*}
		\nabla_{\nabla_Z^Q Z^r}^E (f \cdot s) = & r \sum_{i=1}^r \binom{r}{i} \cdot \frac{i \cdot (r-1)!}{r!} \left( \mathcal{L}_{(\nabla_Z^Q Z) \cdot Z^{i-1}} f \right) \cdot \left( \nabla_{Z^{r-i}}^E s \right)\\
		& + r \sum_{i=0}^{r-1} \binom{r}{i} \cdot \frac{(r-i) \cdot (r-1)!}{r!} \left( \left( \mathcal{L}_{Z^i} f \right) \cdot \left( \nabla_{(\nabla_Z^Q Z) \cdot Z^{r-i-1}}^E s \right) \right)\\
		= & \sum_{i=1}^r \binom{r}{i} \left( \mathcal{L}_{\nabla_Z^Q (Z^i)} f \right) \cdot \left( \nabla_{Z^{r-i}}^E s \right) + \sum_{i=0}^{r-1} \binom{r}{i} \left( \mathcal{L}_{Z^i} f \right) \cdot \left( \nabla_{\nabla_Z^Q (Z^{r-i})}^E s \right),
	\end{align*}
	using the equalities $\nabla_Z^Q (Z^i) = i \cdot (\nabla_Z^Q Z) \cdot Z^{i-1}$ for $i = 0, ..., r$. On the other hand,
	\begin{equation*}
		\nabla_Z^E \circ \nabla_{Z^r}^E (f \cdot s) = \sum_{i=0}^r \binom{r}{i} \left( \left( \mathcal{L}_Z \circ \mathcal{L}_{Z^i} f \right) \cdot \left( \nabla_{Z^{r-i}}^E s \right) + \left( \mathcal{L}_{Z^i} f \right) \cdot \left( \nabla_Z^E \circ \nabla_{Z^{r-i}}^E s \right) \right).
	\end{equation*}
	Under our convention, $\mathcal{L}_Z \circ \mathcal{L}_{Z^0} f = \mathcal{L}_Z f$ and $\nabla_Z^E \circ \nabla_{Z^0}^E s = \nabla_Z^E s$. Therefore,
	\begin{align*}
		\nabla_{Z^{r+1}}^E (f \cdot s) = & \sum_{i=0}^r \binom{r}{i} \left( \mathcal{L}_{Z^{i+1}} f \right) \cdot \left( \nabla_{Z^{r-i}}^E s \right) + \sum_{i=0}^r \binom{r}{i} \left( \mathcal{L}_{Z^i} f \right) \cdot \left( \nabla_{Z^{r-i+1}}^E s \right)\\
		= & \sum_{i=1}^{r+1} \binom{r}{i-1} \left( \mathcal{L}_{Z^i} f \right) \cdot \left( \nabla_{Z^{r+1-i}}^E s \right) + \sum_{i=0}^r \binom{r}{i} \left( \mathcal{L}_{Z^i} f \right) \cdot \left( \nabla_{Z^{r+1-i}}^E s \right)\\
		= & \sum_{i=0}^{r+1} \binom{r+1}{i} \left( \mathcal{L}_{Z^i} f \right) \cdot \left( \nabla_{Z^{r+1-i}}^E s \right).
	\end{align*}
	We have proved our claim by induction.
\end{proof}

\subsection{A decomposition of the Kapranov connection}
\quad\par
\label{Subsection: a decomposition of the Grothendieck connection}
We obtain a flat connection $\nabla^{\operatorname{K}, E}$ on $\widehat{\operatorname{Sym}} Q^* \otimes E$ as the pullback of $\nabla^{\operatorname{G}, E}$ along $S^E$, and call it the \emph{Kapranov connection}. This connection gives a unique flat lift of each $P$-polarized section of $E$, as stated later in Proposition \ref{Proposition: qausi-isomorphism for polarized sections}.

This subsection studies the decomposition of $\nabla^{\operatorname{K},E}$ with respect to a decreasing filtration on $\Omega^*(M,\widehat{\operatorname{Sym}}Q^*\otimes E)$. The decomposition serves two purposes: it enters the proof of Proposition \ref{Proposition: qausi-isomorphism for polarized sections}, and, when $E=L$ is the prequantum line bundle in Setup \ref{Setup: one}, it gives rise to a solution of the Fedosov equation (Proposition \ref{Proposition: alternative desription of Grothendieck connection}), from which our deformation quantization is constructed.

By \eqref{Equation: splitting of transverse jet bundle} and \eqref{Equation: formula of Grothendieck connection}, we have the following explicit formula: for any $s \in \Gamma(M, \widehat{\operatorname{Sym}} Q^* \otimes E)$, $r \in \mathbb{N}$, $\phi \in \Gamma(M, \operatorname{Sym}^{\leq r} Q \otimes E^*)$ and $X \in \Gamma(M, TM_\mathbb{C})$,
\begin{equation}
	\label{Equation: pullback of Grothendieck connections}
	\langle \phi, \nabla_X^{\operatorname{K}, E} s \rangle = \mathcal{L}_X \langle \phi, s \rangle - \langle \mathcal{L}_X \circ \operatorname{pbw}^{E, \underline{\mathbb{C}}} (\phi), S^E(s) \rangle.
\end{equation}

We will show that the `leading term' of $\nabla^{\operatorname{K}, E}$ is
\begin{equation*}
	\nabla^{\widehat{\operatorname{Sym}} Q^* \otimes E} - \delta_Q,
\end{equation*}
where $\nabla^{\widehat{\operatorname{Sym}} Q^* \otimes E}$ denotes the connection on $\widehat{\operatorname{Sym}} Q^* \otimes E$ induced by the pair $\nabla^Q, \nabla^E$. Before stating the proposition, it is convenient to introduce a $\mathcal{C}^\infty(M)$-linear map
\begin{equation}
	\Theta^E := \nabla^{\operatorname{K}, E} - \left( \nabla^{\widehat{\operatorname{Sym}} Q^* \otimes E} - \delta_Q \right): \Gamma(M, \widehat{\operatorname{Sym}} Q^* \otimes E) \to \Omega^1(M, \widehat{\operatorname{Sym}} Q^* \otimes E).
\end{equation}

\begin{proposition}
	\label{Lemma 3.7}
	Let $r \in \mathbb{N}$. Then
	\begin{itemize}
		\item $\nabla^{\operatorname{K}, E}$ sends $\Gamma(M, \widehat{\operatorname{Sym}}^{\geq r} Q^* \otimes E)$ to $\Omega^1(M, \widehat{\operatorname{Sym}}^{\geq (r-1)} Q^* \otimes E)$; and
		\item $\Theta^E$ sends $\Gamma(M, \widehat{\operatorname{Sym}}^{\geq r} Q^* \otimes E)$ to $\Omega^1(M, \widehat{\operatorname{Sym}}^{\geq (r+1)} Q^* \otimes E)$.
	\end{itemize}
\end{proposition}
\begin{proof}
	Let $X \in \Gamma(M, TM_\mathbb{C})$ and $s \in \Gamma(M, \widehat{\operatorname{Sym}}^{\geq r} Q^* \otimes E)$. For any $\phi \in \Gamma(M, \operatorname{Sym}^{\leq r} Q \otimes E^*)$, by Proposition \ref{Lemma 3.4}, $\mathcal{L}_X \circ \operatorname{pbw}^{E, \underline{\mathbb{C}}}(\phi) \in \Gamma(M, \widetilde{D}^{r+1}(E, \underline{\mathbb{C}}))$ and
	\begin{equation*}
		(\operatorname{pbw}^{E, \underline{\mathbb{C}}})^{-1} \circ \mathcal{L}_X \circ \operatorname{pbw}^{E, \underline{\mathbb{C}}}(\phi) - \left( \pi(X) \cdot \phi + \nabla_X^{\operatorname{Sym} Q \otimes E^*} \phi \right) \in \Gamma(M, \operatorname{Sym}^{\leq (r-1)} Q \otimes E^*).
	\end{equation*}
	Then by degree reason,
	\begin{align*}
		\langle \phi, \nabla_X^{\operatorname{K}, E} s \rangle = & \mathcal{L}_X \langle \phi, s \rangle - \langle (\operatorname{pbw}^{E, \underline{\mathbb{C}}})^{-1} \circ \mathcal{L}_X \circ \operatorname{pbw}^{E, \underline{\mathbb{C}}} (\phi), s \rangle\\
		= & \mathcal{L}_X \langle \phi, s \rangle - \langle \pi(X) \cdot \phi + \nabla_X^{\operatorname{Sym} Q \otimes E^*} \phi, s \rangle\\
		= & \left\langle \phi, \nabla_X^{\widehat{\operatorname{Sym}} Q^* \otimes E} s \right\rangle - \left\langle \phi, \iota_X (\delta_Qs) \right\rangle.
	\end{align*}
	Here $\iota_X(\delta_Qs)$ denotes the contraction of the $1$-form $\delta_Qs$ with $X$. Then $\langle \phi, \iota_X\Theta^E s \rangle = 0$. In particular, if $\phi \in \Gamma(M, \operatorname{Sym}^{\leq (r-2)} Q \otimes E^*)$, then by degree reason again, $\langle \phi, \nabla_X^{\operatorname{K}, E} s \rangle = 0$.
\end{proof}

\begin{remark}
	\label{Remark: Kaparanov connection in trivial bundle case}
	When $(E, \nabla^E) = \underline{\mathbb{C}}$, we omit the explicit reference to the trivial bundle and write $\nabla^{\operatorname{K}}$, $\nabla^{\operatorname{G}}$, and $\Theta$ in place of $\nabla^{\operatorname{K}, \underline{\mathbb{C}}}$, $\nabla^{\operatorname{G}, \underline{\mathbb{C}}}$, and $\Theta^{\underline{\mathbb{C}}}$ respectively. Indeed, $\Theta$ coincides with the dual structure of (48) in \cite{LauStiXu2021} for the Lie pair $(TM_\mathbb{C}, P)$.
\end{remark}

According to Lemmas 5.14 and 5.17 in \cite{LauStiXu2021}, $\Theta(1) = 0$ and for all $a, b \in \Gamma(M, \widehat{\operatorname{Sym}} Q^*)$,
\begin{equation*}
	\Theta(a \cdot b) = (\Theta a) \cdot b + a \cdot (\Theta b).
\end{equation*}
Using the above properties of $\Theta$ and Proposition \ref{Lemma 3.7}, there exists a unique element
\begin{equation}
	\label{Equation A.6}
	I \in \Omega^1(M, \operatorname{Hom}(Q^*, \widehat{\operatorname{Sym}}^{\geq 2} Q^*))
\end{equation}
such that $\Theta$ is the extension of $I$ as a derivation on $\widehat{\operatorname{Sym}} Q^*$. We also denote the operator
\begin{equation*}
	\Theta \otimes \operatorname{Id}_E: \Gamma(M, \widehat{\operatorname{Sym}} Q^* \otimes E) \to \Omega^1(M, \widehat{\operatorname{Sym}} Q^* \otimes E).
\end{equation*}
by $I$ by abuse of notation. This operator facilitates a further decomposition of $\Theta^E$.

\begin{theorem}
	There exists a unique element
	\begin{equation}
		\label{Equation A.7}
		F^E \in \Omega^1(M, \widehat{\operatorname{Sym}}^{\geq 1} Q^* \otimes \operatorname{End} E)
	\end{equation}
	such that for all $s \in \Gamma(M, \widehat{\operatorname{Sym}} Q^* \otimes E)$,
	\begin{equation}
		\label{Equation A.8}
		\nabla^{\operatorname{K}, E} (s) = \nabla^{\widehat{\operatorname{Sym}} Q^* \otimes E} s - \delta_Q s + I(s) + F^E (s),
	\end{equation}
	with $F^E (s)$ defined via the multiplication in $\widehat{\operatorname{Sym}} Q^*$ and the natural action of $\operatorname{End} E$ on $E$.
\end{theorem}
\begin{proof}
	We claim that $F^E \in \Omega^1(M, \widehat{\operatorname{Sym}} Q^* \otimes \operatorname{End} E)$ defined as follows satisfies \eqref{Equation A.8}:
	\begin{equation*}
		F^E(s) := (\Theta^E - I)(s) \in \Omega^1(M, \widehat{\operatorname{Sym}} Q^* \otimes E) \quad \text{for any } s \in \Gamma(M, E).
	\end{equation*}
	By Proposition \ref{Lemma 3.7}, both $\Theta^E$ and $I$ send $\Gamma(M, \widehat{\operatorname{Sym}} Q^* \otimes E)$ to $\Omega^1(M, \widehat{\operatorname{Sym}}^{\geq 1} Q^* \otimes E)$. Thus, it follows from our claim that $F^E \in \Omega^1(M, \widehat{\operatorname{Sym}}^{\geq 1} Q^* \otimes \operatorname{End} E)$.\par
	To prove our claim, it suffices to show that
	\begin{equation}
		\label{Equation: theta and I}
		(\Theta^E - I)(a \cdot s) = a \cdot (\Theta^E - I)(s),
	\end{equation}
	for all $a \in \Gamma(M, \widehat{\operatorname{Sym}} Q^*)$ and $s \in \Gamma(M, \widehat{\operatorname{Sym}} Q^* \otimes E)$. Indeed, the following equalities hold:
	\begin{align*}
		\nabla^{\operatorname{K}, E}(a \cdot s) = &(\nabla^{\operatorname{K}} a) \cdot s + a \cdot \nabla^{\operatorname{K}, E} s,\\
		\nabla^{\widehat{\operatorname{Sym}} Q^* \otimes E}(a \cdot s) = & (\nabla^{\widehat{\operatorname{Sym}} Q^*} a) \cdot s + a \cdot \nabla^{\widehat{\operatorname{Sym}} Q^* \otimes E} s,\\
		\delta_Q(a \cdot s) = & (\delta_Q a) \cdot s + a \cdot \delta_Qs.
	\end{align*}
	The second the the third equalities are obvious. For the first one, we need to use a basic fact about Grothendieck connections: for all $\xi \in \Gamma(M, J\underline{\mathbb{C}})$ and $\sigma \in \Gamma(M, JE)$,
	\begin{equation*}
		\nabla^{\operatorname{G}}(\xi \cdot \sigma) = ( \nabla^{\operatorname{G}} \xi ) \cdot \sigma + \xi \cdot ( \nabla^{\operatorname{G}, E} \sigma).
	\end{equation*}
	By this fact and Proposition \ref{Proposition: pbw maps preserve module structure}, the first equality holds. Combining these three equalities, we obtain $\Theta^E(a \cdot s) = I(a) \cdot s + a \cdot \Theta^E(s)$. Considering the special case $E = \underline{\mathbb{C}}$, we also obtain $I(a \cdot s) = I(a) \cdot s + a \cdot I(s)$. Thus, \eqref{Equation: theta and I} holds, implying that our claim holds.
\end{proof}

In the remainder of this subsection, we show that the components of $F^E$ can be determined recursively from the initial term $\delta_Q^{-1}R^E \in \Omega^1(M, Q^* \otimes \operatorname{End} E)$. To this end, we begin by decomposing $I$ and $F^E$ as follows:
\begin{itemize}
	\item $I = \sum_{r=2}^\infty I_{(r)}$, where $I_{(r)} \in \Omega^1(M, \operatorname{Hom}(Q^*, \operatorname{Sym}^r Q^*))$ for each $r \geq 2$;
	\item $F^E = \sum_{r=1}^\infty F_{(r)}^E$, where $F_{(r)}^E \in \Omega^1(M, \operatorname{Sym}^r Q^* \otimes \operatorname{End} E)$ for each $r \geq 1$.
\end{itemize}
For the next proposition, we introduce the following notation for $r \geq 2$ and $r' \geq 1$:
\begin{equation}
	\iota_{I_{(r)}} F_{(r')}^E := \left( \tfrac{1}{r + r'} I_{(r), P}^l + \tfrac{1}{r+r'+1} I_{(r), Q}^l \right) \cdot \iota_{v_l} F_{(r')}^E.
\end{equation}
Here, we write $I_{(r)} = (I_{(r), P}^l + I_{(r), Q}^l) \otimes \partial_{u^l}$, where $I_{(r), P}^l$ and $I_{(r), Q}^l$ are local sections of $\bigwedge^1 P^* \otimes \operatorname{Sym}^r Q^*$ and $\bigwedge^1 Q^* \otimes \operatorname{Sym}^r Q^*$ respectively.

\begin{proposition}
	\label{Proposition 3.13}
	The bundle-valued $1$-form $F^E$ satisfies the recursive relations
	\begin{align*}
		F_{(1)}^E = & \delta_Q^{-1} R^E,\\
		F_{(r+1)}^E = & \delta_Q^{-1} \nabla F_{(r)}^E + \sum_{i=1}^{r-1} \iota_{I_{(i+1)}} F_{(r-i)}^E, \quad \text{for all } r \geq 1,
	\end{align*}
	where $\nabla = \nabla^{\widehat{\operatorname{Sym}} Q^* \otimes E^* \otimes E}$.
\end{proposition}

We first prove two lemmas required for the proof of the above proposition.

\begin{lemma}
	\label{Lemma 3.11}
	Let $Z \in \Gamma(M, Q)$, $A \in \Gamma(M, E^*)$, $X \in \Gamma(M, TM_\mathbb{C})$ and $s \in \Gamma(M, E)$. Then,
	\begin{equation*}
		\langle Z^r \otimes A, \iota_X F_{(r)}^E (s) \rangle = - \left\langle A \circ \nabla_X^E \circ \nabla_{Z^r}^E, S^E(s) \right\rangle, \quad \text{for } r \geq 1.
	\end{equation*}
	Here, $\iota_X F_{(r)}^E (s) \in \Gamma(M, \widehat{\operatorname{Sym}} Q^* \otimes E)$ denotes the contraction of the $\widehat{\operatorname{Sym}} Q^* \otimes E$-valued $1$-form $F_{(r)}^E(s)$ with the vector field $X$.
\end{lemma}
\begin{proof}
	Let $\phi = Z^r \otimes A$. By a degree argument, $\langle \phi, \iota_X F_{(r)}^E (s) \rangle = \langle \phi, \iota_X F^E (s) \rangle$, where
	\begin{equation*}
		\iota_X F^E(s) = \nabla_X^{\operatorname{K}, E} s - \nabla_X^E s + \iota_X (\delta_Qs) - \iota_X I(s) = \nabla_X^{\operatorname{K}, E} s - \nabla_X^E s,
	\end{equation*}
	since $\delta_Q s = I(s) = 0$. Thus it remains to compute $\langle \phi, \nabla_X^{\operatorname{K}, E} s \rangle$ and $\langle \phi, \nabla_X^E s \rangle$. The second term $\langle \phi, \nabla_X^E s \rangle$ vanishes for degree reasons. For the first term, we compute
	\begin{align*}
		\langle \phi, \nabla_X^{\operatorname{K}, E} s \rangle = & \mathcal{L}_X \langle \phi, s \rangle - \langle \mathcal{L}_X \circ A \circ \nabla_{Z^r}^E, S^E(s) \rangle\\
		= & \mathcal{L}_X \langle \phi, s \rangle - \left\langle A \circ \nabla_X^E \circ \nabla_{Z^r}^E, S^E(s) \right\rangle - \langle (\nabla_X^{E^*} A) \circ \nabla_{Z^r}^E, S^E(s) \rangle\\
		= & - \left\langle A \circ \nabla_X^E \circ \nabla_{Z^r}^E, S^E(s) \right\rangle.
	\end{align*}
	In the first line we used formula \eqref{Equation: pullback of Grothendieck connections}. The last equality follows again by degree considerations: $\langle \phi, s \rangle = 0$ and $\langle (\nabla_X^{E^*} A) \circ \nabla_{Z^r}^E, S^E(s) \rangle = \langle Z^r \otimes (\nabla_X^{E^*} A), s \rangle = 0$. This completes the proof.
\end{proof}

\begin{lemma}
	\label{Lemma 3.12}
	The following identity holds: $\delta_Q^{-1}F^E = 0$.
\end{lemma}
\begin{proof}
	Fix $r \geq 1$, $Z \in \Gamma(M, Q)$, $A \in \Gamma(M, E^*)$ and $s \in \Gamma(M, E)$. Then
	\begin{align*}
		\langle Z^{r+1} \otimes A, (\delta_Q^{-1} F_{(r)}^E)(s) \rangle = & \left\langle Z^r \otimes A, (\iota_Z F_{(r)}^E)(s) \right\rangle\\
		= & -\left\langle A \circ \nabla_Z^E \circ \nabla_{Z^r}^E, S^E(s) \right\rangle\\
		= & -\left\langle A \circ \nabla_{Z^{r+1}}^E, S^E(s) \right\rangle + \left\langle A \circ \nabla_{\nabla_Z^Q (Z^r)}^E, S^E(s) \right\rangle\\
		= & -\langle Z^{r+1} \otimes A, s \rangle + \langle (\nabla_Z^Q (Z^r)) \otimes A, s \rangle.
	\end{align*}
	In the second line, we used Lemma \ref{Lemma 3.11}. The whole term in the last line vanishes by degree reason. Therefore, $\delta_Q^{-1}F_{(r)}^E = 0$.
\end{proof}

We introduce a binary operation
\begin{equation*}
	\Omega^*(M, \widehat{\operatorname{Sym}} Q^* \otimes E^* \otimes E) \times \Omega^*(M, \widehat{\operatorname{Sym}} Q^* \otimes E^* \otimes E) \to \Omega^*(M, \widehat{\operatorname{Sym}} Q^* \otimes E^* \otimes E),
\end{equation*}
denoted by $(\Phi_1, \Phi_2) \mapsto \Phi_1 \circ \Phi_2$, and defined as follows. For $\alpha_1, \alpha_2 \in \Omega^*(M)$, $a_1, a_2 \in \Gamma(M, \widehat{\operatorname{Sym}} Q^*)$ and $A_1, A_2 \in \Gamma(M, E^* \otimes E)$,
\begin{equation*}
	(\alpha_1 \otimes a_1 \otimes A_1) \circ (\alpha_2 \otimes a_2 \otimes A_2) := \alpha_1 \wedge \alpha_2 \otimes (a_1 \cdot a_2) \otimes (A_1 \circ A_2).
\end{equation*}
The associated graded commutator is defined by
\begin{equation*}
	[\alpha, \beta] := \alpha \circ \beta - (-1)^{ij} \beta \circ \alpha,
\end{equation*}
for all $\alpha \in \Omega^i(M, \widehat{\operatorname{Sym}} Q^* \otimes E^* \otimes E)$ and $\beta \in \Omega^j(M, \widehat{\operatorname{Sym}} Q^* \otimes E^* \otimes E)$.\par
We now turn to the proof of the above proposition.

\begin{proof}[\myproof{Proposition}{\ref{Proposition 3.13}}]
	Since $\nabla^{\operatorname{K}, E}$ and $\nabla^{\operatorname{K}}$ are flat, we can deduce that
	\begin{equation}
		R^E + (\nabla - \delta_Q + I) (F^E) + \tfrac{1}{2} [F^E, F^E] = 0,
	\end{equation}
	recalling that $\nabla = \nabla^{\widehat{\operatorname{Sym}} Q^* \otimes E^* \otimes E}$. Rearranging terms, we obtain
	\begin{equation}
		\label{Equation 3.10}
		\delta_Q F^E = R^E + (\nabla + I) (F^E) + \tfrac{1}{2} [F^E, F^E].
	\end{equation}
	By Lemma \ref{Lemma 3.12}, $\delta_Q^{-1}F^E = 0$. As $F^E \in \Omega^1(M, \widehat{\operatorname{Sym}}^{\geq 1} Q^* \otimes E^* \otimes E)$, $\pi_Q(F^E) = 0$. By \eqref{Equation: Hodge decomposition for Dolbeault differential},
	\begin{equation*}
		F^E = \delta_Q^{-1} \delta_Q F^E + \delta_Q \delta_Q^{-1} F^E + \pi_Q(F^E) = \delta_Q^{-1} \delta_Q F^E.
	\end{equation*}
	Applying $\delta_Q^{-1}$ to both sides of \eqref{Equation 3.10} yields
	\begin{equation}
		\label{Equation 3.14}
		F^E = \delta_Q^{-1} \left( R^E + (\nabla + I) (F^E) + \tfrac{1}{2} [F^E, F^E] \right).
	\end{equation}
	We claim that $\delta_Q^{-1}[F^E, F^E] = 0$. To prove our claim, it is convenient to define another operator $\delta_Q^* := u^i \iota_{v_i}$ on $\Omega^*(M, \widehat{\operatorname{Sym}} Q^* \otimes E^* \otimes E)$. It is clear that for every $A \in \Omega^*(M, \widehat{\operatorname{Sym}} Q^* \otimes E^* \otimes E)$, $\delta_Q^* A = 0$ if and only if $\delta_Q^{-1} A = 0$. It is also easy to check that
	\begin{equation*}
		\delta_Q^* (F^E \circ F^E) = (\delta_Q^*F^E) \circ F^E - (F^E) \circ (\delta_Q^*F^E),
	\end{equation*}
	which vanishes due to Lemma \ref{Lemma 3.12}. Therefore, our claim holds.\par
	Now, decomposing both sides of \eqref{Equation 3.14} by polynomial degree in $\widehat{\operatorname{Sym}} Q^*$ yields a recursive formula: $F_{(1)}^E = \delta_Q^{-1} R^E$; for all $r \geq 1$,
	\begin{equation*}
		F_{(r+1)}^E = \delta_Q^{-1} \nabla F_{(r)}^E + \sum_{i=1}^{r-1} \delta_Q^{-1} \left( I_{(i+1)} (F_{(r-i)}^E) \right).
	\end{equation*}
	Write $I_{(r)} = I_{(r), P} + I_{(r), Q}$ with $I_{(r), P} = I_{(r), P}^l \otimes \partial_{u^l}$ and $I_{(r), Q} = I_{(r), Q}^l \otimes \partial_{u^l}$. Lemma 5.18 in \cite{LauStiXu2021} implies that $\delta_Q^*I_{(r), P} = \delta_Q^*I_{(r), Q} = 0$ (see Remark \ref{Remark: Kaparanov connection in trivial bundle case}). By direct computations,
	\begin{align*}
		\delta_Q^*(I_{(i+1), P}(F_{(r-i)}^E)) = & (\delta_Q^* I_{(i+1), P}) (F_{(r-i)}^E) - I_{(i+1), P} (\delta_Q^*F_{(r-i)}^E) + I_{(i+1), P}^l \cdot \iota_{v_l} (F_{(r-i)}^E)\\
		= & I_{(i+1), P}^l \cdot \iota_{v_l} (F_{(r-i)}^E),
	\end{align*}
	Similarly, $\delta_Q^*(I_{(i+1), P}(F_{(r-i)}^E)) = I_{(i+1), Q}^l \cdot \iota_{v_l} (F_{(r-i)}^E)$. Note that the $C^{0, 1}(M, \widehat{\operatorname{Sym}} Q^* \otimes E^* \otimes E)$-component of $F_{(r-i)}^E$ is annihilated by the contraction with the vector fields $v_l$. Eventually, we obtain $\delta_Q^{-1} (I_{(i+1)} (F_{(r-i)}^E)) = \iota_{I_{(i+1)}} F_{(r-i)}^E$, completing our proof.
\end{proof}

The following corollary follows from the recursive formula in Proposition \ref{Proposition 3.13} together with the properties of connections and curvatures of tensor products of vector bundles.

\begin{corollary}
	\label{Corollary: additivity of F}
	For $i = 1, 2$, let $E_i$ be a complex vector bundle over $M$ equipped with a connection $\nabla^{E_i}$ which restricts to a flat $P$-connection. Then
	\begin{equation*}
		F^{E_1 \otimes E_2} = F^{E_1} \otimes \operatorname{Id}_{E_2} + F^{E_2} \otimes \operatorname{Id}_{E_1}.
	\end{equation*}
\end{corollary}

\subsection{Chevalley--Eilenberg cohomology of $P$-modules}
\quad\par
\label{Subsection: Chevalley--Eilenberg cohomology}
The goal of this subsection is to state the following proposition.

\begin{proposition}
	\label{Proposition: qausi-isomorphism for polarized sections}
	The map
	\begin{equation*}
		\pi_Q: (\Omega^*(M, \widehat{\operatorname{Sym}} Q^* \otimes E), \nabla^{\operatorname{K}, E}) \to (C^{0, *}(M, E), d_P^E)
	\end{equation*}
	is a quasi-isomorphism. In particular, it induces an isomorphism between the sheaf of $\nabla^{\operatorname{K}, E}$-flat sections of $\widehat{\operatorname{Sym}} Q^* \otimes E$ and the sheaf of $d_P^E$-closed sections of $E$.
\end{proposition}

Since we will later require an analogous result (Lemma \ref{Lemma: quasi-isomorphism for quantizable functions}), whose proof proceeds along the same lines as Proposition \ref{Proposition: qausi-isomorphism for polarized sections}, we defer the detailed argument to Appendix \ref{Sectoin: a degeneracy condition for filtered complexes}.\par
For the moment, it suffices to recall from the proof of Proposition \ref{Proposition: qausi-isomorphism for polarized sections} that for every $s \in \Gamma(M, E)$, there is a unique section $J_s \in \Gamma(M, \widehat{\operatorname{Sym}} Q^* \otimes E)$ solving the equation
\begin{equation*}
	J_s - \delta_Q^{-1} ( \nabla^{\operatorname{K}, E} - d_P^{\widehat{\operatorname{Sym}} Q^* \otimes E} + \delta_Q ) J_s = s.
\end{equation*}
When $s$ is $d_P^E$-closed, $J_s$ is precisely the unique $\nabla^{\operatorname{K}, E}$-flat section of $\widehat{\operatorname{Sym}} Q^* \otimes E$ satisfying $\pi_Q(J_s) = s$. In this situation, the defining property of the Grothendieck connection $\nabla^{\operatorname{G}, E}$ implies that $S^E(J_s)$ is the jet of the $d_P^E$-closed section $s$.

\begin{corollary}
	Suppose $s \in \Gamma(M, E)$. Then, letting $\nabla = \nabla^{\widehat{\operatorname{Sym}} Q^* \otimes E}$, there is
	\begin{equation}
		\label{Equation: classical jet}
		J_s = \sum_{r=0}^\infty (\delta_Q^{-1} \nabla)^r s.
	\end{equation}
	In particular, for all $r \in \mathbb{N}$, $G \in \Gamma(M, \operatorname{Sym}^r Q)$ and $A \in \Gamma(M, E^*)$,
	\begin{equation}
		\label{Equation: pairing for classical jet}
		\langle G \otimes A, J_s \rangle = A \circ \nabla_G^E s.
	\end{equation}
\end{corollary}
\begin{proof}
	Since $d_P^{\widehat{\operatorname{Sym}} Q^* \otimes E} J_s \in C^{0, 1}(M, \widehat{\operatorname{Sym}} Q^* \otimes E)$, we have $\delta_Q^{-1} ( d_P^{\widehat{\operatorname{Sym}} Q^* \otimes E} ) J_s = 0$. On the other hand, $\nabla^{\operatorname{K}, E} + \delta_Q = \nabla + I + F^E$. Observe that $\delta_Q^{-1} (I(J_s)) = (\delta_Q^{-1} I)(J_s) = 0$ since $\delta_Q^{-1} I = 0$; and $\delta_Q^{-1} (F^E J_s) = (\delta_Q^{-1} F^E) J_s = 0$ since $\delta_Q^{-1}F^E = 0$. Therefore,
	\begin{equation*}
		J_s - \delta_Q^{-1} ( \nabla J_s ) = s.
	\end{equation*}
	The desired formula \eqref{Equation: classical jet} is then obtained by solving this equation recursively. Now, fix $r \in \mathbb{N}$, $Z \in \Gamma(M, Q)$ and $A \in \Gamma(M, E^*)$. Observe that $\langle Z^r \otimes A, J_s \rangle = \langle Z^r \otimes A, (\delta_Q^{-1} \nabla)^r s \rangle$. We can then prove by induction on $r$ that $\langle Z^r \otimes A, (\delta_Q^{-1} \nabla)^r s \rangle = A \circ \nabla_{Z^r}^E s$ (c.f. the proof of Proposition 3.16 in \cite{ChaLeuLiYau2025}). We are done.
\end{proof}

\section{Fedosov quantization in the presence of a polarization}
\label{Section: Fedosov quantization in the presence of a polarization}
Let $(M, \omega)$ be a symplectic manifold of real dimension $2n$ equipped with a polarization $P$. In this section we prepare for coupling Fedosov quantization with the polarization $P$ by constructing a curved dgla structure on the space $\Omega^*(M, \mathcal{W})$ for the \emph{Weyl bundle}
\begin{equation}
	\label{Equation: Weyl bundle}
	\mathcal{W} := \widehat{\operatorname{Sym}} T^*M_\mathbb{C}[[\hbar]].
\end{equation}
This is classical when one chooses a symplectic connection and equips $\mathcal{W}$ with the fibrewise Moyal product; Fedosov's quantization is then obtained by solving the weak Maurer--Cartan equation in this curved dgla \cite{Fed1994}.\par
To incorporate the polarization, we fix a (not necessarily involutive) complex Lagrangian subbundle $Q \subset TM_\mathbb{C}$ such that $TM_\mathbb{C} = P \oplus Q$. Such a complement always exists --- for example, $Q = J(\overline{P})$ for any $\omega$-compatible almost complex structure $J$ on $M$. With this setup, Subsection \ref{Subsection: polarized symplectic connection} relaxes the torsion-free assumption in Fedosov’s method, and Subsection \ref{Subsection: fibrwise star product} constructs the desired curved dgla via a fibrewise star product with separation of variables. Drawing on \cite{ChaLeuLi2023}, Subsection \ref{Subsection: fibrewise action} introduces a subalgebra bundle of $\mathcal{W}$ acting fibrewise on $\operatorname{Sym} Q^*$, yielding a corresponding curved dg-module structure. Finally, Subsection \ref{Subsection: Lagrangian splittings} introduces the formal Hodge decomposition on $\mathcal{W}$ adapted to the splitting $TM_\mathbb{C} = P \oplus Q$, completing the preparations for the next section, where Theorems \ref{First main result} and \ref{Second main result} are proved.

\subsection{Connections preserving symplectic forms and polarizations}
\quad\par
\label{Subsection: polarized symplectic connection}
Fedosov's original approach \cite{Fed1994} to the deformation quantization of $(M, \omega)$ requires the choice of a \emph{symplectic connection} $\nabla$, i.e. a torsion-free connection satisfying $\nabla \omega = 0$. To respect the decomposition $TM_\mathbb{C} = P \oplus Q$, one would ideally select a symplectic connection preserving both $P$ and $Q$. For instance, if $(M, \omega)$ is K\"ahler with $P = T^{0, 1}M$ and $Q = T^{1, 0}M$, the Levi--Civita connection fulfills these requirements.\par
In general, however, demanding such a symplectic connection is overly restrictive: its existence would force $Q$ to be involutive. More precisely, for all $Z, Z' \in \Gamma(M, Q)$,
\begin{equation*}
	[Z, Z'] = \nabla_ZZ' - \nabla_{Z'} Z \in \Gamma(M, Q).
\end{equation*}
As we will see, the torsion-free assumption on $\nabla$ is unnecessary for our purposes. It suffices to have a \emph{torsion-free} connection on $Q$ (see Remark \ref{Remark: torsion free connection on quotient bundles}), identifying $Q$ with the quotient $TM_\mathbb{C} / P$ via the decomposition $TM_\mathbb{C} = P \oplus Q$.

\begin{proposition}
	\label{Proposition: polarized symplectic connection}
	For any torsion-free connection $\nabla^Q$ on $Q$, there exists a unique connection $\nabla$ on $TM_\mathbb{C}$ which preserves $\omega$ and the subbundles $P$ and $Q$, and which extends $\nabla^Q$.
\end{proposition}
\begin{proof}
	We identify $P \cong Q^*$ via $Y \mapsto \iota_Y\omega$, and transfer the dual connection on $Q^*$ to $P$ to obtain a connection $\nabla^P$ on $P$ defined as follows: for all $Y \in \Gamma(M, P)$ and $Z \in \Gamma(M, Q)$,
	\begin{equation*}
		\omega(\nabla^PY, Z) = d(\omega(Y, Z)) - \omega(Y, \nabla^QZ).
	\end{equation*}
	Then $\nabla := \nabla^P \oplus \nabla^Q$ on $TM_\mathbb{C} = P \oplus Q$ satisfies the desired properties, giving existence.\par
	For uniqueness, if $\nabla'$ also satisfies the same properties, then $\nabla$ and $\nabla'$ agree on $Q$ and thus
	\begin{equation*}
		\omega(\nabla' Y, Z) = \omega(\nabla Y, Z), \quad \text{for all } Z \in \Gamma(M, Q) \text{ and } Y \in \Gamma(M, P),
	\end{equation*}
	as both $\nabla$ and $\nabla'$ preserve $\omega$. Non-degeneracy of $\omega$ implies $\nabla' Y = \nabla Y$, and hence $\nabla = \nabla'$.
\end{proof}

\begin{remark}
	Note that $\nabla$ is not required to preserve the real tangent bundle $TM$.
\end{remark}

From now on, we fix a torsion-free connection $\nabla^Q$ on $Q$ which, together with the splitting $TM_\mathbb{C} = P \oplus Q$, constitutes the auxiliary data \eqref{Equation: auxiliary data for deformation quantization}. Denote by $\nabla$ the induced connection on $TM_\mathbb{C}$ as in Proposition \ref{Proposition: polarized symplectic connection}. For later reference, we conclude this subsection by expressing the components of $\omega$, $\nabla$ and its curvature relative to a given pair of local frames $(\check{v}_1, ..., \check{v}_n)$ of $P$ and $(v_1, ..., v_n)$ of $Q$. For all $1 \leq i, j \leq n$, define the local functions
\begin{equation}
	\omega_{i\check{j}} = \omega(v_i, \check{v}_j).
\end{equation}
Let $(\omega^{\check{i}j})$ be the inverse matrix of $(\omega_{i\check{j}})$. As $\nabla$ preserves both $P$ and $Q$, we can write
\begin{equation}
	\nabla v_i = \Gamma_i^k v_k \quad \text{and} \quad \nabla \check{v}_j = \Gamma_{\check{j}}^{\check{l}} \check{v}_l,
\end{equation}
where $\Gamma_i^k$'s and $\Gamma_{\check{j}}^{\check{l}}$ are local $1$-forms on $M$. Since $\nabla \omega = 0$, we can verify that
\begin{equation}
	\label{Equation: coordinate description of symplectic form preserving property}
	d \omega_{i\check{j}} = \Gamma_i^k \omega_{k\check{j}} + \Gamma_{\check{j}}^{\check{l}} \omega_{i\check{l}} \quad \text{and} \quad d \omega^{\check{j}i} = -( \omega^{\check{j}k} \Gamma_k^i + \omega^{\check{l}i} \Gamma_{\check{l}}^{\check{j}} ).
\end{equation}
Also, we can write
\begin{equation}
	\nabla^2(v_i) = R_i^k v_k \quad \text{and} \quad \nabla^2(\check{v}_j) = R_{\check{j}}^{\check{l}} \check{v}_l,
\end{equation}
where $R_i^k$'s and $R_{\check{j}}^{\check{l}}$'s are local $2$-forms on $M$. Since $\nabla \omega = 0$ again,
\begin{equation}
	\label{Equation: Curvature and symplectic form}
	\omega_{k\check{j}} R_{i}^k + \omega_{i\check{l}} R_{\check{j}}^{\check{l}} = 0.
\end{equation}

\subsection{Curved dgla structure via a polarization-adapted fibrewise star product}
\quad\par
\label{Subsection: fibrwise star product}
In this subsection, we will construct a curved dgla on $\Omega^*(M, \mathcal{W})$. We first notice that, using a choice of local frames $(\check{v}_1, ..., \check{v}_n)$ of $P$ and $(v_1, ..., v_n)$ of $Q$, every form on $M$ with values in $\mathcal{W} = \widehat{\operatorname{Sym}} T^*M_\mathbb{C}[[\hbar]]$ can be locally expressed as a formal sum of terms of the form
\begin{equation*}
	\hbar^r f v^{b_1} \wedge \cdots \wedge v^{b_q} \wedge \check{v}^{a_1} \wedge \cdots \wedge \check{v}^{a_p} \otimes u^{j_1} \cdots u^{j_m} \check{u}^{i_1} \cdots \check{u}^{i_l},
\end{equation*}
where $f$ is a local smooth $\mathbb{C}$-valued function, $v^i, \check{v}^j, u^i, \check{u}^j$ follows the conventions specified in Subsection \ref{Subsection: notations}, $r \in \mathbb{N}$ and $1 \leq i_1, ..., i_l, j_1, ..., j_m, a_1, ..., a_p, b_1, ..., b_q \leq n$ are index labels.\par
We then define a fibrewise star product $\star^{\operatorname{F}}$ with separation of variables on $\mathcal{W}$.

\begin{definition}
	For sections $a, b$ of $\mathcal{W}$, set
	\begin{equation}
		\label{Equation: star product}
		a \star^{\operatorname{F}} b = \sum_{r=0}^\infty \frac{\hbar^r}{r!} \omega^{\check{j}_1i_1} \cdots \omega^{\check{j}_ri_r} \frac{\partial^r a}{\partial \check{u}^{j_1} \cdots \partial \check{u}^{j_r}} \frac{\partial^r b}{\partial u^{i_1} \cdots \partial u^{i_r}}.
	\end{equation}
\end{definition}
This is well-defined independently of the choice of local frames for $P$ and $Q$, so $(\mathcal{W}, \star^{\operatorname{F}})$ is a non-commutative algebra bundle over $M$, called the \emph{Weyl bundle}. We also extend $\star^{\operatorname{F}}$ to $\Omega^*(M, \mathcal{W})$ by graded linearity, and denote by $[\, \cdot\,, \, \cdot\,]_{\star^{\operatorname{F}}}$ its graded commutator.\par
The connection $\nabla$ on $TM_\mathbb{C}$ induces a connection on $\mathcal{W}$, which is denoted by the same symbol by abuse of notation. We define $R_\omega = -\omega_{k\check{j}} R_{i}^k u^i \check{u}^j \in \Omega^2(M, \operatorname{Sym}^2 T^*M_\mathbb{C})$, which, by \eqref{Equation: Curvature and symplectic form}, can also be written as $R_\omega = \omega_{i\check{l}} R_{\check{j}}^{\check{l}} u^i \check{u}^j$. We now state the following proposition, whose proof relies crucially on the fact that $\nabla$ preserves $\omega$, $P$ and $Q$.

\begin{proposition}
	\label{Proposition 3.6}
	The following quadruple forms a curved differential graded Lie algebra:
	\begin{equation*}
		(\Omega^*(M, \mathcal{W}), R_\omega, \nabla, \tfrac{1}{\hbar}[\, \cdot\,, \, \cdot\,]_{\star^{\operatorname{F}}}).
	\end{equation*}
\end{proposition}
\begin{proof}
	We first claim that for all $a, b \in \Gamma(M, \mathcal{W})$ and $X \in \Gamma(M, TM)$,
	\begin{equation*}
		\nabla_X(a \star^{\operatorname{F}} b) = (\nabla_X a) \star^{\operatorname{F}} b + a \star^{\operatorname{F}} (\nabla_X b).
	\end{equation*}
	Applying the Leibniz rule for $\nabla_X$ to Formula \eqref{Equation: star product} yields terms in which $X$ acts on the structural coefficients $\omega^{\check{j}i}$ and on the fibrewise derivatives of $a$ and $b$. By a direct computation,
	\begin{align*}
		\left[\nabla_X , \frac{\partial}{\partial \check{u}^j} \right] (a)  = \Gamma_{\check{j}}^{\check{l}}(X) \frac{\partial a}{\partial \check{u}^l} \quad \text{and} \quad \left[ \nabla_X, \frac{\partial}{\partial u^i} \right] (b) = \Gamma_i^k(X) \frac{\partial b}{\partial u^k}.
	\end{align*}
	Using these relations to compare $\nabla_X(a \star^{\operatorname{F}} b)$ with $(\nabla_X a) \star^{\operatorname{F}} b + a \star^{\operatorname{F}} (\nabla_X b)$, one finds
	\begin{equation*}
		\nabla_X(a \star^{\operatorname{F}} b) - (\nabla_X a) \star^{\operatorname{F}} b - a \star^{\operatorname{F}} (\nabla_X b) = \sum_{r=0}^\infty \frac{\hbar^r}{r!} \Omega^{i_1 \cdots i_r \check{j}_1 \cdots \check{j}_r} \frac{\partial^r a}{\partial \check{u}^{j_1} \cdots \partial \check{u}^{j_r}} \frac{\partial^r b}{\partial u^{i_1} \cdots \partial u^{i_r}},
	\end{equation*}
	where
	\begin{equation*}
		\Omega^{i_1 \cdots i_r \check{j}_1 \cdots \check{j}_r} := \mathcal{L}_X(\omega^{\check{j}_1i_1} \cdots \omega^{\check{j}_ri_r} )+\sum_{m=1}^r \omega^{\check{j}_1i_1} \cdots \widehat{\omega^{\check{j}_mi_m}} \cdots \omega^{\check{j}_ri_r} (\omega^{\check{l}i_m} \Gamma_{\check{l}}^{\check{j}_m}(X)+\omega^{\check{j}_mk} \Gamma_k^{i_m}(X) ).
	\end{equation*}
	By \eqref{Equation: coordinate description of symplectic form preserving property}: $d \omega^{\check{j}i} = -( \omega^{\check{j}k} \Gamma_k^i + \omega^{\check{l}i} \Gamma_{\check{l}}^{\check{j}} )$. Hence, all $\Omega^{i_1 \cdots i_r \check{j}_1 \cdots \check{j}_r}$ vanish, concluding our claim.\par
	Since $[a, b]_{\star^{\operatorname{F}}} \in \hbar \cdot \Omega^*(M, \mathcal{W})$ for all $a, b \in \Omega^*(M, \mathcal{W})$, $\tfrac{1}{\hbar}[\, \cdot\,, \, \cdot\,]_{\star^{\operatorname{F}}}$ is a well-defined graded Lie bracket on $\Omega^*(M, \mathcal{W})$. The above claim implies that $\nabla$ is compatible with $\tfrac{1}{\hbar}[\, \cdot\,, \, \cdot\,]_{\star^{\operatorname{F}}}$. It remains to show that the following equality holds on $\Omega^*(M, \mathcal{W})$:
	\begin{equation*}
		\nabla^2 = \frac{1}{\hbar} [R_\omega, \, \cdot\,]_{\star^{\operatorname{F}}}.
	\end{equation*}
	This equality is shown by a direct computation as follows:
	\begin{align*}
		\frac{1}{\hbar}[ R_\omega, \, \cdot\,]_{\star^{\operatorname{F}}} = & \omega^{\check{j}i} \left( \frac{\partial R_\omega}{\partial \check{u}^j} \frac{\partial}{\partial u^i} - \frac{\partial R_\omega}{\partial u^i} \frac{\partial}{\partial \check{u}^j} \right) = - \omega^{\check{j}i} \left( \omega_{a\check{j}} R_{k}^a u^k \frac{\partial}{\partial u^i} + \omega_{i\check{b}} R_{\check{l}}^{\check{b}} \check{u}^l \frac{\partial}{\partial \check{u}^j} \right)\\
		= & - R_{k}^i u^k \frac{\partial}{\partial u^i} - R_{\check{l}}^{\check{j}} \check{u}^l \frac{\partial}{\partial \check{u}^j}.
	\end{align*}
	The term in the last line is exactly $\nabla^2$ on $\Omega^*(M, \mathcal{W})$.
\end{proof}

In accordance with Fedosov's method, the next step in the construction of a deformation quantization of $(M, \omega)$ is to find a $1$-form  $\gamma \in \Omega^1(M, \mathcal{W})$ satisfying the weak Maurer–Cartan equation for the deformation of $\nabla$. We defer this construction to the next section.

\subsection{An action of the polarization-adapted fibrewise star product}
\quad\par
\label{Subsection: fibrewise action}
At each point $x \in M$, the fibre of $\mathcal{W}$ can be viewed as a deformation quantization of the symplectic vector space $T_xM$, whereas the fibre of $\widehat{\operatorname{Sym}} Q^*$ can be interpreted as the geometric quantization of $T_xM$ with respect to the polarization $P_x$. This naturally suggests a fibrewise $\mathcal{W}$-module structure on $\widehat{\operatorname{Sym}} Q^*$, given on generators by
\begin{equation*}
	u^i \mapsto u^i \, \cdot\,, \quad \check{u}^j \mapsto \hbar \omega^{\check{j}i} \tfrac{\partial}{\partial u^i},
\end{equation*}
after evaluating at a nonzero value of $\hbar$. However, for a general section of $\mathcal{W}$, convergence issues prevent a straightforward evaluation.\par
To address this, we introduce an increasing filtration $\mathcal{W}^{\leq 0} \subset \cdots \subset \mathcal{W}^{\leq r} \subset \cdots$ by
\begin{equation*}
	\mathcal{W}^{\leq r} := \bigoplus_{i + j \leq r} \widehat{\operatorname{Sym}} Q^* \otimes \operatorname{Sym}^i P^* \cdot \hbar^j,
\end{equation*}
and define the \emph{finite polarized weight subbundle} of $\mathcal{W}$ as the inductive limit of this filtration:
\begin{equation*}
	\mathcal{W}^{<\infty} := \widehat{\operatorname{Sym}} Q^* \otimes \operatorname{Sym} P^*[\hbar].
\end{equation*}
It follows from \eqref{Equation: star product} that $(\mathcal{W}, \star^{\operatorname{F}})$ induces a fibrewise filtered algebra structure on $\mathcal{W}^{<\infty}$. Consequently, one can define a fibrewise action of $\mathcal{W}^{<\infty}$ on $\widehat{\operatorname{Sym}} Q^*$ for each evaluation $\hbar = \tfrac{\sqrt{-1}}{k}$, where $k \neq 0$. In this paper, we focus on $k \in \mathbb{Z}^+$, which can be interpreted as the tensor power of a prequantum line bundle.

\begin{definition}
	Define a $\mathcal{C}^\infty(M)$-bilinear map
	\begin{equation*}
		\Gamma(M, \mathcal{W}^{<\infty}) \times \Gamma(M, \widehat{\operatorname{Sym}} Q^*) \to \Gamma(M, \widehat{\operatorname{Sym}} Q^*), \quad (a, s) \mapsto a \circledast_k^{\operatorname{F}} s,
	\end{equation*}
	as follows: for $a = \hbar^r u^{i_1} \cdots u^{i_p} \check{u}^{j_1} \cdots \check{u}^{j_q}$,
	\begin{equation}
		\label{Equation: fibrewise action}
		a \circledast_k^{\operatorname{F}} s = \left( \frac{\sqrt{-1}}{k} \right)^{r+q} \omega^{\check{j}_1k_1} \cdots \omega^{\check{j}_qk_q} u^{i_1} \cdots u^{i_p} \frac{\partial^q s}{\partial u^{k_1} \cdots \partial u^{k_q}}.
	\end{equation}
\end{definition}

This map equips $\widehat{\operatorname{Sym}} Q^*$ with a $\mathcal{W}^{<\infty}$-module structure:
\begin{equation}
	\label{Equation: fibrewise star product and fibrewise action}
	(a \star^{\operatorname{F}} b) \circledast_k^{\operatorname{F}} s = a \circledast_k^{\operatorname{F}} (b \circledast_k^{\operatorname{F}} s),
\end{equation}
which extends to $\Omega^*(M, \widehat{\operatorname{Sym}} Q^*)$ by graded linearity.\par

\begin{remark}
	In the K\"ahler case (with opposite ordering), Chan--Leung--Li refer to $\circledast_k^{\operatorname{F}}$ as the \emph{fibrewise Bargmann--Fock action} \cite{ChaLeuLi2023}, due to its resemblance to the standard Bargmann--Fock action on holomorphic functions on $\mathbb{C}^n$.
\end{remark}

\begin{proposition}
	\label{Lemma: fibrewise action and connection}
	For each $k \in \mathbb{Z}^+$, $(\Omega^*(M, \widehat{\operatorname{Sym}} Q^*), \circledast_k^{\operatorname{F}})$ is a curved differential graded module over the curved differential graded algebra $(\Omega^*(M, \mathcal{W}^{<\infty}), R_\omega, \nabla, \star^{\operatorname{F}})$.
\end{proposition}
\begin{proof}
	First, for all $X \in \Gamma(M, TM_\mathbb{C})$, $a \in \Gamma(M, \mathcal{W}^{<\infty})$ and $s \in \Gamma(M, \widehat{\operatorname{Sym}} Q^*)$, we have
	\begin{equation*}
		\nabla_X (a \circledast_k^{\operatorname{F}} s) = (\nabla_X a) \circledast_k^{\operatorname{F}} s + a \circledast_k^{\operatorname{F}} (\nabla_X s),
	\end{equation*}
	which follows by the same type of computations as in the proof of Proposition \ref{Proposition 3.6}. It remains to check the curvature condition:
	\begin{equation*}
		\nabla^2 = \tfrac{k}{\sqrt{-1}} R_\omega \circledast_k^{\operatorname{F}} \quad \text{on } \Omega^*(M, \widehat{\operatorname{Sym}} Q^*),
	\end{equation*}
	using the fact that $\tfrac{1}{\hbar}[R_\omega, \, \cdot\,]_{\star^{\operatorname{F}}}$ is independent of $\hbar$ and Lemma \ref{Lemma: commutator of star product and action}.
\end{proof}

\begin{lemma}
	\label{Lemma: commutator of star product and action}
	For $\alpha \in \Omega^*(M, \widehat{\operatorname{Sym}} Q^* \otimes P^*)$, $\left. \tfrac{1}{\hbar} [\alpha, \, \cdot\,]_{\star^{\operatorname{F}}} \right\vert_{\hbar = \frac{\sqrt{-1}}{k}} = \tfrac{k}{\sqrt{-1}} \alpha \circledast_k^{\operatorname{F}}$ on $\Omega^*(M, \widehat{\operatorname{Sym}} Q^*)$.
\end{lemma}
\begin{proof}
	Let $\beta \in \Omega^*(M, \widehat{\operatorname{Sym}} Q^*)$. Then $\left. \tfrac{1}{\hbar} [\alpha, \beta]_{\star^{\operatorname{F}}} \right\vert_{\hbar = \frac{\sqrt{-1}}{k}} = \omega^{\check{j}i} \tfrac{\partial \alpha}{\partial \check{u}^j} \tfrac{\partial \beta}{\partial u^i} = \tfrac{k}{\sqrt{-1}} \alpha \circledast_k^{\operatorname{F}}$.
\end{proof}

\subsection{Polarization-adapted formal Hodge decomposition on the Weyl bundle}
\quad\par
\label{Subsection: Lagrangian splittings}
Recall from Fedosov's quantization that the Weyl bundle $\mathcal{W}$ admits a formal Hodge decomposition: there exist $\mathcal{C}^\infty(M)[[\hbar]]$-linear operators $\delta, \delta^{-1}, \pi_0 \colon \Omega^*(M, \mathcal{W}) \to \Omega^*(M, \mathcal{W})$ satisfying the standard Hodge-type identity:
\begin{equation}
	\label{Equation: formal Hodge decomposition}
	\operatorname{Id} - \pi_0 = \delta \circ \delta^{-1} + \delta^{-1} \circ \delta.
\end{equation}
In local frames adapted to the splitting \(TM_\mathbb{C} = P \oplus Q\), we have $\delta = v^i \wedge \partial_{u^i} + \check{v}^i \wedge \partial_{\check{u}^i}$, and for $\alpha \in \Omega^k(M, \operatorname{Sym}^r T^*M_\mathbb{C})$ with $k + r > 0$,
\begin{equation*}
	\delta^{-1} \alpha = \tfrac{1}{k+r} \big(u^i \cdot \iota_{v_i} + \check{u}^i \cdot \iota_{\check{v}_i}\big) \alpha, \quad
	\pi_0(\alpha) = 0,
\end{equation*}
while for $\alpha \in \mathcal{C}^\infty(M)$,
$\delta^{-1}\alpha = 0$ and $\pi_0(\alpha) = \alpha$.\par
The splitting \(TM_\mathbb{C} = P \oplus Q\) allows us to refine the formal Hodge decomposition. We define six $\mathcal{C}^\infty(M)[[\hbar]]$-linear operators on $\Omega^*(M, \mathcal{W})$ adapted to the polarization:
\begin{equation*}
	\delta_Q, \delta_Q^{-1}, \pi_Q, \quad \delta_P, \delta_P^{-1}, \pi_P
\end{equation*}
as follows. First, $\delta_Q := v^i \wedge \partial_{u^i}$ and $\delta_P := \check{v}^i \wedge \partial_{\check{u}^i}$, giving a decomposition $\delta = \delta_P + \delta_Q$. Then for $\alpha \in C^{q,p}(M, \operatorname{Sym}^l Q^* \otimes \operatorname{Sym}^m P^*)$, we define the homotopy and projection operators:
\begin{align*}
	\delta_Q^{-1} \alpha &:= 
	\begin{cases} 
		\frac{1}{q+l} u^i \cdot \iota_{v_i} \alpha, & q+l>0,\\
		0, & q+l=0,
	\end{cases}
	&
	\pi_Q(\alpha) &:= 
	\begin{cases}
		0, & q+l>0,\\
		\alpha, & q+l=0,
	\end{cases} \\
	\delta_P^{-1} \alpha &:= 
	\begin{cases} 
		\frac{1}{p+m} \check{u}^i \cdot \iota_{\check{v}_i} \alpha, & p+m>0,\\
		0, & p+m=0,
	\end{cases}
	&
	\pi_P(\alpha) &:= 
	\begin{cases}
		0, & p+m>0,\\
		\alpha, & p+m=0.
	\end{cases}
\end{align*}
These operators satisfy Hodge-type identities individually:
\begin{align*}
	\operatorname{Id} - \pi_Q &= \delta_Q \circ \delta_Q^{-1} + \delta_Q^{-1} \circ \delta_Q,\\
	\operatorname{Id} - \pi_P & = \delta_P \circ \delta_P^{-1} + \delta_P^{-1} \circ \delta_P.
\end{align*}
These formal Hodge decompositions are essential for the constructions in the next section.

\section{A construction of star products via their action on polarized sections}
\label{Section: construction of star products via action on polarized sections}
Continuing from Section \ref{Section: Fedosov quantization in the presence of a polarization}, we now consider the following setup:
\begin{maintheorem}
	\label{Setup: two}
	$(M, \omega)$, $(L, \nabla^L)$ and $P$ are as in Setup \ref{Setup: one}, and $\mathbf{L}$ is a complex line bundle over $M$ with a connection $\nabla^{\mathbf{L}}$ extending a flat $P$-connection. Set $L^{(k)} := L^{\otimes k} \otimes \mathbf{L}$ for $k \in \mathbb{Z}^+$.
\end{maintheorem}
In Subsection \ref{Subsection: solution to Fedosov equation}, we construct a solution to Fedosov's equation using the following $1$-forms, which arise from the decomposition of Grothendieck connections via pbw maps:
\begin{enumerate}
	\item $I \in \Omega^1(M, \operatorname{Hom}(Q^*, \widehat{\operatorname{Sym} Q^*}))$ given as in \eqref{Equation A.6};
	\item $F^L \in \Omega^1(M, \widehat{\operatorname{Sym}} Q^* \otimes \operatorname{End} (L)))$ and $F^\mathbf{L} \in \Omega^1(M, \widehat{\operatorname{Sym}} Q^* \otimes \operatorname{End} (\mathbf{L}))$ given as in \eqref{Equation A.7}. 
\end{enumerate}
Subsection \ref{Subsection: Star products} follows a standard procedure to verify that this solution induces a star product on $\mathcal{C}^\infty(M)[[\hbar]]$. Subsection \ref{Subsection: Quantizable functions with respect to a polarization} completes the proof of Theorem \ref{First main result} by introducing \emph{formal quantizable functions}, a generalization of the concept of Chan--Leung--Li from the Kähler polarization case \cite{ChaLeuLi2023} to the general setting. These functions form an algebra that naturally acts on the sheaf of $P$-polarized sections of $L^{(k)}$. Subsection \ref{Subsection: non formal quantizable functions} turns to the non-formal analogue and proves Theorem \ref{Second main result}. Finally, Subsection \ref{Subsection: examples} discusses some examples.

\subsection{Solutions to Fedosov's equation induced by transverse pbw maps}
\label{Subsection: solution to Fedosov equation}
\quad\par
In this subsection, we prove the following proposition. To this end, we extend $\circledast_k^{\operatorname{F}}$ to an action on $\widehat{\operatorname{Sym}} Q^* \otimes L^{(k)}$ by letting it act trivially on the $L^{(k)}$ factor.

\begin{proposition}
	\label{Theorem: flat connection in Fedosov quantization}
	The Fedosov equation
	\begin{equation}
		\label{Equation: Fedosov equation}
		R_\omega + \nabla \gamma + \frac{1}{2\hbar} \left[ \gamma, \gamma \right]_{\star^{\operatorname{F}}} = -\omega - \hbar R^\mathbf{L}
	\end{equation}
	admits a solution $\gamma \in \Omega^1(M, \mathcal{W}^{<\infty})$ such that the induced flat connection
	\begin{equation}
		\label{Equation: flat connection in Fedosov quantization}
		D := \nabla + \tfrac{1}{\hbar} \left[ \gamma, \, \cdot\, \right]_{\star^{\operatorname{F}}}
	\end{equation}
	on $\mathcal{W}$ satisfies the following conditions:
	\begin{enumerate}
		\item for all $a \in \Gamma(M, \mathcal{W}^{<\infty})$, $Da \in \Omega^1(M, \mathcal{W}^{<\infty})$;
		\item for all $a \in \Gamma(M, \mathcal{W}^{<\infty})$, $k \in \mathbb{Z}^+$ and $s \in \Gamma(M, \widehat{\operatorname{Sym}} Q^* \otimes L^{(k)})$,
		\begin{equation*}
			\nabla^{\operatorname{K}, L^{(k)}} (a \circledast_k^{\operatorname{F}} s) = (Da) \circledast_k^{\operatorname{F}} s + a \circledast_k^{\operatorname{F}} \nabla^{\operatorname{K}, L^{(k)}} s,
		\end{equation*}
	\end{enumerate}
	where $\nabla^{\operatorname{K}, L^{(k)}}$ denotes the flat connection on $\widehat{\operatorname{Sym}} Q^* \otimes L^{(k)}$ characterized by \eqref{Equation: pullback of Grothendieck connections}.
\end{proposition}

Our construction of $\gamma$ does not follow from Fedosov's original recursive argument for solving \eqref{Equation: Fedosov equation}. Instead, motivated by the observation of Chan–Leung–Li \cite{ChaLeuLi2022b} (see also \cite{Yau2025}), we aim to reinterpret the connection $\nabla^{\operatorname{K}, L^{(k)}}$ in terms of the fibrewise action $\circledast_k^{\operatorname{F}}$:
\begin{equation*}
	\nabla^{\operatorname{K}, L^{(k)}} = \nabla^{\widehat{\operatorname{Sym}} Q^* \otimes L^{(k)}} + \tfrac{k}{\sqrt{-1}} \gamma \circledast_k^{\operatorname{F}},
\end{equation*}
where a single $1$-form $\gamma \in \Omega^1(M, \mathcal{W})$ is determined so that this equality holds for all $k \in \mathbb{Z}^+$. Once $\gamma$ is constructed, we verify that it satisfies \eqref{Equation: Fedosov equation} by invoking the flatness of $\nabla^{\operatorname{K}, L^{(k)}}$. This approach provides a convenient and conceptual way to obtain $\gamma$, taking advantage of the decomposition of $\nabla^{\operatorname{K}, L^{(k)}}$ established in \eqref{Equation A.8}：
\begin{equation*}
	\nabla^{\operatorname{K}, L^{(k)}} = \nabla^{\widehat{\operatorname{Sym}} Q^* \otimes L^{(k)}} - \delta_Q + I + F^{L^{(k)}}.
\end{equation*}
In other words, our aim is to solve the equation: $\tfrac{k}{\sqrt{-1}} \gamma \circledast_k^{\operatorname{F}} = - \delta_Q + I + F^{L^{(k)}}$ for all $k \in \mathbb{Z}^+$. The following lemma shows that the above equation has at most one solution.

\begin{lemma}
	\label{Lemma: Uniqueness of solution}
	Let $a \in \Gamma(M, \mathcal{W}^{<\infty})$. Suppose that for every $k \in \mathbb{Z}^+$,
	\begin{equation*}
		a \circledast_k^{\operatorname{F}} = 0
	\end{equation*}
	as an operator on $\Gamma(M, \widehat{\operatorname{Sym}} Q^* \otimes L^{(k)})$. Then $a = 0$.
\end{lemma}
\begin{proof}
	Locally write $a = \sum_{l=0}^r \sum_{\lvert \check{J} \rvert \leq r} (-\sqrt{-1}\hbar)^l a_{\check{J}, l} \check{u}^J$ for some $r \in \mathbb{N}$. Here, $\check{J}$ ranges over multi-indices with $\lvert \check{J} \rvert \leq r$, and each $a_{\check{J}, l}$ is a local section of $\widehat{\operatorname{Sym}} Q^*$. Choose local frames $e$ of $L$ and $\mathbf{e}$ of $\mathbf{L}$. The assumption
	\begin{equation*}
		a \circledast_k^{\operatorname{F}} (e^{\otimes k} \otimes \mathbf{e}) = 0 \quad \text{for all } k \in \mathbb{Z}^+
	\end{equation*}
	gives the relation $k^r a_{\emptyset, 0} + k^{r-1} a_{\emptyset, 1} + \cdots + a_{\emptyset, r} = 0$ for all $k \in \mathbb{Z}^+$, hence $a_{\emptyset, 0} = \cdots = a_{\emptyset, r} = 0$. Now fix $1 \leq p < r$ and assume inductively that $a_{\check{J}, 0} = \cdots = a_{\check{J}, r} = 0$ for all $\check{J}$ with $\lvert \check{J} \rvert < p$. Using this hypothesis and the nondegeneracy of $\omega$, the condition
	\begin{equation*}
		a \circledast_k^{\operatorname{F}} ( u^J \otimes e^{\otimes k} \otimes \mathbf{e}) = 0 \quad \text{for all } k \in \mathbb{Z}^+ \text{ and } J \text{ with } \lvert J \rvert = p
	\end{equation*}
	forces that $k^r a_{\check{J}, 0} + k^{r-1} a_{\check{J}, 1} + \cdots + a_{\check{J}, r} = 0$ for all $k \in \mathbb{Z}^+$ and $\check{J}$ with $\lvert \check{J} \rvert = p$. Hence, $a_{\check{J}, 0} = \cdots = a_{\check{J}, r} = 0$. The conclusion $a = 0$ follows by induction on $p$.
\end{proof}

To proceed, we first solve the equation: $\tfrac{k}{\sqrt{-1}} \gamma \circledast_k^{\operatorname{F}} = -\delta_Q$ for all $k \in \mathbb{Z}^+$.

\begin{lemma}
	For all $k \in \mathbb{Z}^+$, $\tfrac{k}{\sqrt{-1}} ( \delta_P^{-1}\omega ) \circledast_k^{\operatorname{F}} = -\delta_Q$.
\end{lemma}
\begin{proof}
	Observing that $\delta_P^{-1} \omega = \delta_P^{-1} ( \omega_{i\check{j}} v^i \wedge \check{v}^j ) = -\omega_{i\check{j}} v^i \otimes \check{u}^j$, we obtain
	\begin{equation*}
		\tfrac{k}{\sqrt{-1}} ( \delta_P^{-1}\omega ) \circledast_k^{\operatorname{F}} = -\tfrac{k}{\sqrt{-1}} \omega_{i\check{j}} v^i \wedge \tfrac{\sqrt{-1}}{k} \omega^{\check{j}i'} \tfrac{\partial}{\partial u^{i'}} = -v^i \wedge \tfrac{\partial}{\partial u^i} = -\delta_Q.
	\end{equation*}
\end{proof}

Next, via the isomorphism $(Q^*)^* \cong Q \to P^*$ induced by $\omega$, we obtain a $1$-form from $I$:
\begin{equation}
	\widetilde{I} \in \Omega^1(M, \widehat{\operatorname{Sym}} Q^* \otimes P^*) \subset \Omega^1(M, \mathcal{W}^{<\infty}).
\end{equation}
By construction, $\gamma = \widetilde{I}$ is a solution to the equation: $\tfrac{k}{\sqrt{-1}} \gamma \circledast_k^{\operatorname{F}} = I$ for all $k \in \mathbb{Z}^+$.\par
Finally, we regard $F^L$ and $F^\mathbf{L}$ as elements of $\Omega^1(M, \widehat{\operatorname{Sym}} Q^*) \subset \Omega^1(M, \mathcal{W}^{<\infty})$ using the canonical trivializations of $\operatorname{End} (L)$ and $\operatorname{End}(\mathbf{L})$. By Corollary \ref{Corollary: additivity of F} and the following lemma, $\gamma = \delta_Q^{-1} \omega + \hbar F^\mathbf{L}$ is a solution to the equation: $\tfrac{k}{\sqrt{-1}} \gamma\circledast_k^{\operatorname{F}} = F^{L^{(k)}}$ for all $k \in \mathbb{Z}^+$.

\begin{lemma}
	\label{Lemma 5.2}
	$F^L = \tfrac{1}{\sqrt{-1}} \delta_Q^{-1} \omega$.
\end{lemma}
\begin{proof}
	By Proposition \ref{Proposition 3.13}, we obtain $F_{(1)}^L = \tfrac{1}{\sqrt{-1}} \delta_Q^{-1} \omega = \tfrac{1}{\sqrt{-1}} \omega_{i\check{j}} \check{v}^j \otimes u^i \in C^{0, 1}(M, Q^*)$ and
	\begin{equation*}
		F_{(2)}^L = \tfrac{1}{\sqrt{-1}} \delta_Q^{-1} \nabla^{\widehat{\operatorname{Sym}} Q^*} \delta_Q^{-1} \omega = \tfrac{1}{\sqrt{-1}} \delta_Q^{-1} \left( \left( d\omega_{i\check{j}} - \Gamma_i^l \omega_{l\check{j}} \right) \wedge \check{v}^j \otimes u^i \right) = 0.
	\end{equation*}
	In the above equality, we used \eqref{Equation: coordinate description of symplectic form preserving property}. Now, assume that $F_{(2)}^L = \cdots = F_{(r)}^L = 0$ for some $r \geq 2$. By Proposition \ref{Proposition 3.13} again, $F_{(r+1)}^L = \iota_{I_{(r)}} F_{(1)}^L = \tfrac{1}{\sqrt{-1}} \iota_{I_{(r)}} \delta_Q^{-1} \omega$, which vanishes since $\iota_{I_{(r)}}$ annihilates $C^{0, 1}(M, Q^*)$. The proof is now complete by induction.
\end{proof}

Combining all the above observations, we have proved the following proposition.

\begin{proposition}
	\label{Proposition: alternative desription of Grothendieck connection}
	There exists a unique $\mathcal{W}^{<\infty}$-valued $1$-form $\gamma$ such that for all $k \in \mathbb{Z}^+$,
	\begin{align*}
		\nabla^{\operatorname{K}, L^{(k)}} = & \nabla^{\widehat{\operatorname{Sym}} Q^* \otimes L^{(k)}} + \tfrac{k}{\sqrt{-1}} \gamma \circledast_k^{\operatorname{F}}.
	\end{align*}
	Explicitly, $\gamma$ is given by
	\begin{equation}
		\label{Equation: twisted solution to Fedosov equation}
		\gamma := \delta_P^{-1}\omega + \delta_Q^{-1}\omega + \widetilde{I} + \hbar F^\mathbf{L} \in \Omega^1(M, \mathcal{W}^{\leq 1}).
	\end{equation}
\end{proposition}

\begin{proof}[\myproof{Proposition}{\ref{Theorem: flat connection in Fedosov quantization}}]
	Define $\gamma$ by \eqref{Equation: twisted solution to Fedosov equation}. Fix $k \in \mathbb{Z}^+$. We have several observations:
	\begin{itemize}
		\item By Proposition \ref{Lemma: fibrewise action and connection}, $\tfrac{k}{\sqrt{-1}} \left( R_\omega + \omega + \hbar R^\mathbf{L} \right) \circledast_k^{\operatorname{F}} = \left( \nabla^{\widehat{\operatorname{Sym}} Q^* \otimes L^{(k)}} \right)^2$.
		\item By the same proposition, $\tfrac{k}{\sqrt{-1}} (\nabla \gamma) \circledast_k^{\operatorname{F}} = \left[ \nabla^{\widehat{\operatorname{Sym}} Q^* \otimes L^{(k)}}, \tfrac{k}{\sqrt{-1}} \gamma \circledast_k^{\operatorname{F}} \right]$.
		\item By \eqref{Equation: fibrewise star product and fibrewise action}, $\tfrac{k}{\sqrt{-1}} \left( \tfrac{1}{2\hbar} [\gamma, \gamma]_{\star^{\operatorname{F}}} \right) \circledast_k^{\operatorname{F}} = \tfrac{1}{2} \left[ \tfrac{k}{\sqrt{-1}} \gamma \circledast_k^{\operatorname{F}}, \tfrac{k}{\sqrt{-1}} \gamma \circledast_k^{\operatorname{F}}\right]$.
	\end{itemize}
	Based on these observations, Proposition \ref{Proposition: alternative desription of Grothendieck connection} and the fact that $(\nabla^{\operatorname{K}, L^{(k)}})^2 = 0$, we obtain
	\begin{equation}
		\label{Equation 3.9}
		\tfrac{k}{\sqrt{-1}} \left( R_\omega + \nabla \gamma + \tfrac{1}{2\hbar} [\gamma, \gamma]_{\star^{\operatorname{F}}} + \omega + \hbar R^\mathbf{L} \right) \circledast_k^{\operatorname{F}} = 0.
	\end{equation}
	We have $R_\omega + \nabla \gamma + \tfrac{1}{2\hbar} [\gamma, \gamma]_{\star^{\operatorname{F}}} + \omega + \hbar R^\mathbf{L} \in \Omega^2(M, \mathcal{W}^{<\infty})$, since $\gamma \in \Omega^1(M, \mathcal{W}^{<\infty})$. Then \eqref{Equation: Fedosov equation} holds according to Lemma \ref{Lemma: Uniqueness of solution}. Note that $-\omega - \hbar R^\mathbf{L} \in \Omega^2(M, \mathbb{C})[\hbar]$ is central in $\tfrac{1}{\hbar}[\, \cdot\,, \, \cdot\,]_{\star^{\operatorname{F}}}$. Thus, $D^2 = \tfrac{1}{\hbar}[-\omega - \hbar R^\mathbf{L}, \, \cdot\,]_{\star^{\operatorname{F}}} = 0$, where $D$ is given as in \eqref{Equation: flat connection in Fedosov quantization}.\par
	Now, fix $a \in \Gamma(M, \mathcal{W}^{<\infty})$ and $s \in \Gamma(M, \widehat{\operatorname{Sym}} Q^* \otimes L^{(k)})$. It is clear that $\nabla a \in \Gamma(M, \mathcal{W}^{<\infty})$. As $\gamma \in \Omega^1(M, \mathcal{W}^{\leq 1})$, $\tfrac{1}{\hbar}[\gamma, a]_{\star^{\operatorname{F}}} \in \Omega^1(M, \mathcal{W}^{<\infty})$. Thus, $D a \in \Omega^1(M, \mathcal{W}^{<\infty})$. For the second condition, due to Proposition \ref{Proposition: alternative desription of Grothendieck connection}, It suffices to show that
	\begin{equation*}
		\gamma \circledast_k^{\operatorname{F}} (a \circledast_k^{\operatorname{F}} s) = [\gamma, a]_{\star_k^{\operatorname{F}}} \circledast_k^{\operatorname{F}} s + a \circledast_k^{\operatorname{F}} (\gamma \circledast_k^{\operatorname{F}} s),
	\end{equation*}
	which is straightforward.
\end{proof}

\subsection{Star products induced by transverse pbw maps}
\quad\par
\label{Subsection: Star products}
In this subsection, we show that the solution $\gamma$ from \eqref{Equation: twisted solution to Fedosov equation} determines a star product $\star$ on $\mathcal{C}^\infty(M)[[\hbar]]$, based on the statement that for every $f \in \mathcal{C}^\infty(M)[[\hbar]]$, there is a unique $D$-flat section $\mathbf{J}_f$ of $\mathcal{W}$ such that $\pi_0(\mathbf{J}_f) = f$. Here, $\pi_0: \Omega^*(M, \mathcal{W}) \to \mathcal{C}^\infty(M)[[\hbar]]$ is the projection map introduced in Subsection \ref{Subsection: Lagrangian splittings}. This statement follows from the following proposition.

\begin{proposition}
	\label{Proposition: qausi-isomorphism for formal functions}
	The following projection is a quasi-isomorphism:
	\begin{equation*}
		\pi_0: (\Omega^*(M, \mathcal{W}), D) \to (\mathcal{C}^\infty(M)[[\hbar]], 0).
	\end{equation*}
	Hence, it induces an isomorphism between the sheaf of $D$-flat sections of $\mathcal{W}$ and $\mathcal{C}_M^\infty[[\hbar]]$.
\end{proposition}

To prove this proposition, it is convenient to rewrite the differential $D$ in an alternative form. Direct computations yield $\tfrac{1}{\hbar} [\delta_P^{-1} \omega, \, \cdot\,]_{\star^{\operatorname{F}}} = -\delta_Q$ and $\tfrac{1}{\hbar} [\delta_Q^{-1} \omega, \, \cdot\,]_{\star^{\operatorname{F}}} = -\delta_P$. Therefore,
\begin{equation*}
	D = \nabla - \delta + \frac{1}{\hbar} \left[ \widetilde{I} + \hbar F^\mathbf{L}, \, \cdot\, \right]_{\star^{\operatorname{F}}},
\end{equation*}
where $\delta = \delta_Q + \delta_P$ is as in Subsection \ref{Subsection: Lagrangian splittings}. We also introduce a \emph{double-weight decomposition}:
\begin{equation}
	\label{Equation: double weight decomposition}
	a = \sum_{r, l \in \mathbb{N}} (a)_{r, l}, \quad (a)_{r, l} \in \Gamma\left( M, \textstyle \bigoplus_{i + j = r} \operatorname{Sym}^l Q^* \otimes \operatorname{Sym}^i P^* \cdot \hbar^j \right),
\end{equation}
for any section $a \in \Gamma(M, \mathcal{W})$, where $r$ is the polarized weight and $l$ is the polynomial degree in $Q^*$; we also call the sum $r + l$ the \emph{total weight}. The components of $D$ act on the double weight $(r, l)$ as follows: $\nabla$ preserves $(r, l)$; $\delta_Q$ sends $(r, l)$ to $(r, l-1)$; $\delta_P$ sends $(r, l)$ to $(r-1, l)$; the commutator term preserves $r$ and does not decrease $l$.

\begin{lemma}
	\label{Lemma: qausi-isomorphism for formal functions}
	There is a cochain isomorphism
	\begin{equation*}
		\Psi: (\Omega^*(M, \mathcal{W}), D) \to (\Omega^*(M, \mathcal{W}), -\delta)
	\end{equation*}
	defined by $\Psi(a) = a - \delta^{-1} (D + \delta) a$.
\end{lemma}
\begin{proof}
	To verify that $\Psi$ is a cochain map, observe that $D^2 = \delta^2 = 0$. Then
	\begin{align*}
		\Psi D + \delta \Psi = & (D + \delta) - \delta^{-1} (D + \delta) D - \delta\delta^{-1} (D + \delta)\\
		= & (D + \delta) - \delta^{-1} \delta (D + \delta) - \delta\delta^{-1} (D + \delta)\\
		= & \pi_0 (D + \delta) = 0.
	\end{align*}
	Next, we claim that, for any $a' \in \Omega^*(M, \mathcal{W})$, there exists a unique $a \in \Omega^*(M, \mathcal{W})$ solving
	\begin{equation*}
		a - \delta^{-1} (D + \delta) a = a'.
	\end{equation*}
	This follows because $\delta^{-1} (D + \delta)$ always increases the total weight. Hence, the solution can be constructed iteratively by successively solving for each total weight component, ensuring both existence and uniqueness. Therefore, $\Psi$ is an isomorphism.
\end{proof}

\begin{proof}[\myproof{Proposition}{\ref{Proposition: qausi-isomorphism for formal functions}}]
	By \eqref{Equation: formal Hodge decomposition}, we know that $\pi_0: (\Omega^*(M, \mathcal{W}), -\delta) \to (\mathcal{C}^\infty(M)[[\hbar]], 0)$ is a quasi-isomorphism. By Lemma \ref{Lemma: qausi-isomorphism for formal functions} and the fact that $\pi_0 \circ \delta^{-1} = 0$,
	\begin{equation*}
		\pi_0 \circ \Psi = \pi_0: (\Omega^*(M, \mathcal{W}), D) \to (\mathcal{C}^\infty(M)[[\hbar]], 0)
	\end{equation*}
	is a quasi-isomorphism.
\end{proof}

Now, we define $\mathbf{J}_f := \Psi^{-1}(f)$, where $\Psi$ is given as in Lemma \ref{Lemma: qausi-isomorphism for formal functions}, and
\begin{equation}
	\label{Equation: formula of star product}
	f \star g := \pi_0 \left( \mathbf{J}_f \star^{\operatorname{F}} \mathbf{J}_g \right).
\end{equation}


\begin{theorem}
	\label{Theorem: star product}
	Under Steup \ref{Setup: two}, $(\mathcal{C}^\infty(M)[[\hbar]], \star)$ is a deformation quantization of $(M, \omega)$.
\end{theorem}
\begin{proof}
	It is evident from \eqref{Equation: formula of star product} that there is a sequence of bi-differential operators $C_0, C_1, ...$ such that for all $f, g \in \mathcal{C}^\infty(M)[[\hbar]]$, $f \star g = \sum_{r=0}^\infty \hbar^r C_r(f, g)$. In particular,
	\begin{align*}
		C_0(f, g) = & (\mathbf{J}_f)_{0, 0} \cdot (\mathbf{J}_g)_{0, 0},\\
		C_1(f, g) = & \mathbf{J}_{f, 2} \cdot (\mathbf{J}_g)_{0, 0} + \omega^{\check{j}i} \frac{\partial \mathbf{J}_{f, 1}}{\partial \check{u}^j} \frac{\partial \mathbf{J}_{g, 1}}{\partial u^i} + (\mathbf{J}_f)_{0, 0} \cdot \mathbf{J}_{g, 2}.
	\end{align*}
	Here, $\mathbf{J}_{f, 1}, \hbar \cdot \mathbf{J}_{f, 2}$ are the $\Gamma(M, T^*M_\mathbb{C})$- and $\hbar \cdot \mathcal{C}^\infty(M)$-components of $\mathbf{J}_f$ respectively, and $\mathbf{J}_{g, 1}, \hbar \cdot \mathbf{J}_{g, 2}$ are their $\mathbf{J}_g$-counterparts. Note that $\mathbf{J}_f$ solves the equation
	\begin{equation}
		\label{Equation: defining equation of quantum jets}
		\mathbf{J}_f - \delta^{-1} (D + \delta) \mathbf{J}_f = f.
	\end{equation}
	Comparing components on both sides, we obtain $(\mathbf{J}_f)_{0, 0} = f$, $\mathbf{J}_{f, 1} = \delta^{-1}df$ and $\mathbf{J}_{f, 2} = 0$ and similar formulas for $\mathbf{J}_{g, i}$'s. Therefore, $C_0(f, g) = f \cdot g$ and
	\begin{align*}
		C_1(f, g) - C_1(g, f) = & \omega^{\check{j}i} \left( \frac{\partial (\delta^{-1} df)}{\partial \check{u}^j} \frac{\partial (\delta^{-1} dg)}{\partial u^i} - \frac{\partial (\delta^{-1} dg)}{\partial \check{u}^j} \frac{\partial (\delta^{-1} df)}{\partial u^i} \right)\\
		= & \omega^{\check{j}i} ( (\mathcal{L}_{\check{v}_j} f) (\mathcal{L}_{v_i} g) - (\mathcal{L}_{\check{v}_j} g) (\mathcal{L}_{v_i} f) ) = \{f, g\},
	\end{align*}
	where $\{\, \cdot\,, \, \cdot\,\}$ is the Poisson bracket induced by $\omega$. We can also clearly see that, for the constant function $1$ on $M$, $\mathbf{J}_1 = 1$, whence $1 \star f = f \star 1 = 1$. The associativity is also immediate.
\end{proof}

\subsection{Formal quantizable functions with respect to a polarization}
\quad\par
\label{Subsection: Quantizable functions with respect to a polarization}
In this subsection, we complete our proof of Theorem \ref{First main result}.\par
Before presenting the argument, we outline the guiding idea. Recall that a $P$\emph{-polarized} section of $L^{(k)}$ is a smooth section $s$ satisfying $\nabla_Y^{L^{(k)}} s = 0$ for all $Y \in \Gamma(M, P)$, where $\nabla^{L^{(k)}}$ is the connection induced by $\nabla^L$ and $\nabla^{\mathbf{L}}$. As $\star^{\operatorname{F}}$ and $\circledast_k^{\operatorname{F}}$ are compatible with the flat connections $D$ and $\nabla^{\operatorname{K}, L^{(k)}}$, it is natural to attempt to descend these operations to flat sections, thereby defining an action of the deformed algebra $(\mathcal{C}^\infty(M)[[\hbar]], \star)$ on the space $\Gamma_P(M, L^{(k)})$ of $P$-polarized sections of $L^{(k)}$. However, such a construction runs into convergence problems when arbitrary formal functions are allowed. To address this, we restrict to a subalgebra for which the construction stabilizes at finite order, ensuring that the resulting action is well defined.

\begin{definition}
	\label{Definition: formal quantizable functions}
	A \emph{formal quantizable function} is a formal function $f \in \mathcal{C}^\infty(M)[\hbar]$ such that $\mathbf{J}_f \in \Gamma(M, \mathcal{W}^{<\infty})$. It is said to be $r$\emph{th order} if $\mathbf{J}_f \in \Gamma(M, \mathcal{W}^{\leq r})$.
\end{definition}

The associated sheaf of formal quantizable functions, denoted $\mathcal{C}_{M, \hbar}^{(<\infty)}$, is a subsheaf of algebras of $(\mathcal{C}_M^\infty[[\hbar]], \star)$. It carries the natural increasing filtration $\mathcal{C}_{M, \hbar}^{(0)} \subset \mathcal{C}_{M, \hbar}^{(1)} \subset \cdots \subset \mathcal{C}_{M, \hbar}^{(r)} \subset \cdots$, where $\mathcal{C}_{M, \hbar}^{(r)}$ is the sheaf of order-$r$ formal quantizable functions on $M$.\par
To describe the induced action on polarized sections, we first introduce:

\begin{definition}
	Let $a \in \Gamma(M, \mathcal{W}^{<\infty})$. Define $a \circledast_k: \Gamma(M, L^{(k)}) \to \Gamma(M, L^{(k)})$ by
	\begin{equation}
		a \circledast_k s := \pi_Q \left( a \circledast_k^{\operatorname{F}} J_s \right), \quad \text{for all } s \in \Gamma(M, L^{(k)}).
	\end{equation}
\end{definition}

\begin{remark}
	As $\pi_Q$ picks out the $\Gamma(M, L^{(k)})$-component of an element of $\Gamma(M, \widehat{\operatorname{Sym}} Q^* \otimes L^{(k)})$, the operator $a \circledast_k$ depends only on the $\Gamma(M, \operatorname{Sym} P^*)[\hbar]$-component of $a$.
\end{remark}

\begin{lemma}
	\label{Lemma: module structure}
	Let $k \in \mathbb{Z}^+$ and $a, b \in \Gamma(M, \mathcal{W}^{<\infty})$ be $D$-closed. Then for all $s \in \Gamma_P(M, L^{(k)})$, $b \circledast_k s \in \Gamma_P(M, L^{(k)})$ and
	\begin{equation*}
		a \circledast_k (b \circledast_k s) = ( a \star^{\operatorname{F}} b ) \circledast_k s.
	\end{equation*}
\end{lemma}
\begin{proof}
	By definition, $\pi_Q(b \circledast_k^{\operatorname{F}} J_s) = b \circledast_k s$. Since $Db = \nabla^{\operatorname{K}, L^{(k)}} J_s = 0$,
	\begin{equation*}
		\nabla^{\operatorname{K}, L^{(k)}} (b \circledast_k^{\operatorname{F}} J_s) = (D b) \circledast_k^{\operatorname{F}} J_s + b \circledast_k^{\operatorname{F}} (\nabla^{\operatorname{K}, L^{(k)}} J_s) = 0.
	\end{equation*}
	By the uniqueness of $J_{b \circledast_k s}$, we obtain $J_{b \circledast_k s} = b \circledast_k^{\operatorname{F}} J_s$ and $d_P^{L^{(k)}} (b \circledast_k s) = 0$. Therefore,
	\begin{equation*}
		a \circledast_k (b \circledast_k s) = \pi_Q (a \circledast_k^{\operatorname{F}} (b \circledast_k^{\operatorname{F}} J_s )) = \pi_Q(( a \star^{\operatorname{F}} b ) \circledast_k^{\operatorname{F}} J_s) = ( a \star^{\operatorname{F}} b ) \circledast_k s.
	\end{equation*}
\end{proof}

We thus obtain an action of $\mathcal{C}_{M, \hbar}^{(<\infty)}(M)$ on $\Gamma_P(M, L^{(k)})$:
\begin{equation}
	\label{Equation: action of formal quantizable functions}
	\mathcal{C}_{M, \hbar}^{(<\infty)}(M) \to \operatorname{End}_\mathbb{C} \Gamma_P(M, L^{(k)}), \quad f \mapsto \mathbf{J}_f \circledast_k,
\end{equation}
and since the construction is local, it extends to an action at the level of sheaves.\par
Finally, to prove our first main result, we focus on the first two filtered pieces of $\mathcal{C}_{M, \hbar}^{(<\infty)}$, namely $\mathcal{C}_{M, \hbar}^{(0)}$ and $\mathcal{C}_{M, \hbar}^{(1)}$. Recall the double weight decomposition $\mathbf{J}_f = \sum_{r, l \in \mathbb{N}} (\mathbf{J}_f)_{r, l}$ in \eqref{Equation: double weight decomposition}.

\begin{lemma}
	\label{Lemma: characterization of formal quantizable functions}
	Let $f = \sum_{i=0}^\infty \hbar^i f_i \in \mathcal{C}^\infty(M)[[\hbar]]$ with each $f_i \in \mathcal{C}^\infty(M)$, and let $r \in \mathbb{N}$. Then $f$ is formal quantizable of order $r \in \mathbb{N}$ if and only if $(\mathbf{J}_f)_{r+1, 0} = 0$ and $f_m = 0$ for all $m > r$.
\end{lemma}
\begin{proof}
	By \eqref{Equation: defining equation of quantum jets}, for all $m \in \mathbb{N}$,
	\begin{equation}
		\label{Equation: recursive formula for double weight components}
		(\mathbf{J}_f)_{m+1, 0} = \hbar^{m+1} \cdot f_{m+1} + \delta_P^{-1} \nabla (\mathbf{J}_f)_{m, 0} - \delta_P^{-1} \pi_Q \left( (\mathbf{J}_f)_{m, 0} \star^{\operatorname{F}} \left( \frac{1}{\hbar} \widetilde{I} + F^\mathbf{L} \right) \right).
	\end{equation}
	Assume $f \in \mathcal{C}_{M, \hbar}^{(r)}(M)$. Fix $m \geq r$. By assumption, $(\mathbf{J}_f)_{m+1, 0} = 0$ (in particular $(\mathbf{J}_f)_{r+1, 0} = 0$). Since the operator $\delta_P^{-1}$ strictly increases the $P^*$-degree, it follows from the R.H.S. of \eqref{Equation: recursive formula for double weight components} that the $\hbar^{m+1} \cdot \mathcal{C}^\infty(M)$-component of $(\mathbf{J}_f)_{m+1, 0}$ is precisely $\hbar^{m+1} \cdot f_{m+1}$, whence $f_{m+1} = 0$.\par
	Conversely, assume $(\mathbf{J}_f)_{r+1, 0} = 0$ and $f_m = 0$ for all $m > r$. We claim, by induction on $l \geq 0$, that $(\mathbf{J}_f)_{m, l} = 0$ for all $m > r$. The case $l = 0$ is justified by our assumption and \eqref{Equation: recursive formula for double weight components}.\par
	Assume our claim holds for all degrees $\leq l$. We inspect the $(m, l)$-component of the flatness equation $D\mathbf{J}_f = 0$ for $m > r$. Among all contributions, the only term involving $(\mathbf{J}_f)_{m, l+1}$ is $-\delta_Q(\mathbf{J}_f)_{m, l+1}$. All remaining terms involve $(\mathbf{J}_f)_{m', l'}$ with $l' \leq l$ and $m' \geq m$, and hence vanish by the induction hypothesis. Therefore, $\delta_Q (\mathbf{J}_f)_{m, l+1} = 0$. Since $\delta_Q$ is injective on the summand of positive $Q^*$-degree, this forces $(\mathbf{J}_f)_{m, l+1} = 0$, completing the induction. Consequently, $\mathbf{J}_f \in \Gamma(M, \mathcal{W}^{\leq r})$. Hence $f$ is formal quantizable of order $r$.
\end{proof}

\begin{theorem}
	\label{Theorem: zeroth order quantizable functions}
	Under Setup \ref{Setup: two}, we have $\mathcal{C}_{M, \hbar}^{(0)}(M) = \mathcal{O}(M)$. Furthermore,
	\begin{enumerate}
		\item for all $f, g \in \mathcal{O}(M)$, $f \star g = f \cdot g$;
		\item for all $f \in \mathcal{O}(M)$, $k \in \mathbb{Z}^+$ and $s \in \Gamma_P(M, L^{(k)})$, $\mathbf{J}_f \circledast_k s = f \cdot s$.
	\end{enumerate}
\end{theorem}
\begin{proof}
	By Lemma \ref{Lemma: characterization of formal quantizable functions}, $f \in \mathcal{C}_{M, \hbar}^{(0)}(M)$ if and only if $f \in \mathcal{C}^\infty(M)$ and $\delta_P^{-1} df = (\mathbf{J}_f)_{1, 0} = 0$. This is precisely the condition that $f \in \mathcal{O}(M)$. Let $f, g \in \mathcal{O}(M)$, $k \in \mathbb{Z}^+$ and $s \in \Gamma_P(M, L^{(k)})$. Then $\mathbf{J}_f, \mathbf{J}_g \in \Gamma(M, \widehat{\operatorname{Sym}} Q^*)$, and we obtain
	\begin{align*}
		f \star g = & \pi_0(\mathbf{J}_f \star^{\operatorname{F}} \mathbf{J}_g) = \pi_0(\mathbf{J}_f \cdot \mathbf{J}_f) = f \cdot g,\\
		\mathbf{J}_f \circledast_k s = & \pi_Q(\mathbf{J}_f \circledast_k J_s) = \pi_Q(\mathbf{J}_f \cdot J_s) = f \cdot s,
	\end{align*}
	using the separation-of-variable property of $\star^{\operatorname{F}}$ (see \eqref{Equation: star product}) and Formula \eqref{Equation: fibrewise action} for $\circledast_k^{\operatorname{F}}$.
\end{proof}

\begin{remark}
	The associated $\nabla^{\operatorname{K}}$-flat section $J_f$ of $\widehat{\operatorname{Sym}} Q^*$ satisfies $\pi_0(J_f) = \pi_Q(J_f) = f$ and
	\begin{align*}
		D J_f = & \nabla^{\widehat{\operatorname{Sym}} Q^*} J_f - \delta_Q f + I(J_f) = \nabla^{\operatorname{K}} J_f = 0,
	\end{align*}
	where we used $\delta_P J_f = \tfrac{1}{\hbar}[\hbar F^\mathbf{L}, J_f]_{\star^{\operatorname{F}}} = 0$. Hence, uniqueness implies $\mathbf{J}_f = J_f$. 
\end{remark}

Recall that $f \in \mathcal{C}^\infty(M)$ is said to \emph{preserve the polarization} $P$ if $[X_f, Y] \in \Gamma(M, P)$ for all $Y \in \Gamma(M, P)$, where $X_f$ is the Hamiltonian vector field associated with $f$.

\begin{theorem}
	\label{Theorem: first order quantizable functions}
	Under Setup \ref{Setup: two}, a formal function $f \in \mathcal{C}^\infty(M)[[\hbar]]$ is formal quantizable of order $1$ if and only if it can be written as $f = f_0 + \hbar f_1$, where $f_0 \in \mathcal{C}^\infty(M)$ preserves $P$ and $f_1 \in \mathcal{C}^\infty(M)$ satisfies $( df_1 - \iota_{X_{f_0}} R^\mathbf{L} ) \vert_P = 0$. In this case, for all $k \in \mathbb{Z}^+$ and $s \in \Gamma_P(M, L^{\otimes k})$,
	\begin{equation*}
		\mathbf{J}_f \circledast_k s = f_0 \cdot s + \tfrac{\sqrt{-1}}{k} f_1 \cdot s + \tfrac{\sqrt{-1}}{k} \nabla_{X_{f_0}}^{L^{(k)}} s.
	\end{equation*}
\end{theorem}
\begin{proof}
	Suppose firs that $f \in \mathcal{C}_{M, \hbar}^{(1)}(M)$. By Lemma \ref{Lemma: characterization of formal quantizable functions}, we must have $f = f_0 + \hbar f_1$ with $f_0, f_1 \in \mathcal{C}^\infty(M)$. Using the recursive formula \eqref{Equation: recursive formula for double weight components}, we compute
	\begin{align*}
		(\mathbf{J}_f)_{0, 0} = & f_0,\\
		(\mathbf{J}_f)_{1, 0} = & \hbar \cdot f_1 + \delta_P^{-1} df_0,\\
		(\mathbf{J}_f)_{2, 0} = & \delta_P^{-1} \nabla (\delta_P^{-1} df_0) + \hbar \cdot \delta_P^{-1} \left( df_1 - \iota_{X_{f_0}} R^\mathbf{L} \right).
	\end{align*}
	Here the Hamiltonian vector field $X_{f_0}$ is decomposed as $X_{f_0} = Y_{f_0} + Z_{f_0}$ with $Y_{f_0} \in \Gamma(M, P)$ and $Z_{f_0} \in \Gamma(M, Q)$. Since $\nabla^\mathbf{L}$ extends a flat $P$-connection, the term $\delta_P^{-1} \iota_{Y_{f_0}} R^\mathbf{L}$ vanishes automatically. The condition $(\mathbf{J}_f)_{2, 0} = 0$ is therefore equivalent to
	\begin{equation*}
		\delta_P^{-1} \nabla (\delta_P^{-1} df_0) = \left. \left( df_1 - \iota_{X_{f_0}} R^\mathbf{L} \right) \right\vert_P = 0.
	\end{equation*}
	Using $\nabla \omega = 0$, one verifies that $\delta_P^{-1} \nabla (\delta_P^{-1} df_0) = 0 \Leftrightarrow d_P^Q (Z_{f_0}) = 0$, which in turn is equivalent to the condition that $f_0$ preserves the polarization $P$ (by torsion-freeness of induced connection on $Q$). Thus the ``only if'' direction is proved, and the ``if'' direction follows from the same computation together with Lemma \ref{Lemma: characterization of formal quantizable functions}. Finally, let $k \in \mathbb{Z}^+$ and $s \in \Gamma_P(M, L^{(k)})$, we obtain
	\begin{equation*}
		\mathbf{J}_{f_0} \circledast_k s = \pi_Q( (f_0 + \hbar \cdot f_1 + \delta_P^{-1} df_0) \circledast_k^{\operatorname{F}} J_s) = f_0 \cdot s + \tfrac{\sqrt{-1}}{k} f_1 \cdot s + \tfrac{\sqrt{-1}}{k} \nabla_{X_{f_0}}^{L^{(k)}} s,
	\end{equation*}
	where we used $\nabla_{Y_{f_0}}^{L^{(k)}} s = 0$ since $s$ is $P$-polarized.
\end{proof}

In conclusion, Theorems \ref{Theorem: star product}, \ref{Theorem: zeroth order quantizable functions} and \ref{Theorem: first order quantizable functions} together imply Theorem \ref{First main result}.

\subsection{Non-formal quantizable functions with respect to a polarization}
\quad\par
\label{Subsection: non formal quantizable functions}
In this subsection we prove Theorem \ref{Second main result}. To do so, we first introduce, for each level $k \in \mathbb{Z}^+$, a non-formal analogue of Definition \ref{Definition: formal quantizable functions}. Let $D_k$ denote the evaluation of $D$ at $\hbar = \tfrac{\sqrt{-1}}{k}$, and let $\star_k^{\operatorname{F}}$ denote the corresponding specialization of $\star^{\operatorname{F}}$ on $\mathcal{W}^{<\infty}$. Unlike in the formal case, a $D_k$-flat section of $\widehat{\operatorname{Sym}} Q^* \otimes \operatorname{Sym} P^*$ is no longer uniquely determined by its $\mathcal{C}^\infty(M)$-component (see \cite{ChaLeuLi2023}). Following \cite{ChaLeuLi2023}, we therefore make the following definition.

\begin{definition}
	A \emph{level-}$k$ \emph{quantizable function} on $M$ is a $D_k$-flat section
	\begin{equation*}
		a \in \Gamma(M, \widehat{\operatorname{Sym}} Q^* \otimes \operatorname{Sym} P^*).
	\end{equation*}
	It is said to be $r$\emph{th order} if $a \in \Gamma(M, \widehat{\operatorname{Sym}} Q^* \otimes \operatorname{Sym}^{\leq r} P^*)$.
\end{definition}

The resulting sheaf $\mathcal{O}_k^{(<\infty)}$ of level-$k$ quantizable functions carries an increasing filtration $\mathcal{O}_k^{(0)} \subset \mathcal{O}_k^{(1)} \subset \cdots \mathcal{O}_k^{(r)} \subset \cdots$, where $\mathcal{O}_k^{(r)}$ consists of all order-$r$ level-$k$ quantizable functions. Evaluation at $\hbar = \tfrac{\sqrt{-1}}{k}$ defines a morphism of sheaves of $\mathbb{C}$-algebras:
\begin{equation*}
	(\mathcal{C}_{M, \hbar}^{(<\infty)}, \star) \to (\mathcal{O}_k^{(<\infty)}, \star_k^{\operatorname{F}}), \quad f \mapsto \left. \mathbf{J}_f \right\vert_{\hbar = \frac{\sqrt{-1}}{k}},
\end{equation*}
whose kernel is the ideal sheaf $(\hbar  - \tfrac{\sqrt{-1}}{k}) \cdot \mathcal{C}_{M, \hbar}^{(<\infty)}$. From the proof of Lemma \ref{Lemma: module structure}, we also know that $\circledast_k$ defines an action of $(\mathcal{O}^{(k)}, \star_k^{\operatorname{F}})$ on $P$-polarized sections of $L^{(k)}$. Let $\widetilde{\mathcal{D}}_{L^{(k)}}$ be the sheaf of transverse differential operators on $L^{(k)}$, i.e. $d_P^{\widetilde{D}(L^{(k)}, L^{(k)})}$-closed sections of $\widetilde{D}(L^{(k)}, L^{(k)})$. It is a sheaf of algebras with multiplication $\circ$ given by \eqref{Equation: composition of transverse differential operators}. We now present (a refinement of) our second main result (Theorem \ref{Second main result}).

\begin{theorem}
	\label{Theorem: isomorphism of transverse differential operators}
	Under Setup \ref{Setup: two}, for each $k \in \mathbb{Z}^+$, the induced action $\circledast_k$ of $\mathcal{O}_k^{(<\infty)}$ on $P$-polarized sections of $L^{(k)}$ defines an isomorphism of sheaves of filtered algebras
	\begin{equation*}
		(\mathcal{O}_k^{(<\infty)}, \star_k^{\operatorname{F}}) \to (\widetilde{\mathcal{D}}_{L^{(k)}}, \circ).
	\end{equation*}
\end{theorem}

The guiding idea is to regard the finite-polarized-weight Weyl bundle as
encoding a prolongation of transverse differential operators. Indeed, by
Theorem \ref{Theorem: pbw theorem}, an order-$r$ transverse differential operator
is represented by a section of $\operatorname{Sym}^{\leq r}Q$.
Prolonging its coefficients in the transverse directions naturally leads to $\widehat{\operatorname{Sym}}Q^*
\otimes \operatorname{Sym}^{\leq r}Q$. At level $k$, the symplectic form identifies
$Q$ with $P^*$ via $Z \mapsto \tfrac{k}{\sqrt{-1}}\iota_Z\omega$, and hence identifies this space with $\widehat{\operatorname{Sym}}Q^*
\otimes \operatorname{Sym}^{\leq r}P^*$.

The proof below implements this intuition at the cochain level. The projection $\pi_Q$ first reduces the complex $\Omega^*(M, \widehat{\operatorname{Sym}} Q^* \otimes \operatorname{Sym}^{\leq r} P^*)$ to $C^{0,*}(M,\operatorname{Sym}^{\leq r}P^*)$, while the
symplectic identification together with the transverse PBW map identifies
the latter with the complex of order-$r$ transverse differential operators. These two steps yield the quasi-isomorphism underlying Theorem \ref{Theorem: isomorphism of transverse differential operators}.

\begin{lemma}
	\label{Lemma: quasi-isomorphism for quantizable functions}
	For all $r \in \mathbb{N}$, the following map is a cochain quasi-isomorphism:
	\begin{equation}
		\label{Equation: quasi-isomorphism for quantizable functions}
		\pi_Q: (\Omega^*(M, \widehat{\operatorname{Sym}} Q^* \otimes \operatorname{Sym}^{\leq r} P^*), D_k) \to (C^{0, *}(M, \operatorname{Sym}^{\leq r} P^*), \pi_Q \circ D_k).
	\end{equation}
\end{lemma}
As mentioned previously, this lemma follows by the same arguments as Proposition \ref{Proposition: qausi-isomorphism for polarized sections}. For completeness, we include the proof in Appendix \ref{Proposition: qausi-isomorphism for polarized sections}.\par
Next, we construct the desired morphism. Define the $\mathcal{C}^\infty(M)$-linear map
\begin{equation*}
	\Pi: \Omega^*(M, \widehat{\operatorname{Sym}} Q^* \otimes \operatorname{Sym} P^*) \to C^{0, *}(M, \widetilde{D}(L^{(k)}, L^{(k)}))
\end{equation*}
by declaring that for $a \in \Omega^*(M, \widehat{\operatorname{Sym}} Q^* \otimes \operatorname{Sym} P^*)$ and $\sigma \in \Gamma(M, \widetilde{J}L^{(k)})$,
\begin{equation}
	\langle \Pi(a), \sigma \rangle := \pi_Q( a \circledast_k^{\operatorname{F}} (S^{L^{(k)}})^{-1} \sigma ).
\end{equation}
In particular, if $a \in \Gamma(M, \widehat{\operatorname{Sym}} Q^* \otimes \operatorname{Sym} P^*)$ and $s \in \Gamma(M, L^{(k)})$ is $d_P^{L^{(k)}}$-closed, then using $\mathfrak{j} s = S^{L^{(k)}}(J_s)$, where $\mathfrak{j} s$ denotes the infinite jet of $s$, we obtain
\begin{equation*}
	\langle \Pi(a), \mathfrak{j} s \rangle = \pi_Q(a \circledast_k^{\operatorname{F}} J_s) = a \circledast_k s.
\end{equation*}
Thus, $\Pi$ extends the action $\circledast_k$.

\begin{lemma}
	\label{Lemma: cochain isomorphism for transverse differential operators}
	For all $r \in \mathbb{N}$, the map $\Pi$ restricts to a cochain isomorphism
	\begin{equation}
		\label{Equation: cochain isomorphism of quantizable functinos}
		(C^{0, *}(M, \operatorname{Sym}^{\leq r} P^*), \pi_Q \circ D_k) \to (C^{0, *}(M, \widetilde{D}^r(L^{(k)}, L^{(k)})), d_P^{\widetilde{D}^r(L^{(k)}, L^{(k)})}).
	\end{equation}
\end{lemma}
\begin{proof}
	Fix $p \in \mathbb{N}$, $a \in C^{0, p}(M, \operatorname{Sym}^{\leq r} P^*)$ and $\sigma \in \Gamma(M, \widetilde{J} L^{(k)})$. Let $s = (S^{L^{(k)}})^{-1} (\sigma)$. Then
	\begin{align*}
		\left\langle \left( d_P^{\widetilde{D}^r(L^{(k)}, L^{(k)})} \circ \Pi \right)(a), \sigma \right\rangle = & d_P^{L^{(k)}} \langle \Pi(a), \sigma \rangle - (-1)^p \cdot \langle \Pi(a), d_P^{\widetilde{J}L^{(k)}} \sigma \rangle\\
		= & \pi_Q( \nabla^{\operatorname{K}, L^{(k)}} ( a \circledast_k^{\operatorname{F}} s ) - (-1)^p \cdot a \circledast_k^{\operatorname{F}} \nabla^{\operatorname{K}, L^{(k)}} s )\\
		= & \pi_Q ( ( D_k a ) \circledast_k^{\operatorname{F}} s )\\
		= & \langle (\Pi \circ \pi_Q \circ D_k) (a), \sigma \rangle.
	\end{align*}
	The second line uses $d_P^{L^{(k)}} \langle \Pi(a), \sigma \rangle = d_P^{L^{(k)}}\pi_Q(a \circledast_k^{\operatorname{F}} s) = \pi_Q( \nabla^{\operatorname{K}, L^{(k)}} (a \circledast_k^{\operatorname{F}} s))$ by Proposition \ref{Proposition: qausi-isomorphism for polarized sections}, and $\langle \Pi(a), d_P^{\widetilde{J}L^{(k)}} \sigma \rangle = \langle \Pi(a), \nabla^{\operatorname{G}, L^{(k)}} \sigma \rangle = \pi_Q( a \circledast_k^{\operatorname{F}} \nabla^{\operatorname{K}, L^{(k)}} s)$ by the property of $\pi_Q$. The last line uses $\Pi \circ \pi_Q = \Pi$ by the property of $\pi_Q$ again. Thus, \eqref{Equation: cochain isomorphism of quantizable functinos} is a cochain map.\par
	To show that \eqref{Equation: cochain isomorphism of quantizable functinos} is an isomorphism, we claim that the following diagram commutates:
	\begin{equation}
		\label{Equation: commutative diagram for action}
		\begin{tikzcd}
			C^{0, *}(M, \operatorname{Sym}^{\leq r} P^*) \ar[r, "\Pi"] \ar[d] & C^{0, *}(M, \widetilde{D}^r(L^{(k)}, L^{(k)}))\\
			C^{0, *}(M, \operatorname{Sym}^{\leq r} Q) \ar[ur, "\operatorname{pbw}^{L^{(k)}, L^{(k)}}"']
		\end{tikzcd}
	\end{equation}
	where the vertical arrow, denote by $a \mapsto \widetilde{a}$, is the isomorphism induced by (the inverse of) the isomorphism $Q \to P^*$, $Z \mapsto \tfrac{k}{\sqrt{-1}} \iota_Z \omega$. Our claim is proved by the observation that, for all $\alpha \in C^{0, *}(M)$, $a \in \Gamma(M, \operatorname{Sym}^{\leq r} P^*)$ and $s \in \Gamma(M, \operatorname{Sym}^{\leq r} Q^* \otimes L^{(k)})$,
	\begin{equation*}
		\pi_Q((\alpha \otimes a) \circledast_k^{\operatorname{F}} s) = \alpha \otimes \langle \widetilde{a}, s \rangle.
	\end{equation*}
	Since the other arrow $\operatorname{pbw}^{L^{(k)}, L^{(k)}}$ is also an isomorphism, the proof is complete.
\end{proof}

\begin{proof}[\myproof{Theorem}{\ref{Theorem: isomorphism of transverse differential operators}}]
	Recall that $\Pi \circ \pi_Q = \Pi$ by the property of $\pi_Q$. Then by Lemmas \ref{Lemma: quasi-isomorphism for quantizable functions} and \ref{Lemma: cochain isomorphism for transverse differential operators}, for each $r \in \mathbb{N}$, $\Pi$ restricts to a quasi-isomorphism
	\begin{equation*}
		(\Omega^*(M, \widehat{\operatorname{Sym}} Q^* \otimes \operatorname{Sym}^{\leq r} P^*), D_k) \to (C^{0, *}(M, \widetilde{D}^r(L^{(k)}, L^{(k)})), d_P^{\widetilde{D}^r(L^{(k)}, L^{(k)})}).
	\end{equation*}
	Hence, $\Pi$ induces an isomorphism of sheaves of filtered $\mathcal{C}_M^\infty$-modules $\mathcal{O}_k^{(<\infty)} \to \widetilde{\mathcal{D}}_{L^{(k)}}$. Arguing exactly as in Lemma \ref{Lemma: module structure}, one verifies that whenever $a$ and $b$ are level-$k$ quantizable functions, $\Pi (a) \circ \Pi(b) = \Pi(a \star_k^{\operatorname{F}} b)$, concluding the proof.
\end{proof}

Recall $\mathcal{W}^{\leq r} = \bigoplus_{i + j \leq r} \widehat{\operatorname{Sym}} Q^* \otimes \operatorname{Sym}^i P^* \cdot \hbar^j$. We end this subsection by stating a formula of the action $a \circledast_k$, for $a \in \Gamma(M, \mathcal{W}^{\leq r})$, in terms of higher covariant derivatives.
\begin{corollary}
	\label{Corollary: action in terms of higher covariant derivatives}
	Let $r \in \mathbb{N}$ and $a \in \Gamma(M, \mathcal{W}^{\leq r})$. Suppose $G \in \Gamma(M, \operatorname{Sym}^{\leq r} Q)$ is the image of $\left. \pi_Q(a) \right\vert_{\hbar = \frac{\sqrt{-1}}{k}} \in \Gamma(M, \operatorname{Sym}^{\leq r} P^*)$ under the vertical arrow in the diagram \eqref{Equation: commutative diagram for action}. Then
	\begin{equation}
		a \circledast_k s = \nabla_G^{L^{(k)}} s, \quad \text{for all } s \in \Gamma(M, L^{(k)}).
	\end{equation}
\end{corollary}
\begin{proof}
	It follows from the proof of Lemma \ref{Lemma: cochain isomorphism for transverse differential operators} and \eqref{Equation: pairing for classical jet} that $a \circledast_k s = \langle G, J_s \rangle = \nabla_G^{L^{(k)}} s$.
\end{proof}

\subsection{Examples}
\quad\par
\label{Subsection: examples}
In this subsection we examine several examples to which our framework applies, focusing in particular on how quantizable functions are characterized in each case.

\subsubsection{K\"ahler manifolds}
\quad\par
Consider a prequantizable K\"ahler manifold $(M, \omega)$ equipped with a prequantum line bundle $(L, \nabla^L)$. Let the polarization be $P = T^{0, 1}M$ and let $Q = T^{1, 0}M$. The Levi--Civita connection $\nabla$ of $(M, \omega)$ satisfies the conditions stated in Proposition \ref{Proposition: polarized symplectic connection}. Consequently, $\nabla$ can be used to construct the star product $\star$ in Theorem \ref{First main result} for any Hermitian holomorphic line bundle $\mathbf{L}$ endowed with its Chern connection $\nabla^\mathbf{L}$.\par
This construction agrees with that appearing in the second author’s work \cite{Yau2024, Yau2025}. The oppositely ordered analogue was established earlier in \cite{ChaLeuLi2022b}. In this setting, a necessary and sufficient condition for functions to be first order formal quantizable is given in \cite{LeuLiMa2024}. For additional background and related results, see \cite{BorWal1997, ChaLeuLi2023}, among others.

\subsubsection{Cotangent bundles}
\quad\par
The study of deformation quantization on cotangent bundles has its roots in classical Weyl quantization and pseudodifferential operator theory. A major advance specific to this case was achieved by Bordemann--Neumaier--Waldmann \cite{BorNeuWal1998}, who used a torsion-free connection $\nabla^B$ on a smooth manifold $B$ to construct a star product on $M = T^*B$ together with a differential operator representation. In this setting, $M$ is equipped with the symplectic form $\omega = d\theta$, where $\theta$ is the Liouville $1$-form on $M$, together with a canonical prequantum line bundle $(L, \nabla^L) = (M \times \mathbb{C}, d - \sqrt{-1}\theta)$ and the vertical real polarization $P$. We now revisit this case within our framework; for simplicity, we take $(\mathbf{L}, \nabla^\mathbf{L}) = \underline{\mathbb{C}}$.\par
Let us explain how a torsion-free connection on $B$ determines the auxiliary data in \eqref{Equation: auxiliary data for deformation quantization}. The dual connection of $\nabla^B$ on $T^*B$ defines a horizontal distribution $Q$ on the total space $M$, complementary to the vertical polarization $P$. Since $\nabla^B$ is torsion-free, this horizontal distribution is Lagrangian. Pulling back $\nabla^B$ along the bundle projection $M \to B$ yields a torsion-free connection on $Q$ (in the sense of Remark \ref{Remark: torsion free connection on quotient bundles}).\par
To characterize formal quantization functions, Lemma \ref{Lemma: characterization of formal quantizable functions} shows that it suffices to examine the $\Gamma(M, \widehat{\operatorname{Sym}} P^*)[[\hbar]]$-component of $\mathbf{J}_f$. In the present setting this component is given by $\sum_{r=0}^\infty (\delta_P^{-1} \nabla)^r f$, as follows from the recursive formula \eqref{Equation: recursive formula for double weight components} and the fact that
\begin{equation}
	\label{Equation 5.10}
	\widetilde{I} \in C^{1, 0}(M, \widehat{\operatorname{Sym}} Q^* \otimes P^*).
\end{equation}
The latter fact is justified as follows. For the Lie pair $(TM_\mathbb{C}, P)$ equipped with the auxiliary data \eqref{Equation: auxiliary data for deformation quantization} and the connection $\nabla$ on $TM_\mathbb{C}$ provided by Proposition \ref{Proposition: polarized symplectic connection}, the Atiyah cocycle $\mathcal{Z}^\nabla$ defined in (45) of Subsection 5.2 in \cite{LauStiXu2021} vanishes. Then Lemma 5.22 (2) of \cite{LauStiXu2021} implies \eqref{Equation 5.10}.\par
Applying Lemma \ref{Lemma: characterization of formal quantizable functions}, we conclude that a formal quantizable function on $M$ is precisely a formal function on $M$ polynomial in fibrewise variables and in $\hbar$. Equivalently,
\begin{equation}
	\mathcal{C}_{M, \hbar}^{(<\infty)} (M) \cong \Gamma(B, \operatorname{Sym} TB)[\hbar].
\end{equation}

\subsubsection{Symplectic tori with transversal real polarizations of compact leaves}
\quad\par
\label{Subsection: symplectic tori}
Consider $(M, \omega) = (\mathbb{R}^{2n} / \mathbb{Z}^{2n}, 2\pi \sum_{i=1}^n dx^i \wedge dy^i)$, where $(x^1, ..., x^n, y^1, ..., y^n)$ are periodic coordinates. This symplectic torus admits a prequantum line bundle $(L, \nabla^L)$, obtained as the quotient of the following Hermitian line bundle with unitary connection:
\begin{equation}
	\label{Equation: Unitary connection in symplectic torus case}
	(\mathbb{R}^{2n} \times \mathbb{C}, d - 2\pi \sqrt{-1} \textstyle\sum_{i=1}^n x^i dy^i)
\end{equation}
by the following factor of automorphy: for all $\mu, \nu \in \mathbb{Z}^n$ and $(x, y, c) \in \mathbb{R}^{2n} \times \mathbb{C}$,
\begin{equation*}
	(\mu, \nu) \cdot (x, y, c) = (x + \mu, y + \nu, e^{2\pi \sqrt{-1} \langle \mu, y \rangle} c).
\end{equation*}
Let $P$ and $Q$ be the real polarizations of $(M, \omega)$ generated respectively by the global frames $(\check{v}_1, ..., \check{v}_n) = (\tfrac{\partial}{\partial y^1}, ..., \tfrac{\partial}{\partial y^n})$ and $(v_1, ..., v_n) = (\tfrac{\partial}{\partial x^1}, ..., \tfrac{\partial}{\partial x^n})$. We take $(\mathbf{L}, \nabla^\mathbf{L}) = \underline{\mathbb{C}}$ and equip $TM_\mathbb{C} \cong M \times \mathbb{R}^{2n}$ with the trivial connection $d$. This connection satisfies the conclusion of Proposition \ref{Proposition: polarized symplectic connection} for an appropriate torsion-free connection on the quotient bundle $TM_\mathbb{C} / P$. In this setup the structural tensors $\widetilde{I}, F^\mathbf{L}$ vanish. Hence, for every $f \in \mathcal{C}^\infty(M)[[\hbar]]$,
\begin{equation}
	\label{Equation: quantum jet in symplectic torus case}
	\mathbf{J}_f = \sum_{p, q \geq 0} \frac{1}{(p+q)!} \frac{\partial^{p+q} f}{\partial x^{i_1} \cdots \partial x^{i_q} \partial y^{j_1} \cdots \partial y^{j_p}} u^{i_1} \cdots u^{i_q} \check{u}^{j_1} \cdots \check{u}^{j_p},
\end{equation}
as follows directly from the defining equation \eqref{Equation: defining equation of quantum jets}. For all $f, g \in \mathcal{C}^\infty(M)[[\hbar]]$, the associated star product (c.f. \cite{LeuYau2023}) takes the form
\begin{equation*}
	f \star g = \sum_{r=0}^\infty \frac{1}{r!} \left( \frac{\hbar}{2\pi} \right)^r \sum_{1 \leq i_1, ..., i_r \leq n} \frac{\partial^r f}{\partial y^{i_1} \cdots \partial y^{i_r}} \frac{\partial^r g}{\partial x^{i_1} \cdots \partial x^{i_r}}.
\end{equation*}
The explicit formula for $\mathbf{J}_f$ shows that a formal function $f \in \mathcal{C}^\infty(M)[[\hbar]]$ is formal quantizable if and only if it is independent of the variables $y$ and polynomial in $\hbar$. Consequently,
\begin{equation}
	\mathcal{C}_{M, \hbar}^{(<\infty)} (M) = \mathcal{C}^\infty(\mathbb{R}^n / \mathbb{Z}^n)[\hbar].
\end{equation}
This suggests that, when the projection $M \to M / P_{\operatorname{r}}$ onto the leaf space is a smooth submersion with compact leaves --- where $P_{\operatorname{r}}$ is the real distribution on $M$ defined by $P_{\operatorname{r}} \otimes \mathbb{C} = P \cap \overline{P}$ --- restricting attention to global formal quantizable functions may be overly limiting.

\section{Asymptotics of Toeplitz-type operators in real polarization}
\label{Section: asymptotics of Toeplitz operators in real polarization}
Toeplitz-type operators in the real polarization setting were first introduced in \cite{LeuYau2022}. This section explains how an asymptotic expansion of these operators arises naturally from the action $\circledast_k$. For real or mixed polarizations, one often need to work with distributional sections (see \cite{BaiFloMouNun2011, BaiMouNun2010, LeuWan2023}). Subsection \ref{Subsection: a lift to an action on distributional sections} describes how the action $\circledast_k$ extends to distributional sections of $L^{(k)} = L^{\otimes k} \otimes \mathbf{L}$. Subsection \ref{Subsection: asymptotics of Toeplitz-type operators in real polarization} then uses this extension to establish Theorem \ref{Fourth main result}.

\subsection{A lift to an action on distributional sections}
\quad\par
\label{Subsection: a lift to an action on distributional sections}
Let $k \in \mathbb{Z}^+$. Via the Liouville volume form on $(M, \omega)$, we identify the space $\Gamma^{-\infty}(M, L^{(k)})$ of distributional sections on $L^{(k)}$ with $(\Gamma_{\operatorname{c}}(M, L^{(-k)}))'$. Here, the prime $'$ denotes the continuous dual with respect to the usual locally convex topology on the space $\Gamma_{\operatorname{c}}(M, L^{(-k)})$ of compactly supported smooth sections of $L^{(-k)} := (L^*)^{\otimes k} \otimes \mathbf{L}^*$. There are a natural pairing
\begin{equation*}
	\langle \, \cdot\,, \, \cdot\, \rangle: \Gamma^{-\infty}(M, L^{(k)}) \times \Gamma_{\operatorname{c}}(M, L^{(-k)}) \to \mathbb{C},
\end{equation*}
and a natural embedding $\iota: \Gamma(M, L^{(k)}) \to \Gamma^{-\infty}(M, L^{(k)})$ given by $\langle \iota(s), \tau \rangle = \int_M \tau(s) \cdot \tfrac{1}{n!} \omega^n$. For each $X \in \Gamma(M, TM_\mathbb{C})$, the covariant derivative $\nabla_X^{L^{(k)}}$ acts on a distributional section $s \in \Gamma^{-\infty}(M, L^{(k)})$ in the weak sense: for any test section $\tau \in \Gamma_{\operatorname{c}}(M, L^{(-k)})$,
\begin{equation*}
	\langle \nabla_X^{L^{(k)}} s, \tau \rangle := -\langle s, \nabla_X^{L^{(-k)}} \tau + \operatorname{div}(X) \tau \rangle.
\end{equation*}
Here, $\operatorname{\operatorname{div}}(X)$ is the divergence of $X$, i.e. $\operatorname{div}(X) \cdot \omega^n = \mathcal{L}_X (\omega^n)$. Higher covariant derivatives $\nabla_G^{L^{(k)}}$, for $G \in \Gamma(M, \operatorname{Sym} Q)$, extends to distributional sections in the same weak sense; for instance, when $f \in \mathcal{C}^\infty(M)$, $\langle \nabla_f^{L^{(k)}} s, \tau \rangle := \langle f \cdot s, \tau \rangle = \langle s, f \cdot \tau \rangle$. Then Corollary \ref{Corollary: action in terms of higher covariant derivatives} allows us to extend, for each $a \in \Gamma(M, \mathcal{W}^{\leq r})$, the operator $a \circledast_k$ to a continuous linear operator
\begin{equation*}
	\Gamma^{-\infty}(M, L^{(k)}) \to \Gamma^{-\infty}(M, L^{(k)}),
\end{equation*}
which we denote by the same symbol. We can now introduce the following operators.

\begin{definition}
	\label{Definition: components of actions}
	Let $f \in \mathcal{C}^\infty(M)$, $k \in \mathbb{Z}^+$ and $r \in \mathbb{N}$. We define
	\begin{equation}
		\label{Equation: action of weight component of quantum jets}
		T_{f, r, k} := (\mathbf{J}_f)_r \circledast_k: \Gamma^{-\infty}(M, L^{(k)}) \to \Gamma^{-\infty}(M, L^{(k)}),
	\end{equation}
	where $(\mathbf{J}_f)_r$ denotes the weight-$r$ component of $\mathbf{J}_f$ with respect to the polarized weight.
\end{definition}

Recall that $s \in \Gamma^{-\infty}(M, L^{(k)})$ is called $P$-\emph{polarized} if $\nabla_Y^{L^{(k)}} s = 0$ for all $Y \in \Gamma(M, P)$. We focus on the action of $T_{f, r, k}$ on the space $\Gamma_P^{-\infty}(M, L^{(k)})$ of $P$-polarized distributional sections of $L^{(k)}$, which will be analyzed in a real-polarized case in the next subsection.

\subsection{Asymptotics of Toeplitz-type operators in real polarization}
\quad\par
\label{Subsection: asymptotics of Toeplitz-type operators in real polarization}
In the setting of Subsection \ref{Subsection: symplectic tori}, we identify $\mathcal{C}^\infty(M)$ with the space of periodic functions on $\mathbb{R}^{2n}$. Then every $f \in \mathcal{C}^\infty(M)$ admits a fibrewise Fourier expansion
\begin{equation}
	f(x, y) = \sum_{p \in \mathbb{Z}^n} \widehat{f}_p(x) e^{2\pi\sqrt{-1} \langle p, y \rangle},
\end{equation}
where $\widehat{f}_p$ is the $p$th fibrewise Fourier coefficient of $f$ (see Definition 4 in \cite{LeuYau2023}). Similarly, for $k \in \mathbb{Z}^+$, $\Gamma(M, (L^*)^{\otimes k})$ corresponds to functions $\tau \in \mathcal{C}^\infty(\mathbb{R}^{2n})$ satisfying
\begin{equation*}
	\tau(x + \mu, y + \nu) = e^{-2\pi k\sqrt{-1} \langle \mu, y \rangle} \tau(x, y) \quad \text{for all } \mu, \nu \in \mathbb{Z}^n.
\end{equation*}
Every section $\tau \in \Gamma(M, (L^*)^{\otimes k})$ can then be expressed as
\begin{equation}
	\tau(x, y) = \sum_{[m] \in (\mathbb{Z} / k\mathbb{Z})^n} \sum_{q \in \mathbb{Z}} \widetilde{\tau}_m (x + q) e^{2\pi\sqrt{-1} \langle -m + kq, y \rangle},
\end{equation}
where for each $m \in \mathbb{Z}^n$, $\widetilde{\tau}_m$ is the $m$th Weil--Brezin coefficient of $\tau$ (c.f. Definition 5 in \cite{LeuYau2023}), which is a Schwartz function on $\mathbb{R}^n$ satisfying the condition that
\begin{equation*}
	\widetilde{\tau}_{m+kq}(x) = \widetilde{\tau}_m(x - q) \quad \text{for all } x \in \mathbb{R}^n \text{ and } q \in \mathbb{Z}^n.
\end{equation*}
The space $\Gamma_P^{-\infty}(M, L^{\otimes k})$ admits a basis $\{ \sigma_k^{[m]} \}_{[m] \in (\mathbb{Z} / k \mathbb{Z})^n}$ over $\mathbb{C}$ \cite{BaiMouNun2010}, where $[m] := m + k\mathbb{Z}^n$ denotes the equivalence class of $m \in \mathbb{Z}^n$ in $\mathbb{Z}^n / (k\mathbb{Z})^n$ and
\begin{equation}
	\langle \sigma_k^{[m]}, \tau \rangle := \int_{[0, 1]^n} e^{2\pi\sqrt{-1} \langle m, y \rangle} \tau(\tfrac{m}{k}, y) d^n y = \widetilde{\tau}_m (\tfrac{m}{k}).
\end{equation}
In \cite{LeuYau2023}, Leung and the second author studied Toeplitz-type operators associated with the pair of real polarizations $(P, Q)$ and analyzed the resulting asymptotic action of the corresponding star product. Building on that framework, we now consider, for each $f \in \mathcal{C}^\infty(M)$, the operator
\begin{equation*}
	T_{f, k}: \Gamma_P^{-\infty}(M, L^{\otimes k}) \to \Gamma_P^{-\infty}(M, L^{\otimes k}),
\end{equation*}
defined by the same formula as in (6.6) of \cite{LeuYau2023}. This operator is completely determined by its action on the basis $\{ \sigma_k^{[m]} \}_{[m] \in (\mathbb{Z} / k \mathbb{Z})^n}$: for all $m \in \mathbb{Z}^n$ and $\tau \in \Gamma(M, (L^*)^{\otimes k})$,
\begin{equation}
	\label{Equation: formula of Toeplitz operator}
	\left\langle T_{f, k} \sigma_k^{[m]}, \tau \right\rangle = \sum_{p \in \mathbb{Z}^n} \left( \widehat{f}_p \cdot \widetilde{\tau}_{m+p} \right) \left( \frac{m+p}{k} \right).
\end{equation}
This yields the opposite-ordered analogue of the Toeplitz-type operator associated with $(P, Q)$ appearing in Definition 6.1 of the same reference (without half-form correction).\par
Our final goal is to compare this operator $T_{f, k}$ with the operators $T_{f, r, k}$ introduced in \eqref{Equation: action of weight component of quantum jets}.

\begin{proposition}
	Let $k \in \mathbb{Z}^+$ and $r \in \mathbb{N}$. For all $m \in \mathbb{Z}^n$ and $\tau \in \Gamma(M, (L^*)^{\otimes k})$,
	\begin{equation}
		\label{Equation: formula of asymptotic expansion of Toeplitz operator}
		\left\langle T_{f, r, k} \sigma_k^{[m]}, \tau \right\rangle = \sum_{p \in \mathbb{Z}^n} \frac{1}{r!} \cdot \frac{\partial^r (\widehat{f}_p \cdot \widetilde{\tau}_{m+p})}{\partial x^{i_1} \cdots \partial x^{i_r}} \left( \frac{m}{k} \right) \cdot \frac{p_{i_1} \cdots p_{i_r}}{k^r}.
	\end{equation}
\end{proposition}
\begin{proof}
	From \eqref{Equation: Unitary connection in symplectic torus case} we know that $\nabla_{\partial_{x^i}}^{L^{\otimes k}} = \tfrac{\partial}{\partial x^i}$. Thus, by \eqref{Equation: quantum jet in symplectic torus case} and Corollary \ref{Corollary: action in terms of higher covariant derivatives},
	\begin{equation*}
		(\mathbf{J}_f)_r \circledast_k = \sum_{p \in \mathbb{Z}^n} \frac{1}{r!} \cdot \frac{p_{i_1} \cdots p_{i_r}}{(-k)^r} \cdot \widehat{f}_p(x) e^{2\pi\sqrt{-1} \langle p, y \rangle} \frac{\partial^r}{\partial x^{i_1} \cdots \partial x^{i_r}} \quad \text{on } \Gamma(M, L^{\otimes k}).
	\end{equation*}
	Applying $(\mathbf{J}_f)_r \circledast_k$ to $\sigma_k^{[m]}$ and pairing with $\tau$, integration by parts (in the weak sense) yields
	\begin{align*}
		\left\langle T_{f, r, k} \sigma_k^{[m]}, \tau \right\rangle = & \sum_{p \in \mathbb{Z}^n} \left\langle \sigma_k^{[m]}, \frac{1}{r!} \cdot \frac{p_{i_1} \cdots p_{i_r}}{k^r} \frac{\partial^r (\widehat{f}_p(x) e^{2\pi\sqrt{-1} \langle p, y \rangle} \cdot \tau)}{\partial x^{i_1} \cdots \partial x^{i_r}} \right\rangle\\
		= & \sum_{p \in \mathbb{Z}^n} \frac{1}{r!} \cdot \frac{p_{i_1} \cdots p_{i_r}}{k^r} \int_{[0, 1]^n} e^{2\pi\sqrt{-1} \langle m + p, y \rangle} \frac{\partial^r (\widehat{f}_p \cdot \tau)}{\partial x^{i_1} \cdots \partial x^{i_r}} \left(\frac{m}{k}, y\right) d^ny\\
		= & \sum_{p \in \mathbb{Z}^n} \frac{1}{r!} \cdot \frac{\partial^r (\widehat{f}_p \cdot \widetilde{\tau}_{m+p})}{\partial x^{i_1} \cdots \partial x^{i_r}} \left( \frac{m}{k} \right) \cdot \frac{p_{i_1} \cdots p_{i_r}}{k^r}.
	\end{align*}
\end{proof}

To interpret this formula, consider a single Fourier mode $f(x, y) = \widehat{f}_p(x) e^{2\pi\sqrt{-1}\langle p, y \rangle}$. By the above computation, the sum $\sum_{r=0}^N T_{f, r, k} \sigma_k^{[m]}$ applied to $\tau$ reproduces the Taylor expansion of $\left\langle T_{f, k} \sigma_k^{[m]}, \tau \right\rangle$ in the $x$-variable at the point $\tfrac{m}{k}$ up to order $N$.\par
To refine this observation, we seek an explicit bound on the remainder term, analogous to the K\"ahler case treated in \cite{ChaLeuLiYau2025}. A technical obstacle arises because there is no canonical $L^2$-norm on distributional sections, so the norm estimates used in the K\"ahler setting cannot be applied directly.\par
To bypass this, recall that the Weil-Brezin transformation identifies the topological vector space $\Gamma(M, (L^*)^{\otimes k})$ with a product of $k^n$ copies of the Schwartz space $\mathcal{S}(\mathbb{R}^n)$, each equipped with its standard Fr\'echet topology \cite{BaiMouNun2010}. Accordingly, for each $N \in \mathbb{N}$, we can define a norm on $\Gamma(M, (L^*)^{\otimes k})$ by
\begin{equation}
	\label{Equation: seminorm}
	\lVert \tau \rVert_N := \sup \left\{ \left\lvert \frac{\partial^r \widetilde{\tau}_m}{\partial x^{i_1} \cdots \partial x^{i_r}}(x) \right\rvert: [m] \in (\mathbb{Z} / k\mathbb{Z})^n, r \leq N, 1 \leq i_1, ..., i_r \leq n, x \in \mathbb{R}^n \right\}.
\end{equation}
This supremum is finite because each component $\widetilde{\tau}_m$ is a Schwartz function. Dually, on $\Gamma_P^{-\infty}(M, L^{\otimes k})$, we may introduce an $\ell^1$-norm $\lVert \, \cdot \, \rVert_{\ell^1}$, normalized with respect to the basis $\{ \sigma_k^{[m]} \}_{[m] \in (\mathbb{Z} / k\mathbb{Z})^n}$: for $s = \sum_{[m] \in (\mathbb{Z} / k\mathbb{Z})^n} s_{[m]} \sigma_k^{[m]} \in \Gamma_P^{-\infty}(M, L^{\otimes k})$, we set
\begin{equation}
	\label{Equation: l1 norm}
	\lVert s \rVert_{\ell^1} := \sum_{[m] \in (\mathbb{Z} / k \mathbb{Z})^n} \lvert s_{[m]} \rvert.
\end{equation}
These norms give uniform control over distributional pairings: for every $N \in \mathbb{N}$,
\begin{equation*}
	\lvert \langle s, \tau \rangle \rvert \leq \lVert s \rVert_{\ell^1} \cdot \lVert \tau \rVert_N,
\end{equation*}
which follows directly from H\"older’s inequality. We can now state our final main result.

\begin{theorem}[Theorem \ref{Fourth main result}]
	\label{Theorem: symplectic tori}
	Let $M = \mathbb{R}^{2n} / \mathbb{Z}^{2n}$ be the standard symplectic torus. Suppose $f \in \mathcal{C}^\infty(M)$ and $N \in \mathbb{N}$. Then there exists a constant $C_{f, N} > 0$ such that for all $k \in \mathbb{Z}^+$, $s \in \Gamma_P^{-\infty}(M, L^{\otimes k})$ and test section $\tau \in \Gamma(M, (L^*)^{\otimes k})$,
	\begin{equation*}
		\left\lvert \left\langle T_{f, k} s - \sum_{r=0}^N T_{f, r, k} s, \tau \right\rangle \right\rvert \leq C_{f, N} \cdot \lVert s \rVert_{\ell^1} \cdot \lVert \tau \rVert_{N+1} \cdot \frac{1}{k^{N+1}}.
	\end{equation*}
\end{theorem}
\begin{proof}
	We first recall an a priori estimate from Lemma 2 of \cite{LeuYau2023}: there exists $\tilde{C}_{f, N} > 0$ such that for all $p \in \mathbb{Z}^n$, $1 \leq j \leq N+1$, $1 \leq i_1, ..., i_{l+1}, i'_1, ..., i'_j \leq n$ and $x \in \mathbb{R}^n$,
	\begin{equation}
		\label{Equation: a priori estimate}
		\left\lvert \frac{\partial^j \widehat{f}_p}{\partial x^{i'_1} \cdots \partial x^{i'_j}} (x) \cdot p_{i_1} \cdots p_{i_{N+1}} \right\rvert \leq \widetilde{C}_{f, N} \cdot N_p^2,
	\end{equation}
	where $N_p := \prod_{i \in \{1, ..., n\}: p_i \neq 0} \tfrac{1}{p_i}$. Now, fix $k \in \mathbb{Z}^+$, $m \in \mathbb{Z}^n$ and $\tau \in \Gamma(M, (L^*)^{\otimes k})$. Define $R_{f, N, k} := T_{f, k} - \sum_{r=0}^N T_{f, r, k}$. Using the explicit formulas \eqref{Equation: formula of Toeplitz operator} and \eqref{Equation: formula of asymptotic expansion of Toeplitz operator}, we obtain
	\begin{align*}
		\left\lvert \left\langle R_{f, N, k} \sigma_k^{[m]}, \tau \right\rangle \right\rvert \leq & \sum_{p \in \mathbb{Z}^n} \left\lvert \left( \widehat{f}_p \cdot \widetilde{\tau}_{m+p} \right) \left( \frac{m+p}{k} \right) - \sum_{r=0}^N \frac{1}{r!} \cdot \frac{\partial^r (\widehat{f}_p \cdot \widetilde{\tau}_{m+p})}{\partial x^{i_1} \cdots \partial x^{i_r}} \left( \frac{m}{k} \right) \cdot \frac{p_{i_1} \cdots p_{i_r}}{k^r} \right\rvert.
	\end{align*}
	By Taylors theorem, for each $p$ there exists $t_p \in [0, 1]$ such that the inner difference equals
	\begin{equation}
		\label{Taylor remainder}
		\frac{1}{(N+1)!} \cdot \frac{\partial^{N+1} (\widehat{f}_p \cdot \widetilde{\tau}_{m+p})}{\partial x^{i_1} \cdots \partial x^{i_{N+1}}} \left( \frac{m + t_p p}{k} \right) \cdot \frac{p_{i_1} \cdots p_{i_{N+1}}}{k^{N+1}}.
	\end{equation}
	Let $\mathcal{S}_{j, N+1-j}$ denote the set of $(j, N+1-j)$-shuffles. The Leibniz rule, the definition \eqref{Equation: seminorm} and the estimate \eqref{Equation: a priori estimate} show that the absolute value of \eqref{Taylor remainder} is bounded by
	\begin{align*}
		& \sum_{j=0}^{N+1} \sum_{\sigma \in \mathcal{S}_{j, N+1-j}} \left\lvert \frac{p_{i_1} \cdots p_{i_{N+1}}}{(N+1)! \cdot k^{N+1}} \cdot \left( \frac{\partial^j \widehat{f}_p}{\partial x^{i_{\sigma(1)}} \cdots \partial x^{i_{\sigma(j)}}} \cdot \frac{\partial^{N+1-j} \widetilde{\tau}_{m+p}}{\partial x^{i_{\sigma(j+1)}} \cdots \partial x^{i_{\sigma(N+1)}}} \right) \left( \frac{m + t_p p}{k} \right) \right\rvert\\
		\leq & \frac{1}{k^{N+1}} \cdot \frac{(2n)^{N+1}}{(N+1)!} \cdot \widetilde{C}_{f, N} \cdot N_p^2 \cdot \lVert \tau \rVert_{N+1}.
	\end{align*}
	Here, $2^{N+1} = \sum_{j=0}^{N+1} \lvert \mathcal{S}_{j, N+1-j} \rvert$ comes from the binomial theorem, and $n^{N+1}$ accounts for the number of index sequences of length $N+1$. Summing over $p \in \mathbb{Z}^n$ gives
	\begin{equation*}
		\lvert \langle R_{f, N, k} \sigma_k^{[m]}, \tau \rangle \rvert \leq \tfrac{1}{k^{N+1}} \cdot C_{f, N} \cdot \lVert \tau \rVert_{N+1} \quad \text{with} \quad C_{f, N} := \tfrac{(2n)^{N+1}}{(N+1)!} \cdot \widetilde{C}_{f, N} \cdot \textstyle\sum_{p \in \mathbb{Z}^n} N_p^2 < \infty.
	\end{equation*}
	For $s \in \Gamma_P^{-\infty}(M, L^{\otimes k})$, $\left\lvert \langle R_{f, N, k} (s), \tau \rangle \right\rvert \leq \frac{1}{k^{N+1}} \cdot C_{f, N} \cdot \lVert s \rVert_{\ell^1} \cdot \lVert \tau \rVert_{N+1}$ by Hölder's inequality.
\end{proof}

\begin{remark}
	As discussed in Subsection \ref{Subsection: symplectic tori}, if $f \in \mathcal{C}^\infty(M)$ is formal quantizable, then $f(x, y) = \widehat{f}_0(x)$. In this case, $T_{f, r, k} = 0$ for all $r > 0$, so that $T_{f, k} = T_{f, 0, k} = \sum_{r=0}^\infty T_{f, r, k}$.
\end{remark}

In future work, we hope to study analogous relationships for other --- potentially even singular --- polarizations, whenever a suitable Toeplitz-type calculus becomes available.

\appendix

\section{Proof of Proposition \ref{Proposition: transverse jet bundles are smooth subbundles}}
\label{Section: Proof of locally freeness of transverse jet bundles}
To prove Proposition \ref{Proposition: transverse jet bundles are smooth subbundles}, it is convenient to work on a $P$-\emph{foliated chart} on an open subset $U$ of $M$, which means real coordinates $(x^1, ..., x^n)$ on $U$ for which
\begin{equation*}
	\left( \frac{\partial}{\partial \check{u}^1}, ..., \frac{\partial}{\partial \check{u}^{k+l}} \right)
\end{equation*}
is a local complex frame of $P$. Here, $k$ and $k + l$ are the ranks of $P \cap \overline{P}$ and $P$ respectively. Also, $\check{u}^i := x^i - \sqrt{-1} x^{l+i}$ for each $1 \leq i \leq l$ and $\check{u}^{l+j} := x^{2l+j}$ for each $1 \leq j \leq k$. Indeed, $M$ always has an open cover in which each open subset has a $P$-foliated chart on it \cite{Nir1957}.\par
Let $d_P^{E^*}$ be the dual $P$-connection of $d_P^E$. A $P$-\emph{foliated fibred chart} on $E \vert_U$ means a tuple
\begin{equation*}
	(x^1, ..., x^n, \xi^1, ..., \xi^m)
\end{equation*}
for which
\begin{enumerate}
	\item $(x^1, ..., x^n)$ is a $P$-foliated chart on $U$; and
	\item $(\xi^1, ..., \xi^m)$ is a local $d_P^{E^*}$-closed frame of $E^*$ over $U$.
\end{enumerate}
In what follows, we define $u^i := x^i + \sqrt{-1} x^{l+i}$ for each $1 \leq i \leq l$ and $u^{l + j} := x^{2l+k+j}$ for each $1 \leq j \leq n-2l-k$. Also, we always decompose a multi-index $\mathcal{I}$ on $\{1, ..., n\}$ into a pair $(I, J)$, where $I$ is a multi-index on $\{1, ..., l+k\}$ and $J$ is a multi-index on $\{1, ..., n-l-k\}$. Then such a $P$-foliated fibred chart induces a family of local functions $\xi_{\mathcal{I}}^\alpha \in \mathcal{C}^\infty( J^rE \vert_U )$ parametrized by an index $1 \leq \alpha \leq m$ and an multi-index $\mathcal{I} = (I, J)$ on $\{1, ..., n\}$ with $\lvert \mathcal{I} \rvert \leq r$, which are defined as follows. For all $x \in U$ and $\sigma \in J_x^rE$,
\begin{equation*}
	\xi_{\mathcal{I}}^\alpha(\sigma) := \frac{\partial^{\lvert \mathcal{I} \rvert} (\xi^\alpha \circ s)}{\partial \check{u}^I \partial u^J} (x),
\end{equation*}
where $s$ is any local smooth section of $E$ around $x$ representing $\sigma$. They form a local trivialization of $J^rE$ over $U$:
\begin{equation}
	\label{Equation: local trivialization of jet bundle}
	\Gamma(U, J^rE) \to \mathcal{C}^\infty(U)^{\oplus m \binom{n+r}{r}}, \quad \sigma \mapsto ( \xi_\mathcal{I}^\alpha \circ \sigma )_{1 \leq \alpha \leq m, \lvert \mathcal{I} \rvert \leq r}.
\end{equation}

\begin{lemma}
	\label{Lemma: coordinate description of transverser jet bundle}
	Let $r \geq 1$ and $\sigma \in \Gamma(U, J^rE)$. Then the following conditions are equivalent:
	\begin{enumerate}
		\item $\sigma \in \Gamma(U, \widetilde{J}^rE)$.
		\item For all $1 \leq \alpha \leq m$ and multi-indices $\mathcal{I} = (I, J)$ on $\{1, ..., n\}$ with $\lvert \mathcal{I} \rvert \leq r$ and $\lvert I \rvert > 0$,
		\begin{equation*}
			\xi_{\mathcal{I}}^\alpha \circ \sigma = 0.
		\end{equation*}
		\item For all $x \in U$, $\sigma(x)$ is represented by the germ of a local $d_P^E$-closed section of $E$ at $x$.
	\end{enumerate}
\end{lemma}
\begin{proof}
	Let $(e_1, ..., e_m)$ denote the local $d_P^E$-closed frame of $E$ whose dual frame is $(\xi^1, ..., \xi^m)$.\par
	(1) $\Rightarrow$ (2): Fix $x \in U$. Then $\sigma(x)$ is represented by the germ of a local smooth section $s_x$ of $E$ over a sufficiently small neighbourhood of $x$ defined by
	\begin{equation}
		\label{Equation: local representative of a jet}
		s_x(y) = \sum_{\lvert \mathcal{I} \rvert \leq r} \frac{(\xi_\mathcal{I}^\alpha \circ \sigma) (x)}{\mathcal{I}!} (\check{u}(y) - \check{u}(x))^I (u(y) - u(x))^J e_\alpha(y).
	\end{equation}
	Consider a particular index $1 \leq \alpha \leq m$ and a multi-index $\mathcal{I} = (I, J)$ with $\lvert \mathcal{I} \rvert \leq r$ and $\lvert I \rvert > 0$. Pick $i \in \{1, ..., l+k\}$ such that $I(i) > 0$. Then we can define a new multi-index $\mathcal{I}' = (I', J')$ by $I'(j) = I(j) - \delta_{ij}$ for all $1 \leq j \leq l + k$ and $J' = J_0$, where $\delta_{ij}$ denotes the Driac delta function. Then we can define a differential operator $\Phi: \Gamma(U, E) \to \mathcal{C}^\infty(U)$ of order $r-1$ by
	\begin{equation*}
		\Phi(f^\beta e_\beta) = \frac{\partial^{\lvert \mathcal{I}' \rvert} f^\alpha}{\partial \check{u}^{I'} \partial u^{J'}}, \quad \text{for all } f^1, ..., f^m \in \mathcal{C}^\infty(U).
	\end{equation*}
	Since $e_1, ..., e_m$ are $d_P^E$-closed, we obtain
	\begin{equation*}
		(\xi_\mathcal{I}^\alpha \circ \sigma)(x) = \left( \left( \Phi \circ \iota_{\partial_{\check{u}^i}} \circ d_P^E \right) s_x \right) (x) = \left\langle \Phi \circ \iota_{\partial_{\check{u}^i}} \circ d_P^E, \sigma \right\rangle (x) = 0.
	\end{equation*}
	
	(2) $\Rightarrow$ (3): Fix $x \in U$. We claim that the local section $s_x$ of $E$ defined as in \eqref{Equation: local representative of a jet}, whose germ at $x$ represents $\sigma(x)$, is $d_P^E$-closed. Fix $1 \leq i \leq l + k$. Since $e_1, ..., e_m$ are $d_P^E$-closed,
	\begin{equation*}
		\iota_{\partial_{\check{u}^i}} (d_P^E s_x) (y) = \sum_{\lvert J \rvert \leq r} \left( \frac{\partial}{\partial \check{u}^i} \left( \frac{(\xi_{\emptyset, J}^\alpha \circ \sigma)(x)}{J!} (u(y) - u(x))^J \right) \right) e_\alpha(y) = 0,
	\end{equation*}
	for all $y$ in a neighbourhood of $x$. Therefore, our claim holds.\par
	(3) $\Rightarrow$ (1): This is trivial.
\end{proof}

\begin{proof}[\myproof{Proposition}{\ref{Proposition: transverse jet bundles are smooth subbundles}}]
	Identify local smooth sections of $J^rE$ with tuples of local functions $(\zeta_\mathcal{I}^\alpha)_{1 \leq \alpha \leq m, \mathcal{I} \leq r}$ under the local trivialization given in \eqref{Equation: local trivialization of jet bundle}. We also take a similar identification of local smooth sections of $J^{r-1}E$. Under these identifications, the projection map
	\begin{equation*}
		\Gamma(U, J^rE) \to \Gamma(U, J^{r-1}E)
	\end{equation*}
	is given by forgetting the functions $\zeta_{\mathcal{I}}^\alpha$'s with $1 \leq \alpha \leq m$ and $\lvert \mathcal{I} \rvert = r$.\par
	Correspondingly, the kernel $\Gamma(U, \operatorname{Sym}^r T^*M_\mathbb{C} \otimes E)$ of this projection map is identified with the subspace of $\mathcal{C}^\infty(U)^{\oplus m \binom{n+r}{r}}$ defined by the equations
	\begin{equation*}
		\zeta_{\mathcal{I}}^\alpha = 0, \quad \text{for } 1 \leq \alpha \leq m \text{ and } \lvert \mathcal{I} \rvert < r.
	\end{equation*}
	Furthermore, the subspace $\Gamma(U, \operatorname{Sym} Q^* \otimes E)$ of $\Gamma(U, \operatorname{Sym}^r T^*M_\mathbb{C} \otimes E)$ is identified with the space of tuples $(\zeta_{\mathcal{I}}^\alpha)_{1 \leq \alpha \leq m, \lvert \mathcal{I} \rvert \leq r} \in \mathcal{C}^\infty(U)^{\oplus m \binom{n+r}{r}}$ in which all functions but those with $\lvert J \rvert = r$ are zero. Then it clearly follows from Lemma \ref{Lemma: coordinate description of transverser jet bundle} that
	\begin{center}
		\begin{tikzcd}
			0 \ar[r] & \Gamma(U, \operatorname{Sym} Q^* \otimes E) \ar[r] & \Gamma(U, \widetilde{J}^rE) \ar[r] & \Gamma(U, \widetilde{J}^{r-1}E) \ar[r] & 0
		\end{tikzcd}
	\end{center}
	is a short exact sequence of $\mathcal{C}^\infty(U)$-modules. Consequently, we obtain a short exact sequence of $\mathcal{C}_M^\infty$-modules:
	\begin{center}
		\begin{tikzcd}
			0 \ar[r] & \Gamma(-, \operatorname{Sym}^r Q^* \otimes E) \ar[r] & \Gamma(-, \widetilde{J}^rE) \ar[r] & \Gamma(-, \widetilde{J}^{r-1}E) \ar[r] & 0
		\end{tikzcd}
	\end{center}
	It implies that, if $\Gamma(-, \widetilde{J}^{r-1}E)$ is locally free of finite rank, then so is $\Gamma(-, \widetilde{J}^rE)$. We conclude our proof by induction that $\Gamma(-, \widetilde{J}^rE)$ is a locally free $\mathcal{C}_M^\infty$-module for all $r \in \mathbb{N}$.
\end{proof}

\section{Proofs of Proposition \ref{Proposition: qausi-isomorphism for polarized sections} and Lemma \ref{Lemma: quasi-isomorphism for quantizable functions}}
\label{Sectoin: a degeneracy condition for filtered complexes}
In this appendix, we prove Proposition \ref{Proposition: quasi isomorphism}, from which Proposition \ref{Proposition: qausi-isomorphism for polarized sections} and Lemma \ref{Lemma: quasi-isomorphism for quantizable functions} follow as corollaries. We begin by describing the setting.\par
Let $M$ be a smooth manifold endowed with a Nirenberg integrable complex distribution $P$, a complementary subbundle $Q \subset TM_\mathbb{C}$, and a complex vector bundle $E \to M$. Assume that
\begin{itemize}
	\item $D$ is a flat connection on $\widehat{\operatorname{Sym}} Q^* \otimes E$, and
	\item $D + \delta_Q$ preserves $\widehat{\operatorname{Sym}}^{\geq l} Q^* \otimes E$ for all $l \in \mathbb{N}$.
\end{itemize}
We then consider the $\mathbb{C}$-linear map
\begin{equation}
	\label{Equation: truncation of flat connection}
	d_P: \Omega^*(M, \widehat{\operatorname{Sym}} Q^* \otimes E) \to \Omega^*(M, \widehat{\operatorname{Sym}} Q^* \otimes E)
\end{equation}
defined by: for any $p, q, l \in \mathbb{N}$, the restriction $d_P$ to $C^{q, p}(M, \operatorname{Sym}^l Q^* \otimes E)$ is the composition
\begin{center}
	\begin{tikzcd}
		C^{q, p}(M, \operatorname{Sym}^l Q^* \otimes E) \ar[r, "D"] & \Omega^{p+q+1}(M, \widehat{\operatorname{Sym}} Q^* \otimes E) \ar[r] & C^{q, p+1}(M, \operatorname{Sym}^l Q^* \otimes E)
	\end{tikzcd}
\end{center}
where the second arrow denotes the canonical projection.\par
Let us examine some basic properties of $d_P$. Note that $P$ is involutive. Then for a section $s \in C^{q, p}(M, \operatorname{Sym}^l Q^* \otimes E)$ with $q, p, l \in \mathbb{N}$, our assumptions imply that $Ds$ lies in
\begin{equation*}
	C^{q, p + 1}(M, \widehat{\operatorname{Sym}}^{\geq l} Q^* \otimes E) \oplus C^{q + 1, p}(M, \widehat{\operatorname{Sym}}^{\geq (l-1)} Q^* \otimes E) \oplus C^{q + 2, p - 1}(M, \widehat{\operatorname{Sym}}^{\geq l} Q^* \otimes E).
\end{equation*}
From this, we deduce that $d_P^2 s$ is precisely the $C^{q, p+2}(M, \operatorname{Sym}^l Q^* \otimes E)$-component of $D^2s$. Since $D$ is flat, $D^2s = 0$, hence $d_P^2 s = 0$.\par
Under the natural identification
\begin{equation*}
	\Omega^*(M, \widehat{\operatorname{Sym}} Q^* \otimes E) \cong C^{0, *}(M, \textstyle\bigwedge Q^* \otimes \widehat{\operatorname{Sym}} Q^* \otimes E),
\end{equation*}
the operator $d_P$ can be viewed as a flat $P$-connection on $\bigwedge Q^* \otimes \widehat{\operatorname{Sym}} Q^* \otimes E$. It satisfies
\begin{equation*}
	d_P(\alpha \otimes s) = ( d_P^{\bigwedge Q^*} \alpha ) \otimes s + (-1)^q \alpha \wedge d_P(s),
\end{equation*}
for all $\alpha \in \Gamma(M, \bigwedge^q Q^*)$ and $s \in \Gamma(M, \widehat{\operatorname{Sym}} Q^* \otimes E)$, where $d_P^{\bigwedge Q^*}$ denotes the flat $P$-connection on $\bigwedge Q^*$ induced by the Bott $P$-connection on $Q$. Finally, denote by $d_P^{\widehat{\operatorname{Sym}} Q^*}$ the flat $P$-connection on $\widehat{\operatorname{Sym}} Q^*$ induced by the Bott $P$-connection on $Q$, and by $d_P^E$ the flat $P$-connection on $E$ obtained from restricting $d_P$. We can then state the following proposition.

\begin{proposition}
	\label{Proposition: quasi isomorphism}
	Let $D$ be a flat connection on $\widehat{\operatorname{Sym}} Q^* \otimes E$ such that for all $l \in \mathbb{N}$, $D + \delta_Q$ preserves $\widehat{\operatorname{Sym}}^{\geq l} Q^* \otimes E$. Suppose that the map $d_P$ defined as in \eqref{Equation: truncation of flat connection} satisfies
	\begin{equation*}
		d_P(a \otimes s) = (d_P^{\widehat{\operatorname{Sym}} Q^*} a) \otimes s + a \otimes d_P^E s,
	\end{equation*}
	for all $a \in \Gamma(M, \widehat{\operatorname{Sym}} Q^*)$ and $s \in \Gamma(M, E)$. Then the map
	\begin{equation*}
		\pi_Q: (\Omega^*(M, \widehat{\operatorname{Sym}} Q^* \otimes E), D) \to (C^{0, *}(M, E), d_P^E)
	\end{equation*}
	is a quasi-isomorphism. In particular, it induces an isomorphism between the sheaf of $D$-flat sections of $\widehat{\operatorname{Sym}} Q^* \otimes E$ and the sheaf of $d_P^E$-closed sections of $E$.
\end{proposition}

In what follows, we always assume the hypothesis in the above proposition.

\begin{lemma}
	\label{Lemma D.16}
	The following equality holds on $\Omega^*(M, \widehat{\operatorname{Sym}} Q^* \otimes E)$:
	\begin{equation*}
		d_P \delta_Q^{-1} + \delta_Q^{-1} d_P = 0.
	\end{equation*}
\end{lemma}
\begin{proof}
	Since $\delta_Q^{-1}$ annihilates $C^{0, *}(M)$, the graded Leibniz rule implies that for any $\beta \in C^{0, p}(M)$, $\alpha \in \Gamma(M, \bigwedge Q^*)$ and $s \in \Gamma(M, \widehat{\operatorname{Sym}} Q^* \otimes E)$, 
	\begin{align*}
		d_P ( \delta_Q^{-1} (\beta \wedge \alpha \otimes s))
		= & (-1)^p ( d_P \beta ) \wedge \delta_Q^{-1} (\alpha \otimes s) + \beta \wedge  d_P \delta_Q^{-1} (\alpha \otimes s),\\
		\delta_Q^{-1} (d_P (\beta \wedge \alpha \otimes s))
		= & -(-1)^p ( d_P \beta ) \wedge \delta_Q^{-1} (\alpha \otimes s) + \beta \wedge \delta_Q^{-1} d_P (\alpha \otimes s).
	\end{align*}
	The terms involving $d_P\beta$ cancel, so it remains to prove that $(d_P \circ \delta_Q^{-1} + \delta_Q^{-1} \circ d_P) (\alpha \otimes s) = 0$. Choose local frames $(\check{v}_j)$ of $P$ and $(v_i)$ of $Q$. Write $\partial_{\check{v}_j} = \iota_{\check{v}_j} d_P$, $\partial_{\check{v}_j}^Q = \iota_{\check{v}_j} d_P^Q$, etc. Then
	\begin{align*}
		d_P ( \delta_Q^{-1} (\alpha \otimes s)) = & \check{v}^j \wedge \partial_{\check{v}_j} ( u^i \cdot \iota_{v_i} (\alpha \otimes s)) = \check{v}^j \wedge \partial_{\check{v}_j} \left( \left( \iota_{v_i} \alpha \right) \otimes u^i \cdot s \right),\\
		\delta_Q^{-1} (d_P (\alpha \otimes s)) = & u^i \cdot \iota_{v_i} \left( \check{v}^j \wedge \partial_{\check{v}_j} \left( \alpha \otimes s \right) \right) = \check{v}^j \wedge \left( - u^i \cdot \iota_{v_i} \partial_{\check{v}_j} \left( \alpha \otimes s \right) \right).
	\end{align*}
	By the Leibniz rule for each fixed $\check{j}$, the sum of the index-$\check{j}$ components is
	\begin{align*}
		& \left( \partial_{\check{v}_j}^{\bigwedge Q^*} \left( \iota_{v_i} \alpha \right) \right) \otimes u^i \cdot s + \left( \iota_{v_i} \alpha \right) \otimes \left( \partial_{\check{v}_j}^{\widehat{\operatorname{Sym}} Q^*} u^i \right) \cdot s - \left( \iota_{v_i} \partial_{\check{v}_j}^{\bigwedge Q^*} \alpha \right) \otimes u^i \cdot s\\
		= & \iota_{\partial_{\check{v}_j}^Q v_i} \alpha \otimes u^i \cdot s + \left( \iota_{v_i} \alpha \right) \otimes \left( \partial_{\check{v}_j}^{Q^*} u^i \right) \cdot s.
	\end{align*}
	The natural duality relation $\partial_{\check{v}_j}^Q v_i \otimes u^i = -v_i \otimes \partial_{\check{v}_j}^{Q^*} u^i$ shows that this sum vanishes.
\end{proof}

Now, we define a decreasing filtration on $\Omega^*(M, \widehat{\operatorname{Sym}} Q^* \otimes E)$ as follows: for $r \in \mathbb{N}$,
\begin{equation}
	\label{Equation: decreasing filtration}
	\operatorname{F}_r = \operatorname{F}_r \left( \Omega^*(M, \widehat{\operatorname{Sym}} Q^* \otimes E) \right) := \bigoplus_{q+l \geq r} C^{q, *}(M, \operatorname{Sym}^l Q^* \otimes E).
\end{equation}
The $r$th graded piece of this filtration is given by
\begin{equation*}
	\operatorname{Gr}_r = \operatorname{Gr}_r \left( \Omega^*(M, \widehat{\operatorname{Sym}} Q^* \otimes E) \right) := \bigoplus_{q+l = r} C^{q, *}(M, \operatorname{Sym}^l Q^* \otimes E).
\end{equation*}

\begin{lemma}
	\label{Lemma D.15}
	The following equalities hold on $\Omega^*(M, \widehat{\operatorname{Sym}} Q^* \otimes E)$:
	\begin{equation*}
		\pi_Q ( D - (d_P - \delta) ) = (d_P - \delta_Q)^2 = 0.
	\end{equation*}
\end{lemma}
\begin{proof}
	Let $p, q, l \in \mathbb{N}$ and $s \in C^{q, p}(M, \operatorname{Sym}^l Q^* \otimes E)$. We have $(d_P - \delta_Q) s \in \operatorname{Gr}_{q+l}$. By our assumptions, $P$ is involutive and $D + \delta_Q$ preserves $\widehat{\operatorname{Sym}}^{\geq (q+l)} Q^* \otimes E$. Hence,
	\begin{equation}
		D s - (d_P - \delta_Q) s \in \operatorname{F}_{q+l+1}.
	\end{equation}
	It implies that $\pi_Q (D - (d_P - \delta))= 0$. On the other hand,
	\begin{equation*}
		(d_P - \delta_Q)^2s \in \operatorname{Gr}_{q+l} \quad \text{and} \quad D^2s - (d_P - \delta_Q)^2s \in \operatorname{F}_{q+l+1}.
	\end{equation*}
	Since $D^2 s = 0$, we obtain $(d_P - \delta_Q)^2s = 0$.
\end{proof}

\begin{lemma}
	\label{Proposition 3.7}
	The map
	\begin{equation*}
		\Psi: (\Omega^*(M, \widehat{\operatorname{Sym}} Q^* \otimes E), D) \to (\Omega^*(M, \widehat{\operatorname{Sym}} Q^* \otimes E), d_P - \delta_Q)
	\end{equation*}
	defined by $\Psi(s) = s - \delta_Q^{-1} (D - d_P + \delta_Q) s$, where , is a cochain isomorphism.
\end{lemma}
\begin{proof}
	Set $\underline{D} := D - d_P + \delta_Q$. We first show that $\Psi$ is a cochain map. Compute:
	\begin{align*}
		\Psi D - (d_P - \delta_Q) \Psi = & D - \delta_Q^{-1} \underline{D} D - (d_P - \delta_Q) + (d_P - \delta_Q) \delta_Q^{-1} \underline{D}\\
		= & \underline{D} + \delta_Q^{-1} (d_P - \delta_Q) D + (d_P - \delta_Q) \delta_Q^{-1} \underline{D}\\
		= & (1 - \delta_Q^{-1} \delta_Q - \delta_Q \delta_Q^{-1}) \underline{D} + (d_P \delta_Q^{-1} + \delta_Q^{-1} d_P) \underline{D}\\
		= & \pi_Q \underline{D} = 0.
	\end{align*}
	The second and third equalities used $D^2 = 0$ and $(d_P - \delta_Q)^2 = 0$ respectively. The fourth equality used the identity \eqref{Equation: Hodge decomposition for Dolbeault differential} and Lemma \ref{Lemma D.16}. Note that $\delta_Q^{-1} \underline{D}$ sends $\operatorname{F}_r$ to $\operatorname{F}_{r+1}$ for every $r \in \mathbb{N}$. The same argument as in Lemma \ref{Lemma: qausi-isomorphism for formal functions} shows that $\Psi$ is an isomorphism.
\end{proof}

\begin{proof}[\myproof{Proposition}{\ref{Proposition: quasi isomorphism}}]
	It follows from Lemma \ref{Lemma D.15} that $d_P$ is a flat connection. Now, suppose that $d_P \delta_Q^{-1} + \delta_Q^{-1} d_P = 0$. It is clear that $d_P \circ \pi_Q = \pi_Q \circ d_P$ and $\pi_Q \circ \delta_Q = 0$. Thus,
	\begin{equation}
		\label{Equation: quasi isomorphism for polarized sections}
		\pi_Q: (\Omega^*(M, \widehat{\operatorname{Sym}} Q^* \otimes E), d_P - \delta_Q) \to (C^{0, *}(M, E), d_P)
	\end{equation}
	is a cochain map. By our hypothesis together with the identity \eqref{Equation: Hodge decomposition for Dolbeault differential},
	\begin{equation*}
		\operatorname{Id} - \pi_Q = (\delta_Q - d_P) \circ \delta_Q^{-1} + \delta_Q^{-1} \circ (\delta_Q - d_P)
	\end{equation*}
	holds on $\Omega^*(M, \widehat{\operatorname{Sym}} Q^* \otimes E)$. Also note that $\delta_Q = 0$ on $C^{0, *}(M, E)$. We then see that \eqref{Equation: quasi isomorphism for polarized sections} is a quasi-isomorphism. Finally, by the fact that $\pi_Q \circ \delta_Q^{-1} = 0$,
	\begin{equation*}
		\pi_Q = \pi_Q \circ \Psi: (\Omega^*(M, \widehat{\operatorname{Sym}} Q^* \otimes E), D) \to (C^{0, *}(M, E), d_P^E)
	\end{equation*}
	is a quasi-isomorphism, where $\Psi$ is the cochain isomorphism in Lemma \ref{Proposition 3.7}.
\end{proof}

Now, we can apply Proposition \ref{Proposition: quasi isomorphism} to prove Proposition \ref{Proposition: qausi-isomorphism for polarized sections} and Lemma \ref{Lemma: quasi-isomorphism for quantizable functions}.

\begin{proof}[\myproof{Proposition}{\ref{Proposition: qausi-isomorphism for polarized sections}}]
	It suffices to have the following observations, which are direct consequences of Proposition \ref{Lemma 3.7}. First, the operator $\nabla^{\operatorname{K}, E} + \delta_Q$ preserves the subbundle $\widehat{\operatorname{Sym}}^{\geq l} Q^* \otimes E$ for all $l \in \mathbb{N}$. Second, for the operator $d_P$ in \eqref{Equation: truncation of flat connection}, constructed from $\nabla^{\operatorname{K}, E}$, its restriction to each $C^{q, p}(M, \widehat{\operatorname{Sym}} Q^* \otimes E)$ can be written as the composition
	\begin{center}
		\begin{tikzcd}
			C^{q, p}(M, \widehat{\operatorname{Sym}} Q^* \otimes E) \ar[rr, "\nabla^{\widehat{\operatorname{Sym}} Q^* \otimes E}"] && \Omega^{p+q+1}(M, \widehat{\operatorname{Sym}} Q^* \otimes E) \ar[r] & C^{q, p+1}(M, \widehat{\operatorname{Sym}} Q^* \otimes E),
		\end{tikzcd}
	\end{center}
	where the second arrow denotes the canonical projection. Consequently, $d_P$ coincides with the flat $P$-connection on $\bigwedge Q^* \otimes \widehat{\operatorname{Sym}} Q^* \otimes E$ induced by the Bott $P$-connection on $Q$ and the flat $P$-connection $d_P^E$ on $E$. Hence, Proposition \ref{Proposition: qausi-isomorphism for polarized sections} follows directly from Proposition \ref{Proposition: quasi isomorphism}.
\end{proof}

\begin{proof}[\myproof{Lemma}{\ref{Lemma: quasi-isomorphism for quantizable functions}}]
	Fix $r \in \mathbb{N}$. From the explicit formula
	\begin{align*}
		D_k = \nabla - \delta_Q + \frac{k}{\sqrt{-1}} \left[ \widetilde{I} + \delta_Q^{-1}\omega + \frac{\sqrt{-1}}{k} F^\mathbf{L}, \, \cdot\, \right]_{\star_k^{\operatorname{F}}},
	\end{align*}
	we extract two structural properties:
	\begin{itemize}
		\item $D_k$ preserves $\Omega^*(M, \widehat{\operatorname{Sym}} Q^* \otimes \operatorname{Sym}^{\leq r} P^*)$; and
		\item for every $l \in \mathbb{N}$, $D_k + \delta_Q$ preserves the subspace $\Omega^*(M, \operatorname{Sym}^{\geq l} Q^* \otimes \operatorname{Sym}^{\leq r} P^*)$.
	\end{itemize}	
	Construct the operator $d_{P, k}$ from $D_k$ as in \eqref{Equation: truncation of flat connection}. Let $d_P^{\mathcal{W}}$ be the $P$-connection on $\mathcal{W}$ obtained from restricting $\nabla$. Then we can verify that for all $a \in \Gamma(M, \widehat{\operatorname{Sym}} Q^*)$ and $b \in \Gamma(M, \operatorname{Sym}^{\leq r} P^*)$,
	\begin{align*}
		d_{P, k} (a \otimes b) = & d_P^{\mathcal{W}} (a \otimes b) - \frac{k}{\sqrt{-1}} a \otimes \pi_Q\left( b \star_k^{\operatorname{F}} \left( \widetilde{I} + \delta_Q^{-1}\omega + \frac{\sqrt{-1}}{k} F^\mathbf{L} \right) \right)\\
		= & (d_P^{\mathcal{W}} a) \otimes b + a \otimes d_{P, k} b.
	\end{align*}
	In particular, $d_{P, k} = \pi_Q \circ D_k$.\par
	The connection $\nabla$ restricts to a torsion-free connection $\nabla^Q$ on $Q$, which extends the Bott $P$-connection $d_P^Q$ on $Q$. Hence, $d_P^{\mathcal{W}}$ restricts to the flat $P$-connection on $\widehat{\operatorname{Sym}} Q^*$ induced by $d_P^Q$. Thus the result follows from Proposition \ref{Proposition: quasi isomorphism}.
\end{proof}

\bibliographystyle{amsplain}
\bibliography{References}

\end{document}